
\vsize=23.5truecm
\hsize=6.25truein
\nopagenumbers

\parindent=0pt
\parskip=5pt
\baselineskip=13pt
\hsize=6.25truein
\mathsurround=1.5pt
\nopagenumbers
\headline={\ifnum\pageno=1{}\else\hfill{\rm\folio.}\fi}

\newcount\resultcount
\resultcount=0

\newcount\Resultcount
\Resultcount=0

\newcount\sectioncount
\sectioncount=-1

\newcount\eqcount
\eqcount=0

\newcount\myeqcount
\myeqcount=0

\newcount\refcount
\refcount=0
\def\r{\advance\refcount by 1{\number\refcount. }}

\def\result#1#2#3\par{\advance\resultcount by 1%
\setbox0=\hbox{\bf #2 
\number\sectioncount.\number\resultcount.~~~}%
\hangindent=\wd0 \hangafter=-2 \hskip-\parindent\hskip-\wd0\box0
\hskip-.8pt#3%
\advance\refcount by 1\par}

\def\Result#1#2#3\par{\advance\Resultcount by 1%
\setbox0=\hbox{\bf #2 
\number\Resultcount.~~~}%
\hangindent=\wd0 \hangafter=-2 \hskip-\parindent\hskip-\wd0\box0
\hskip-.8pt#3%
\advance\refcount by 1}

\def\pf{\setbox0=\hbox{\bf Proof:~~~}%
\hangindent=\wd0 \hangafter=-2 \hskip-\parindent\hskip-\wd0\box0\rm \hskip-.8pt}

\def\lem#1 #2\par{\result#1{Lemma}{\sl#2}\par\rm}
\def\prop#1 #2\par{\result#1{Proposition}{\sl#2}\par\rm}
\def\thm#1 #2\par{\result#1{Theorem}{\sl#2}\par\rm}
\def\cor#1 #2\par{\result#1{Corollary}{\sl#2}\par\rm}
\def\ex#1 #2\par{\result#1{Example}#2\par}
\def\exer#1 #2\par{\result#1{Exercise}#2\par}
\def\defn#1 #2\par{\result#1{Definition}{\sl#2}\par\rm}

\def\Thm#1 #2\par{\Result#1{Theorem}{\sl#2}\par\rm}

\def\rem#1 #2\par{%
\setbox0=\hbox{\bf Remark.~~~}%
\hangindent=\wd0 \hangafter=-2 \hskip-\parindent\hskip-\wd0\box0
\hskip-.8pt#2%
\advance\refcount by 1\par}
\def\cite#1{[#1]}

\def\target#1{}

\def\paper#1#2#3#4#5#6#7{%
\item{[#1]~~}{\rm#2, }{\sl#3\/}, {\rm#4}~{\bf#5} {\rm(#6), #7}.\target{#1}}
\def\book#1#2#3#4#5#6{%
\item{[#1]~~}{\rm#2, } {\it#3.\/} {\rm#4}: {\rm#5}, {\rm#6}.\target{#1}}
\def\chap#1#2#3#4#5#6#7#8{%
\item{[#1]~~}{\rm#2, }{\sl#3\/}. In: {\rm#4}, %
{\rm#5} {\rm#6} {#7}, {\rm#8}.\target{#1}}

\def\rems#1#2\par{%
\setbox0=\hbox{\bf Remarks.~~~}%
\hangindent=\wd0 \hangafter=-2 \hskip-\parindent\hskip-\wd0\box0
\hskip-.8pt#2%
\advance\refcount by 1\par}

\def\proclaim#1{\setbox0=\hbox{\bf #1~~~}%
\hangindent=\wd0 \hangafter=-2 \hskip-\parindent\hskip-\wd0\box0\hskip-.8pt}

\def\eqtag#1{}

\font\bm=msbm10
\font\bmm=msbm7
 at 9truept

\def\m#1{\hbox{\bm #1}}
\def\mm#1{\hbox{\bmm #1}}

\font\gothfont=eufm10
\def\frak#1{\hbox{\gothfont#1}}


\def\ul#1{$\underline{\hbox{#1}}$}
\def\w{\wedge}

\def\fract#1#2{{\textstyle#1\strut\over\textstyle#2\strut}}
\def\sqr#1#2{{\vbox{\hrule height.#2pt \hbox{\vrule width.#2pt
  height#1pt \kern#1pt \vrule width.#2pt}\hrule height.#2pt}}}
\def\qed{\qquad \sqr 7 3\rm}
\def\:{\;\colon\;}
\def\<{\langle}
\def\>{\rangle}
\def\-{\backslash}
\def\Int{\displaystyle\int}

\def\strutt#1#2{\vrule height #1pt depth #2pt width 0pt}
\def\bmatrix#1{\left[\matrix{#1}\right]}
\def\subsett{\subseteq}

\def\limm{\lim\limits}
\def\d{\partial}
\def\db{\kern1pt\d\kern-5pt\raise 8.25pt%
\hbox{\vrule height .4pt depth .05pt width 5pt\kern.25pt}}
\def\dbs{\db^*\kern-2pt}
\def\dbd{\db\d}
\def\ddb{\d\db}
\def\tr{{\rm tr}\,}

\def\im{{\rm im}\,}

\def\norm#1{\|#1\|}
\def\nnorm#1{\big\|#1\big\|}

\def\tto{\raise1.5pt\hbox{${\scriptstyle\to}\kern-6pt{\scriptstyle\to}$}}

\def\dV{\,dV}
\def\dv{\dV}
\def\lap{\triangle}

\def\Re{{\rm Re}\,}

\def\rk{{\rm rk}\,}

\def\wh{\widehat}

\def\bar{\overline}

\def\et{{\rm E}_{\rm h}}
\def\eet{{\rm End}\,\et}

\def\upin{\cup\hskip-3.55pt plus0pt minus0pt 
\vrule height 5.5pt depth0pt width .5pt }

\def\eqn{\global\advance\eqcount by 1%
(\nums.\number\eqcount)}

\def\myeq{\eqno{\eqn}}
\def\myeqn{\myeq}


\def\nums{\number\sectioncount}

\def\secn#1{\eqcount=0\resultcount=0\advance\sectioncount by 1%
\bigskip{\bf \S\number\sectioncount. \quad #1\medskip}}

\def\lasteq#1{\myeqcount=\eqcount\advance\myeqcount by -#1(\nums.\number\myeqcount)}


\def\source#1#2{#1}

\def \one {\source {Theorem~1}{1}}
\def \two {\source {Theorem~2}{2}}
\def \three {\source {Theorem~3}{3}}
\def \four {\source {Theorem~4}{4}}
\def \eAA {\source {(1.1)}{5}}
\def \eAAb {\source {(1.1)}{6}}
\def \eAt {\source {(1.2)}{7}}
\def \rAA {\source {Proposition~2.1}{8}}
\def \rAB {\source {Corollary~2.2}{9}}
\def \eAB {\source {(2.1)}{10}}
\def \rAC {\source {Lemma~2.3}{11}}
\def \rAD {\source {Corollary~2.4}{12}}
\def \rAE {\source {Lemma~2.5}{13}}
\def \rAG {\source {Proposition~2.6}{14}}
\def \kBx {\source {remarks}{15}}
\def \ePsi {\source {(2.4)}{16}}
\def \rCA {\source {Lemma~3.1}{17}}
\def \rCB {\source {Proposition~3.2}{18}}
\def \kCa {\source {remark}{19}}
\def \rCC {\source {Lemma~4.1}{20}}
\def \rCCa {\source {Corollary~4.2}{21}}
\def \rDequivariant {\source {Corollary~4.3}{22}}
\def \rCE {\source {Proposition~4.4}{23}}
\def \rCD {\source {Proposition~4.5}{24}}
\def \rDhomeo {\source {Corollary~4.6}{25}}
\def \rCK {\source {Theorem~4.7}{26}}
\def \kDa {\source {remark}{27}}
\def \rCL {\source {Corollary~4.8}{28}}
\def \rFA {\source {Lemma~5.1}{29}}
\def \rEeq {\source {Lemma~5.2}{30}}
\def \rEmoment {\source {Corollary~5.3}{31}}
\def \rDB {\source {Theorem~5.4}{32}}
\def \kEa {\source {remark}{33}}
\def \eFtwo {\source {(6.2)}{34}}
\def \eSixthree {\source {(6.3)}{35}}
\def \eFaz {\source {(6.6)}{36}}
\def \eEaa {\source {(6.9)}{37}}
\def \rFagn {\source {Lemma~6.1}{38}}
\def \rEA {\source {Proposition~6.2}{39}}
\def \eEaz {\source {(6.11)}{40}}
\def \eEk {\source {(6.13)}{41}}
\def \eEC {\source {(6.14)}{42}}
\def \rFzz {\source {Lemma~6.3}{43}}
\def \eFvar {\source {(6.15)}{44}}
\def \eEg {\source {(6.16)}{45}}
\def \eEp {\source {(6.17)}{46}}
\def \eEcc {\source {(6.18)}{47}}
\def \rEC {\source {Lemma~6.4}{48}}
\def \eEdd {\source {(6.19)}{49}}
\def \eFfour {\source {(6.21)}{50}}
\def \rED {\source {Proposition~7.1}{51}}
\def \eEx {\source {(7.1)}{52}}
\def \eDeriv {\source {(7.2)}{53}}
\def \rClose {\source {Theorem~7.2}{54}}
\def \rEint {\source {Corollary~7.3}{55}}
\def \kGw {\source {remarks}{56}}
\def \rEuniq {\source {Proposition~7.4}{57}}
\def \kGbb {\source {remark}{58}}
\def \rThere {\source {Theorem~8.1}{59}}
\def \rFgauge {\source {Lemma~8.2}{60}}
\def \eHa {\source {(8.1)}{61}}
\def \eHb {\source {(8.2)}{62}}
\def \eHd {\source {(8.3)}{63}}
\def \rHc {\source {Lemma~8.3}{64}}
\def \eHz {\source {(8.4)}{65}}
\def \rHskew {\source {Lemma~8.4}{66}}
\def \rHatlast {\source {Lemma~8.5}{67}}
\def \kHz {\source {remark}{68}}
\def \Approx {\source {Lemma~9.1}{69}}
\def \nonpoly {\source {Theorem~9.2}{70}}
\def \zero {\source {Proposition~9.3}{71}}
\def \conforms {\source {(9.1)}{72}}


\centerline{\bf POLYSTABILITY AND THE HITCHIN-KOBAYASHI CORRESPONDENCE}

\bigskip

\centerline{Nicholas Buchdahl and Georg Schumacher}

\bigskip
\parshape=1 45pt 5.1in
Using a quasi-linear version of Hodge theory,
holomorphic vector bundles in a neighbourhood
of a given polystable bundle $E_0$ on a compact K\"ahler manifold
are shown to be (poly)stable if and only if their corresponding
classes are (poly)stable in the sense of geometric invariant
theory with respect to the action of the
automorphism group of $E_0$
on the (finite-dimensional vector) space of infinitesimal 
deformations.



\bigskip
\secn{Introduction.}

In the first edition of his book \cite{MFK}, Mumford introduced
the notion of stability in algebraic geometry, giving a 
condition under which the set of orbits of an algebraic
group action can be endowed with a ``good" structure, 
a moduli space.
A numerical characterisation of stability derived from the
Hilbert criterion leads to the notion of stability
for holomorphic vector bundles, with which
Narasimhan and Seshadri \cite{NS} proved that stable
holomorphic vector bundles on a compact Riemann surface
correspond precisely to irreducible projective unitary
representations of the fundamental group. Such representations
in turn correspond to irreducible projectively flat
unitary connections on a given hermitian vector bundle,
and in this form Donaldson \cite{Do1} gave an analytical
(as opposed to a representation-theoretic) proof of the
theorem of Narasimhan and Seshadri.

Donaldson \cite{Do3} subsequently generalised his results to the
case of algebraic surfaces equipped with Hodge metrics
by proving that
a holomorphic vector bundle is stable if and only if
it admits an irreducible Hermite-Einstein connection
(see \S1 for definitions). In the concluding
section of his paper, following the seminal work of
Atiyah and Bott \cite{AB}, he explained how his result could be seen
in a broader context, namely as an infinite-dimensional version
of a theorem of Kempf and Ness \cite{KN} giving a metric
characterisation of stability and relating
categorical quotients of Mumford's geometric invariant
theory with the symplectic quotients of Marsden and Weinstein \cite{MW}.
In addition to their remarkable consequences for the
differential topology of smooth four-manifolds, Donaldson's results
have inspired a wide range of generalisations. In particular,
Uhlenbeck and Yau \cite{UY} proved that a holomorphic
vector bundle on a compact K\"ahler manifold is stable
if and only if it admits an irreducible Hermite-Einstein
connection, this (and other generalisations) known as the
Hitchin-Kobayashi correspondence.

In this paper, we demonstrate that the relationship between
the two notions of stability---that in geometric invariant theory
and that in K\"ahler geometry---is even more intimate than it
might appear at first sight. 
For a holomorphic vector bundle $E$ on a compact complex manifold $X$ 
the space of
infinitesimal deformations is the (finite-dimensional) cohomology 
group $H^1(X,End\,E)$, and
Kuranishi theory describes the parameter space of a semi-universal
deformation of $E$ as the zero set of a holomorphic function $\Psi$ 
defined in a neighbourhood ${\cal N}$ of zero in
$H^1(X,End\,E)$ with values in $H^2(X,End\,E)$. 
If there are non-trivial automorphisms of $E$, different points
in $\Psi^{-1}(0)$ may correspond to isomorphic bundles.
The group $Aut(E)$ of holomorphic automorphisms of $E$ acts on
the space of infinitesimal deformations,
but in general it does not leave $\Psi^{-1}(0)$ invariant.
Even if it does, it cannot be expected that points of the
(set-theoretic/topological) quotient $\Psi^{-1}(0)/Aut(E)$ should
represent isomorphism classes of holomorphic bundles near $E$ unless
$E$ is simple, in which case the group action is trivial.

For a simple holomorphic 
vector bundle $E$ on a compact complex manifold $X$, 
universal deformations exist and the base spaces of these can be 
glued together in a unique way to give a (generally non-Hausdorff) coarse 
moduli space. Using the methods of \S6 of \cite{AHS} 
and the 
Hitchin-Kobayashi correspondence, 
coarse moduli spaces for stable holomorphic vector bundles 
in the category of reduced complex analytic spaces were constructed 
in \cite{Kim} 
and in \cite{FS} 
by a construction of 
local slices for the action of the complex gauge group,
here viewing holomorphic structures as integrable semi-connections,
with isomorphic such structures being (by definition) those in the same
orbit of that group.

If $X$ is equipped with a K\"ahler metric and if $E_0$
is a polystable vector bundle, several remarkable features of the
local description of holomorphic bundles near $E_0$ emerge,
these features otherwise failing to exist:

\medskip

\Thm{\one} Let $E_0$ be a polystable holomorphic vector bundle
equipped with an Hermite-Einstein connection $d_0$ 
on a compact K\"ahler manifold $(X,\omega)$. Then the holomorphic
function $\Psi\: {\cal N}\subsett 
H^1(X,End\,E_0)\to H^2(X,End\,E_0)$ is 
equivariant with respect to the
action of $Aut\,E_0$ on $\db_0$-harmonic representatives. \rm 

\medskip
\Thm{\two} With the same objects as in Theorem 1 and restricting
${\cal N}$ if necessary, every holomorphic 
bundle $E$ corresponding to a point of $\Psi^{-1}(0)$ is
semi-stable, and any destabilising subsheaf of $E$ is a subbundle. \rm

\smallskip

Restricting ${\cal N}$ further if 
necessary, the set of isomorphism classes of holomorphic vector bundles
near to $E_0$ is in 1--1 correspondence with the set-theoretic quotient
$\Psi^{-1}(0)/Aut(E_0)$; indeed, when each of these spaces is 
equipped with the quotient topology, the natural map 
from the latter into the space of isomorphism classes of 
holomorphic bundles near $E_0$ is a homeomorphism onto its 
image. This topological space is non-Hausdorff in general, 
but when attention is restricted to those bundles 
that are 
polystable---corresponding to an $Aut(E_0)$-invariant subset
of $\Psi^{-1}(0)$---the
quotient is Hausdorff. The group $Aut(E_0)$ is a 
complex reductive Lie group, acting linearly on the finite-dimensional
vector spaces $H^q(X,End\,E_0)$ and extending the isometric
actions of the unitary automorphisms when restricted to 
the $\db_0$-harmonic spaces equipped with the $L^2$ 
norm. 

At this point, the description 
of a neighbourhood of $E_0$ has, through the Hitchin-Kobayashi
correspondence and its GIT analogue the Kempf-Ness theorem,
an interpretation in terms of geometric invariant theory:

\medskip
\Thm{\three} A class $\alpha\in \Psi^{-1}(0)\subsett H^1(X,End\,E_0)$
is (poly)stable with respect to the action of $Aut(E_0)$
if and only if the corresponding bundle $E_{\alpha}$
is (poly)stable with respect to $\omega$. \rm

\medskip
For those bundles that are semistable but not polystable, the
following analogue of the standard GIT result holds:

\medskip
\Thm{\four} If $\bar\alpha$ is a polystable point in the closure
of the orbit of $\alpha\in \Psi^{-1}(0)$ under the
action of
$Aut(E_0)$, the polystable
bundle $E_{\bar\alpha}$ is isomorphic to the graded
object $Gr(E)$ associated to a Seshadri filtration of
$E$. \rm

\medskip

In addition to these theorems, 
a number of other results of a quantitative but 
technical nature appear
throughout the paper, many of these not restricted
to integrable connections alone and being of
considerable independent interest---indeed, the role played
by integrability, apart from that of the connection
defining $E_0$, turns out to be relatively minor
throughout.

\smallskip
The results of this  paper shed 
further light on the close
relationship between stability
in the context of geometric invariant theory and
stability for holomorphic
vector bundles on compact K\"ahler manifolds, providing a
framework for the further analysis of moduli spaces of
stable bundles on such manifolds. Although some of these
results may appear
unsurprising at first sight, the fact that polystability
lies at the borderline between stable and unstable renders
inapplicable the standard techniques---notably, the implicit function
theorem---so that new approaches lying
at the interface between infinite and finite-dimensional
analysis are required.

The application of methods from geometric invariant theory 
to problems in geometric analysis is not new, and indeed
some of the methods used here have appeared in the literature
in the context of $K$-(poly)stability and the existence
of constant scalar curvature K\"ahler metrics; see for example,
\cite{S}, \cite{CS},  and also \cite{Th}. 
The (yet to be fully established) 
Yau-Tian-Donaldson conjecture 
is the analogue of the Hitchin-Kobayashi 
correspondence, postulating a one-to-one correspondence
between objects satisfying an algebro-geometric 
condition and the zeros of an appropriate moment map,
manifested as the solution of a particular differential
equation. The results here are of a different nature, describing
an actual correspondence between classical GIT and 
the theory of holomorphic vector bundles, as opposed to 
describing an analogue; the results are not
subsumed by the (established) validity of the Hitchin-Kobayashi
correspondence. Apart from providing potentially useful insight
into an analogous theory, it is not the intention here to contribute to the
developing theory of cscK metrics.

\medskip

The construction of slices for the action of 
the complex gauge group is given in \rAE\ below, this in 
the general setting of a review of the Kuranishi theory of 
deformations of holomorphic vector
bundles. Following a brief study of a neighbourhood of 
a semi-stable holomorphic bundle on a compact K\"ahler manifold in \S3,
the discussion is specialised in \S4 to the case of 
a polystable bundle,
with the bundle represented by an Hermite-Einstein connection.
The equivariance of the function $\Psi$ described earlier is 
proved in \rDequivariant, and the claim that the natural
map mentioned above is a homeomorphism follows from the first statement of
\rCL.  In \S5 there is a general discussion
of some of the basic ideas associated with geometric invariant theory
that will be well-known to experts, and the ``if" part of the
main theorem (\three\ above) 
is established in \rDB. The more difficult ``only if" part of the
main theorem is proved in \S\S6-8, culminating in \rThere.
Provided that polystable bundles 
considered are ``sufficiently close" to $E_0$, 
a quantitative estimate of the distance between the 
corresponding 
Hermite-Einstein connections is provided by 
\rClose; see also the second of the \kGw\ following that result. 
The case of bundles close to $E_0$ that are not polystable
is considered in \S9, this corresponding \four\ above. 
The paper concludes with some general 
remarks concerning the contents and possible extensions of the
paper, and it commences in the next section, which is
devoted to establishing terminology
and notation and providing references.

\bigskip

\secn{Preliminaries.}

Let $X$ be a compact complex manifold.
In this paper, as in \cite{Do3}, holomorphic vector bundles
on $X$ are viewed
most often from the
perspective of integrable $\db$-operators on fixed smooth
bundles.
Following Kobayashi \cite{Ko3},
a {\it semi-connection\/} on a complex
vector bundle $\rm E$ is
a $\m C$-linear map $\db$ on differentiable local sections of $\et$
taking values in $\Lambda^{0,1}\otimes \et$ and
satisfying the $\db$-Leibniz rule; here $\Lambda^{p,q}$ is the
space of $(p,q)$-forms on $X$. Any such operator has natural
prolongations $\Lambda^{0,q}\otimes{\rm E}\to \Lambda^{0,q+1}\otimes
\rm E$,
and the operator
is integrable if $\db \circ \db=0$.
As for connections, the set of
semi-connections on $\rm E$ is an affine space with all the familiar
properties of connections, including the ways in which they can
be induced naturally on associated bundles. The semi-connections
induced by integrable connections are again integrable.

Let $\db_0$ be an integrable semi-connection on $\rm E$.
Then every semi-connection
$\db$ on the bundle can be written $\db=\db_0+a''$ for some
unique $(0,1)$-form $a''$ with coefficients in
$\rm End\,E$, and the integrability condition for $\db$ is
$\db_0a'' + a''\w a''=0$. The notation, which will be used
throughout, is inspired by the often-used convention
to denote by $a'$ and $a''$ respectively the $(1,0)$-
and $(0,1)$-components of a $1$-form $a$, where that
$1$-form may take values in some vector bundle.

The group ${\cal G}$ of complex automorphisms of $\rm E$ acts on
the affine space of semi-connections as a ``complex gauge group".
This action, which preserves the integrability condition,
is denoted  by $g\cdot\db := g\circ \db \circ g^{-1}$.
A holomorphic structure is
defined by an integrable semi-connection, and two such structures
are isomorphic if and only if they lie in the same orbit of
${\cal G}$. By virtue of the Newlander-Nirenberg theorem,
this view of holomorphic structures is equivalent to the more
usual one of holomorphic vector bundles being described by
systems of holomorphic transition functions.

Denote by $A^{p,q}(\rm E)$ the global smooth $(p,q)$-forms
with coefficients in $\rm E$. For an integrable
semi-connection $\db_0$ defining a holomorphic structure
$E_0$, the Dolbeault cohomology groups
$$
H^q_{\overline\partial_0}(X,\rm E)
={\ker\db_0\: A^{0,q}(\rm E)\to A^{0,q+1}({\rm E})\over
\im\db_0\: A^{0,q-1}({\rm E})\to A^{0,q}({\rm E})}
$$
are denoted by $H^q(X,E_0)$, these being finite-dimensional
spaces with $H^1(X,End\,E_0)$ being by definition the
space of infinitesimal deformations of $E_0$.

Analysis of the small deformations of $E_0$ is achieved 
with the introduction of metrics on both the bundle
$\rm E$ and the manifold $X$.
A hermitian metric ${\rm h}$
is fixed once and for all on the bundle $\rm E$, which
is henceforth denoted $\et$. The
group ${\cal G}$ is the complexification of the group
${\cal U}$ of unitary gauge transformations. The hermitian
structure on $\et$ gives a one-to-one correspondence
between semi-connections $\db$ and hermitian connections $d=\d+\db$
on $\et$, and if $d=d_0+a$ for some skew-adjoint
$a\in A^1(\eet)$, then $\d=\d_0+a'$ and $\db=\db_0+a''$ where
$a=a'+a''$ and $a'=-(a'')^*$. Henceforth, all connections are
taken to be hermitian.

The action of ${\cal G}$ on the space of semi-connections
extends to an action on the space of connections via
$g\cdot d := g^*{}^{-1}\circ \d \circ g^* +
g \circ \db \circ g^{-1} = d+g^*{}^{-1}\d g^*
-\db g\,g^{-1}$. The curvature $F(d)=d\circ d
\in A^2(\eet)$ of a connection is a skew-adjoint
$2$-form with coefficients in $\eet$, and the
connection (i.e., the
associated semi-connection) is
integrable if and only if $F(d)$ is of type $(1,1)$.
Since $F(d_0+a)=F(d_0)+d_0a+a\w a$, for $g\in {\cal G}$
it follows that
$$
\eqalignno{
F(g\cdot d) = g\,F^{0,2}(d)\,g^{-1} + g^*{}^{-1}\,F^{2,0}(d)\,g^*
+ \big( F^{1,1}&(d) + \db(g^*{}^{-1}\d g^*)-\d(\db g\,g^{-1})+
& \eqn \cr
& (g^*{}^{-1}\d g^*)\w (\db g\,g^{-1}) +
(\db g\,g^{-1})\w(g^*{}^{-1}\d g^*)\big)\;.
}
$$\eqtag{\eAA}

Now fix a positive $(1,1)$-form $\omega$ on $X$. If $\dim X=n$,
the associated volume form is $dV = \omega^n$, where the
convention is adopted throughout that
$\omega^q := (1/q!)\,\omega\w\cdots \w \omega$ ($q$ times).
If $d$ is an integrable connection on $\et$, standard Hodge
theory on compact manifolds gives a
unique $\db$-harmonic representative in each Dolbeault
cohomology class, where $\tau\in A^{0,q}(\et)$
is $\db$-harmonic if $\db\tau=0=\db^*\tau$ for
$\db^*=-*\d*$, the formal adjoint of $\db$. So
$\tau$ is $\db$-harmonic if and only if it lies in the
kernel of the $\db$-Laplacian $\lap''=\db\db^*
+\db^*\db$. In general,
there is an $L^2$-orthogonal decomposition
$$
A^{0,q}(\et)=(\ker\db)^{\perp}\oplus (\ker\db^*)^{\perp}\oplus H^{0,q}
= \overline{\im \db^*} \oplus \overline{\im\db}\oplus H^{0,q}
$$\eqtag{\eAAb}
where $H^{0,q}=H^{0,q}(\db)$ is the space of
$\db$-harmonic $(0,q)$-forms. Here, notation has been
abused in that $A^{0,q}(\et)$ is no longer denoting the
space of smooth sections of $\et$, but rather the space
of global sections of $\Lambda^{0,1}\otimes \et$ that
are square integrable, and the closures on the right
are the closures in $L^2$ of the images under $\db$
and $\db^*$ of the spaces of smooth global sections.
Standard elliptic regularity implies that the $\db$-harmonic
sections are smooth, at least if the connections are.

This abuse of notation will be employed throughout, so that
$A^{0,q}(\et)$ will always denote a space of global sections
of $\Lambda^{0,q}\otimes \et$ but with the degree of
differentiability and/or integrability to be specified in
the respective context. Sobolev spaces of functions
are denoted by $L^p_k$, meaning all weak derivatives up to
and including those of order $k$ lie in $L^p$.
Having fixed a base connection on $\et$ once and
for all, the spaces of $L^p_k$ elements of $A^{0,q}(\et)$
acquire norms that make them Banach spaces.

Henceforth, a number $p>2n$ (for $n=\dim X$) will
be fixed, so by the Sobolev
embedding theorem there are compact embeddings
$L^p_1\subset C^0$ and $L^p_2\subset C^1$.
By standard elliptic theory on compact manifolds,
for an integrable connection $d=\d+\db$ there
is a constant $C>0$ (depending
upon $d$) such that
$$
\norm{\tau}_{L^p_1}\le C\big(\norm{\db\tau}_{L^p}+
\norm{\db^*\tau}_{L^p}+\norm{\Pi^{0,q}\tau}_{L^2}\big)\quad{\rm for}\;
\tau\in A^{0,q}(\et)\;,\myeq
$$\eqtag{\eAt}
where $\Pi^{0,q}\tau$ is the $L^2$-orthogonal projection of
$\tau\in A^{0,q}(\et)$ in $H^{0,q}(\db)$.
Connections on $\et$ will be permitted to have
coefficients in $L^p_1$, and the fact that an
integrable $L^p_1$ connection defines a holomorphic structure
in the usual way follows from Lemma~8 of \cite{Bu2}.
When not indicated by a subscript on the norm symbol,
$\norm\tau$ will always mean $\norm{\tau}_{L^2}$.

If $d$ is an integrable semi-connection and $g\in {\cal G}$,
the Dolbeault cohomology groups defined by $\db$ and by $g\cdot \db$
are isomorphic, the isomorphism induced by mapping
a $\db$-closed $(0,q)$-form $\tau\in A^{0,q}(\et)$ to
the $(g\cdot\db)$-closed $(0,q)$-form $g\,\tau$. This isomorphism
does not preserve harmonic representatives in general, unless
$g\in {\cal U}$ in which case it also preserves $L^2$ norms.

Subsequently in this paper it will be useful to consider connections
that are not integrable, in which case the Dolbeault cohomology
groups are not defined. However, one can still define the
spaces $H^{0,q}(\db)$ of $\db$-harmonic $(0,q)$-forms
as null spaces of the appropriate Laplacians,
these still being finite-dimensional spaces consisting only
of smooth forms (if $\db$ is itself smooth).

\smallskip

If $d\omega=0$, the formal adjoints $\d^*$ and
$\db^*$ have alternative expressions in terms of the
K\"ahler identities:
$$
\d^*=i(\Lambda\db -\db\Lambda)\,,\quad
\db^*=-i(\Lambda\d-\d\Lambda)\;,
$$
where $\Lambda\: \Lambda^{p+1,q+1}\to \Lambda^{p,q}$ is
the adjoint of $\omega$, so $\Lambda\omega=n$. For
an integrable connection $d$ on $\et$ with curvature
$F=F(d)$, the Bianchi and K\"ahler identities imply
that the Yang-Mills equations $d^*F=0$ are equivalent to
the equation $d\wh F=0$, where $\wh F := \Lambda F$, the
central component of the curvature.
Under these circumstances, the bundle and connection
split into eigenspaces of the covariantly constant
self-adjoint endomorphism $i\wh F$, and when restricted
to any such eigenspace, the curvature of the restricted
connection has central component that is a constant multiple
of the identity. That is, it is an {\it Hermite-Einstein\/}
connection, a term coined by Kobayashi in \cite{Ko1}.

If $d_0$ is an Hermite-Einstein connection on
$\et$ with $i\wh F(d_0)= \lambda\,1$, the
constant $\lambda$ is a fixed positive multiple of the
{\it slope\/} $\mu(\et)=deg(\et)/\rk(\et)$, where
$deg(\et)=(n-1)!\,[\omega^{n-1}]\cdot c_1(\et)$.
Kobayashi \cite{Ko2} proved that a holomorphic
bundle $E$ defined by an Hermite-Einstein connection
is semi-stable in the sense of Mumford-Takemoto
\cite{Ta},
meaning that any coherent subsheaf
$A\subset E$
with non-zero torsion-free quotient should
satisfy $\mu(A)\le \mu(B)$, and if the connection
is irreducible (meaning that the unitary bundle-with-connection
does not split), then the bundle is stable, meaning
$\mu(A) <\mu(B)$. These results were proved independently
by L\"ubke \cite{L\"u}. Kobayashi had earlier mooted in
\cite{Ko1} that there should exist a relationship between
the algebro-geometric notion of stability and the
existence of Hermite-Einstein connections, this
also conjectured by Hitchin \cite{Hi} in the broader
setting of arbitrary compact manifolds equipped with
Gauduchon metrics \cite{G}. Following the work of Donaldson,
the affirmation of that conjecture was established
by Uhlenbeck and Yau \cite{UY} in the K\"ahler case
and by Li and Yau \cite{LY} in the general case. Many
other authors have made significant contributions to the
theory, extending the correspondence in a variety of ways.
Of relevance here is the semi-stable version
of the correspondence due to Jacob \cite{J} and
a version for torsion-free semi-stable sheaves due to
Bando and Siu \cite{BS}.

\bigskip

\secn{A neighbourhood of a holomorphic bundle.}

As in the previous section, let $X$ be a compact complex
manifold
and let $E_0$ be a holomorphic vector bundle on $X$.
A vector bundle version of Kuranishi's fundamental results \cite{Kur}
on deformations of compact complex manifolds has been
given by Forster and Knorr \cite{FK} using power series methods
and by Miyajima \cite{Miy} using more analytical methods; (cf.\ also \cite{FS}). 
Either way, there is a holomorphic function $\Psi$ defined
in a neighbourhood of $0\in H^1(X,End\,E_0)$ such that
$\Psi^{-1}(0)$ is a complete family of small deformations of $E_0$.
In this section, the construction of the function $\Psi$
will be presented in a manner to suit the purposes of the remainder
of the paper. The entire discussion is essentially an
$n$-dimensional version of the $2$-dimensional case
presented in \S6.4 of \cite{DK}.

\smallskip
Fix a positive $(1,1)$-form $\omega$ on $X$, which at
this stage is not assumed to be $d$-closed. Let $\et$ be the
complex bundle underlying $E_0$ equipped with a
fixed hermitian structure, and let $d_0$ be a connection on
$\et$ inducing the holomorphic structure $E_0$.

\smallskip
\prop{\rAA}
There is a number $\epsilon>0$ depending on $d_0$ with the
property that for
any integrable hermitian connection
$d_a = d_0+a$ with $\norm{a}_{L^p_1}
<\epsilon$ the $L^2$ orthogonal projection
$
H^{0,q}(\db_a) \owns \tau \mapsto \Pi^{0,q}\tau\in H^{0,q}(\db_0)
$
is injective.

\pf Write $a=a'+a''$ for $a''\in A^{0,1}(\eet)$.
Suppose $\tau\in A^{0,q}(\et)$ satisfies
$\db_0\tau+a''\w \tau=0=\db_0^*\tau-i\Lambda(a'\w \tau)$, so
$\tau\in H^{0,q}(\db_a)$. If $\Pi^{0,q}\tau=0
\in H^{0,q}(d_0)$, then $\tau=\db_0\mu+\db_0^*\nu$
for some $\mu\in A^{0,q-1}(\et)$
and some $\nu\in A^{0,q+1}(\et)$, with $\mu$ and $\nu$
respectively orthogonal to the kernels of $\db_0$
and $\db_0^*$. Then $-a''\w\tau=\db_0\tau
= \db_0^{}\db_0^*\nu$ and $i\,\Lambda(a'\w\tau)
=\db_0^*\tau=\db_0^*\db_0^{}\mu$, so on taking inner
products with $\nu$ and $\mu$ respectively it follows that
$\norm{\db_0^*\nu}^2\le \norm{\nu}\,\norm{a''\w \tau}$
and $\norm{\db_0\mu}^2\le \norm{\mu}\,\norm{a'\w\tau}$,
where $\norm{\cdot}$ is the $L^2$ norm. Since $p>2n$,
the Sobolev embedding theorem gives
$\sup|a| \le C\,\norm{a}_{L^p_1}$ for some constant $C$
independent of $a$ and $d_0$, so after adding the last
two inequalities it follows that $\norm{\tau}^2
=\norm{\db_0\mu}^2+\norm{\db_0^*\nu}^2
\le  C\,\norm{a}_{L^p_1}\,\norm{\tau}\big(\norm{\mu}+\norm{\nu}\big)$.
Since $\mu$ and $\nu$ are orthogonal to $\ker\db_0$ and
$\db_0^*$ respectively, ellipticity of $\lap_0''$ implies
that $\norm{\mu}\le C_1\,\norm{\db_0\mu}$ and
$\norm{\nu}\le C_2\,\norm{\db_0^*\nu}$ for some constants
$C_1,\,C_2$, so
$\norm{\tau}^2\le C_3\norm{a}_{L^p_1}\,\norm{\tau}^2$ for
some new constant $C_3=C_3(d_0)$, giving the stated result.
\quad\qed

\medskip
This proposition clearly implies the well-known and
standard semi-continuity
of cohomology. It also has the following useful consequence,
resulting from equality of dimensions of cohomology groups:

\smallskip
\cor{\rAB} Under the hypotheses of {\rm\rAA}, suppose in addition
that $d_0$ and $d_a$ lie in the same ${\cal G}$-orbit. Then
the map  $H^{0,q}(\db_a)\to H^{0,q}(\db_0)$ induced by orthogonal
projection is an isomorphism. \quad \qed
\rm

\medskip
Regardless of the integrability or otherwise of
$\db_a=\db_0+a''$, if $\norm{a''}_{L^p_1}$ is sufficiently small, the
Sobolev embedding theorem combined with \eAt\ gives
the following perturbed version of that estimate:

\smallskip

\lem{\rAC} There exists $\epsilon>0$ and $C>0$ depending on $d_0$
with the properties that if $a''\in A^{0,1}(\eet)$ satisfies
$\norm{a''}_{L^p_1} < \epsilon$, then
$$
\norm{\tau}_{L^p_1} \le C\big(\norm{\db_0^*\tau}_{L^p}
+ \norm{\db_0\tau+a''\w\tau}_{L^p} +
 \norm{\Pi^{0,q}\tau}_{L^2}\big)\;, \qquad \tau\in A^{0,q}(\et)\;,
\myeq
$$\eqtag{\eAB}
where $\Pi^{0,q}\:A^{0,q}(\eet)\to H^{0,q}(\db_0)$ is $L^2$
orthogonal projection.
\quad \qed

\medskip

Replacing $\et$ by $\eet$, $d_0$ by the induced connection
$d_0$ on $\eet$, taking $q=1$ and $\tau=a''$ gives
\smallskip

\cor{\rAD} There exists $\epsilon>0$ and $C>0$ with the properties that
each $a''\in A^{0,1}(\eet)$ with $\norm{a''}_{L^p_1} < \epsilon$
satisfies
$$
\norm{a''}_{L^p_1} \le C\big(\norm{\db_0^*a''}_{L^p}
+ \norm{\db_0a''+a''\w a''}_{L^p} + \norm{\Pi^{0,1}a''}_{L^2}\big)\;,
\myeq
$$
where $\Pi^{0,1}\:A^{0,1}(\eet)\to H^{0,1}(\db_0)$ is  $L^2$ orthogonal
projection.
\quad \qed

\rm

\medskip

If $a''\in A^{0,1}(\eet)$ and $g\in {\cal G}$, then
$g\cdot(\db_0+a'') = \db_0-\db_0g\,g^{-1}+ga'' g^{-1}$. The
map
$$
{\cal G}\times A^{0,1}(\eet)\owns
(g,a'')\mapsto \db_0^*(-\db_0g\,g^{-1}+ga'' g^{-1})
\in A^{0,0}(\eet)
$$
maps into the subspace of $A^{0,0}(\eet)$ orthogonal to
the kernel of $\db_0$, and its linearization at $(1,0)$
in the ${\cal G}$-direction is $A^{0,0}(\eet)
\owns \gamma\mapsto - \lap_0''\,\gamma$.  This is an isomorphism
from the space of $L^p_2$ sections in $A^{0,0}(\eet)$ orthogonal
to $\ker\db_0=H^{0,0}$ to the space of $L^p$ sections in
$A^{0,0}(\eet)$ orthogonal to $\ker\db_0$. The implicit function
theorem for Banach spaces now implies:

\smallskip

\lem{\rAE} There exist $\epsilon,\, C>0$ with the property that for
each $a'' \in A^{0,1}(\eet)$ with
$\norm{a''}_{L^p_1}<\epsilon$ there is a unique
$\varphi\in A^{0,0}(\eet)\cap (\ker\db_0)^{\perp}$
with $\norm{\varphi}_{L^p_2}
\le C\,\norm{\db_0^*a''}_{L^p}$
such that $\db_0^*(-\db_0g\,g^{-1}+g\,a''\,g^{-1})=0$,
where $g = \exp(\varphi)$.  \quad\qed

\rm
\smallskip

This is a complex analogue of fixing a unitary gauge for
hermitian connections near a
given such connection, corresponding to the linear operation
of projecting
a $(0,1)$-form orthogonal to the range of $\db_0$.
If
$\db_0+a''$ is an integrable semi-connection with
$\norm{a''}_{L^p_1}$ sufficiently small as dictated by
this lemma, after applying an appropriate complex gauge
transformation so that the new semi-connection
$\db_0+\tilde a''$ satisfies $\db_0^*\tilde a''=0$,
\rAD\ gives an estimate of the form
$\norm{\tilde a''}_{L^p_1} \le
C\,\big\|\Pi^{0,1}\tilde a''\big\|_{L^2}$.
Consequently, the well-known result that if  $H^1(X,End\,E_0)=0$
then every small deformation of $E_0$ is
isomorphic to $E_0$ follows immediately.
A simple but pertinent example is given by a
trivial bundle on $\m P_1$.

\medskip
It follows from \rAD\ that for integrable semi-connections
in the ``good" complex gauge lying in a sufficiently small
neighbourhood of $\db_0$ in $L^p_1$, the projection onto the
$\db_0$-harmonic component is a homeomorphism onto a closed subset
of an open neighbourhood of $0$ in $H^{0,1}$, where
$H^{0,q}$ denotes $H^{0,q}(\db_0,\eet)$ in this section. This closed subset
is the zero-set of the holomorphic function $\Psi$ mentioned in the
introduction, as
will now be discussed.

\smallskip
If $a''\in A^{0,1}(\eet)$ is $\db_0$-closed, the semi-connection
$\db_0+a''$ is not integrable in general, since $a''\w a''$
is not zero in general. But one might attempt to perturb
$a''$ in such a way that the corresponding perturbed
semi-connection is integrable.

The derivative of the map
$$
\eqalignno{
\big(A^{0,1}\times A^{0,2}\big)&(\eet) \quad\longrightarrow \quad
 A^{0,2}(\eet)\cr
&\upin \hskip1.4in \upin & \eqn \cr
 (a''&,\beta) \quad \mapsto \quad
\db_0(a''+\db_0^*\beta)+(a''+\db_0^*\beta)\w (a''+\db_0^*\beta)
}
$$
in the $A^{0,2}$-direction at $(0,0)$ is
$
A^{0,2}(\eet)\owns b\mapsto \db_0^{}\db_0^*\,b
\in A^{0,2}(\eet)
$,
which is an isomorphism from the closed subspace of
the $L^p_2$ sections in $A^{0,2}(\eet)$ that are
orthogonal in $L^2$ to the kernel of $\db_0^*$
onto the closed subspace of the $L^p$ sections in  $A^{0,2}(\eet)$
orthogonal in $L^2$ to the kernel of $\db_0^*$.
The non-linear mapping \lasteq0 does not map
into the latter closed subspace in general, but if $\Pi^{0,2}$ is the
$L^2$ orthogonal projection of $A^{0,2}(\eet)$
onto the closed subspace $\ker\db_0^*$ and $\Pi^{0,2}_{\;\perp}=
1-\Pi^{0,2}$,
the map \lasteq0\ can be replaced by
$$
(a'',\beta)\mapsto
\Pi^{0,2}_{\perp}\big[\,\db_0(a''+\db_0^*\beta)+
(a''+\db_0^*\beta)\w (a''+\db_0^*\beta)\,\big]
$$
to obtain a map with the same linearization in the $A^{0,2}$-direction
at $(0,0)$.

Given $\db_0$-closed $a''\in A^{0,1}(\eet)$ in $L^p_1$
and  $\beta\in A^{0,2}(\eet)$ in $L^p_2$,
the form $\tau := \db_0(a''+\db_0^*\beta)+(a''+\db_0^*\beta)
\w (a''+\db_0^*\beta)$ is an $L^p$ section in
$A^{0,2}(\eet)$. The weak derivative $\db_0\tau$
of $\tau$ is
$\db_0^{}\db_0^*\beta\w(a''+\db_0^*\beta) - (a''+\db_0^*\beta)\w
\db_0^{}\db_0^*\beta$, which lies in $L^p$ (again using $p>2n$
and the Sobolev embedding theorem).
This implies that
$\Pi^{0,2}_{\perp}\tau$ lies $L^p$, and
that the composition $(a'',\beta) \mapsto
\Pi^{0,2}_{\perp}\,\tau$ is continuous with respect to the $L^p_2$
topology
on the domain and the $L^p$ topology on the codomain.

By ellipticity of $\lap_0''$ on $A^{0,2}(\eet)$,
there is a constant
$K$ such that $\sup_X|e|\le K$ for any $e\in H^{0,2}$ with
$\norm{e}_{L^2}=1$. Since
$\Pi^{0,2}(\db_0a''+a''\w a'')=\Pi^{0,2}(a''\w a'')$,
it follows that
$\norm{\Pi^{0,2}(\db_0a''+a''\w a'')}_{L^p}
\le Const.\norm{a''}_{L^{4}}^2$,
so another application of the implicit function theorem yields
\medskip

\prop{\rAG} There exist $\epsilon>0$ and $C>0$
with the properties that for
any $a''\in A^{0,1}(\eet)$
satisfying $\norm{a''}_{L^p_1}
<\epsilon$ there is a unique
$\beta\in (\ker\db_0^*)^{\perp}\subsett A^{0,2}(\eet)$
satisfying $\norm{\beta}_{L^p_2}\le C\,\norm{a''}_{L^4}^2$
such that $\db_0(a''+\db_0^*\beta)+(a''+\db_0^*\beta)
\w (a''+\db_0^*\beta)$ is $\db_0^*$-closed.
\quad \qed

\bigskip
\rems{\kBx}

\proclaim{1. }
If $\beta$ is as in this lemma, then
$\norm{a''+\db_0^*\beta}_{L^p_1}$ is uniformly
bounded by a constant multiple of $\norm{a''}_{L^p_1}$. In
fact, if $\tilde a'' := a''+\db_0^*\beta$, then
$\tau := \db_0\tilde a''+\tilde a''\w \tilde a''$
satisfies
$\db_0\tau + \tilde a''\w \tau-\tau\w \tilde a''=0=\db_0^*\tau$ weakly,
so by elliptic regularity it follows that $\tau$ in fact lies in
$L^p_1$. Hence by \rAC, if
$\norm{a''}_{L^p_1}$ (and hence $\norm{\tilde a''}_{L^p_1}$)
is sufficiently small,
there is a constant $C = C(d_0)$ such that $\norm{\tau}_{L^p_1}
\le C\,\norm{\Pi^{0,2}\tau}_{L^2}$.

\proclaim{2. }
\rAG\ remains valid even if $d_0$ is not integrable---all that is
required is a uniform $C^0$ bound on $F^{0,2}(d_0)$.
The proof as given only needs modification by noting that
$\db_0\tau$ involves an extra term
$F^{0,2}(d_0)\w a''-a''\w F^{0,2}(d_0)$, this
lying in $L^p$ if $|F(d_0)|$ is bounded in $C^0$.

\proclaim{3. } If $\gamma$ is a complex automorphism of $\et$
satisfying $d_0\gamma=0$, then by the uniqueness statement of
the lemma, $\beta(\gamma a''\gamma^{-1})=\gamma\beta(a'')\gamma^{-1}$,
at least if $\norm{\gamma a''\gamma^{-1}}_{L^4}$ is sufficiently small.

\smallskip

From \eAt\ there is a constant $c>0$ depending only on $d_0$ such that
any $\alpha\in A^{0,1}$ with $\db_0\alpha=0=\db_0^*\alpha$
satisfies $\norm{\alpha}_{L^p_1}\le c\,\norm{\alpha}$. Thus
there is a number $\delta>0$ depending only on $d_0$
such that $\norm{\alpha}_{L^p_1} < \epsilon$ if
$\norm{\alpha}<\delta$, and for such $\alpha$
there is a unique $\beta\in A^{0,2}(\eet)$ orthogonal to
$\ker\db_0^*$ with $\norm{\beta}_{L^p_2}\le C'\norm{\alpha}^2$
for which the form $\tau = \db_0a''+a''\w a''$ is
$\db_0^*$-closed, where $a'' := \alpha+\db_0^*\beta$.
Moreover, $\norm{\tau}_{L^p_1}\le C\norm{\Pi^{0,2}\tau}$ for
some constant $C=C(d_0)$.

Define the function $\Psi$ on the set of $\db_0$-harmonic
forms $\alpha\in A^{0,1}(\eet)$ with $\norm{\alpha}<\delta$
that takes values in $H^{0,2}(\db_0,\eet)$ by
$$
\Psi(\alpha) := \Pi^{0,2}(a''\w a'') \quad
\hbox{for $a''=\alpha+\db_0^*\beta$ with $\beta=\beta(\alpha)$
as in \rAG.} \myeqn
$$\eqtag{\ePsi}
Then the zero set of $\Psi$ parameterises precisely the 
integrable connections
in this $L^p_1$ neighbourhood of $d_0$, in the sense of defining a 
semi-universal deformation.
By fixing orthonormal bases for each of $H^{0,1}$ and $H^{0,2}$ 
and working in
these bases, it is immediate that $\Psi$ is holomorphic. As will be
shown later (\rDequivariant), in the special case that $d_0$ is 
Hermite-Einstein, 
$\Psi$ is equivariant with respect to the action of
$Aut(E_0)$ on these spaces, a consequence of the third
remark above.

\bigskip

\secn{A neighbourhood of a semi-stable bundle.}

With the same notation and conventions as in the previous section,
we assume from now on that $(X,\omega)$ is K\"ahler.
The following results on bounds of slopes of subsheaves
may also be of independent interest.

\medskip

\lem{\rCA} Let $E$ be a holomorphic bundle defined by an integrable
semi-connection $\db$ on $\et$ and let $d=\db+\d$ be the
associated hermitian connection. Then
there is a constant $C$ depending only on $\omega$
such that
$$
\deg(A) \le C\,\norm{\widehat F(d)}_{L^1}\quad{\rm and}\quad
\norm{c_1(A)} \le C\,\big(\norm{F(d)}_{L^1}+|\deg(A)|\big)
$$
for any coherent analytic sheaf $A\subset E$ with torsion-free
quotient. Here $\norm{c_1(A)}$ denotes the $L^2$ norm of the harmonic
$(1,1)$-form representing the image of $c_1(A)$ in
$H^2_{dR}(X)$.

\pf Suppose first that $A\subset E$ is a holomorphic
sub-\ul{bundle} of $E$ and let $B := E/A$ be the quotient
bundle.
In a unitary frame for $\et = A_{\rm h}\oplus B_{\rm h}$,
the connection $d=: d_E$ has the form
$$
d_E=\bmatrix{d_A& \beta\cr -\beta^*& d_B}
$$
where $d_A$ and $d_B$ are the induced hermitian connections
on $A$ and $B$, and $\beta\in A^{0,1}(Hom(B,A))$ is a
$\db$-closed $(0,1)$-form representing the extension
$0\to A\to E\to B\to 0$. The curvature $F_E=F(d_E)$ has the form
$$
F_E= \bmatrix{F_A-\beta\w\beta^* & d_{BA}\beta\cr -d_{AB}\beta^*
& F_B-\beta^*\!\w \beta}\;,
$$
where $d_{BA}$ here is the
connection on $Hom(B,A)$ induced by $d_A$ and $d_B$.
So if $\Pi_A$ is pointwise-orthogonal projection $E\to A$,
it follows that $F_A = \Pi_A\,F_E\,\Pi_A +\beta\w\beta^*$.
Since $\beta$ is a $(0,1)$-form, $i\,\tr\beta\w \beta^*$ is
a {\it non-positive\/} $(1,1)$-form and therefore
$i\,\tr F_A \le \tr (\Pi_A\,i F_E\,\Pi_A)$. Applying
$\omega^{n-1}\w$ and integrating over $X$, it follows that
$c_1(A)\cdot[\omega^{n-1}]$ is
bounded above by a fixed multiple of $\norm{\widehat F_E}_{L^1}$.

Now if $A$ is only a subsheaf of $E$ of rank $a>0$
with torsion-free quotient $B$,  replace
$E$ with $\Lambda^aE$, $A$ by the maximal normal extension
of $\Lambda^aA$ in $\Lambda^aE$ (which is a line bundle)
and then after blowing up the zero set of the induced
section of $Hom(\det A,\Lambda^aE)$
and resolving singularities,
there is again an upper bound on the degree of the
desingularised subsheaf in terms of the $L^1$ norm of
$F_E$ on the blowup. This
upper bound depends on the metric used
on the blowup, but as in \S\S2, 3 of \cite{Bu2}, there is	
a family $\omega_{\epsilon}$ of such
metrics converging to the pullback of $\omega$, and the
resulting limit then gives the same bound: $\deg(A)$ is
bounded above by a fixed multiple of $\norm{\widehat F(d)}_{L^1}$
for any subsheaf $A\subset E$ with torsion-free quotient.

To obtain  uniform bounds on $\norm{c_1(A)}$,
suppose again initially that $A$
is a subbundle of $E$. For notational simplicity,
let $f := i\,\tr F_A$ and $g := i\,\tr(\Pi_A F_E\Pi_A)$,
so from the preceding arguments, $g-f\ge 0$ as hermitian forms on $TX$.
The space $H^{1,1}_{\;\mm R}(X)$ of real harmonic $(1,1)$-forms
is finite dimensional, so by picking an orthonormal basis,
it is apparent that there is a constant $C>0$ depending
on $\omega$ such that $-C\,\omega\le \varphi\le C\,\omega$
for any $\varphi\in H^{1,1}_{\;\mm R}(X)$ with
$\norm{\varphi}=1$; equivalently, $C\omega\pm\varphi\ge 0$. Therefore
$
(C\omega\pm\varphi)\w(g-f)\ge 0$, implying that
$\pm\varphi\w f\le (C\omega\pm\varphi)\w g -C\omega\w f$ as
real $(2,2)$-forms pointwise on
$X$. Applying $\omega^{n-2}\w$ and integrating over $X$, it
follows that
$$
\pm\int_X\omega^{n-2}\w\varphi\w f
\le \int_X\omega^{n-2}\w (C\omega\pm \varphi)\w g-C\int_X\omega^{n-1}\w f\;,
$$
and by allowing $\varphi$
to vary over the harmonic $(1,1)$-forms of norm $1$, it
follows that $\norm{c_1(A)}_{L^2} \le C'\norm{F(d_E)}_{L^1} \le
C'\big(\norm{F(d_E)}+|\deg(A)|\big)$
for some new constant $C'$ independent of $d_E$.
This implies the second statement of the
lemma when $A$ is a  sub-bundle.

In the case that $A$ is only a subsheaf rather than a subbundle,
the same method as earlier can be used to reduce to
the case of a line subbundle on a blowup of $X$. It need only be
checked that for metrics of the kind $\omega_{\epsilon}$ mentioned earlier,
$\lim_{\epsilon\to 0}\norm{c_1(L)}_{L^2(\omega_{\epsilon})}
= 0$ for any line bundle $L$ on the blowup that is trivial off
the exceptional divisor, which is straight-forward to verify.
\quad \qed

\bigskip
\rCA\ implies the following result, showing that
semi-stability is an open condition.

\medskip

\prop{\rCB} Let $E_0$ be a semi-stable bundle on a compact K\"ahler
manifold $(X,\omega)$, defined by an integrable connection
$d_0$. Then there exists
$\epsilon>0$ such that any integrable connection $d_0+a$
with $\norm{a}_{L^2_1}+\norm{a}_{L^4}<\epsilon$ defines
a semi-stable holomorphic structure.

\pf
If not, there is a sequence $(a_j)\in A^{1}(\eet)$
with $\norm{a_j}_{L^2_1}+\norm{a_j}_{L^4}\to 0$ and $d_j := d_0+a_j$
integrable such that the holomorphic bundle $E_j$ defined by
$d_j$ is not semi-stable, so there is a subsheaf $A_j\subset E_j$
with torsion-free quotient that strictly destabilises
$E_j$. Passing to a subsequence, it can be assumed that the
ranks of the sheaves $A_j$ are constant, $a$ say.
The hypotheses on $a_j$ imply that $\norm{F(d_j)}_{L^2}$ is
uniformly bounded and therefore so too is $\norm{F(d_j)}_{L^1}$,
and consequently \rCA\ yields a uniform upper bound on $\deg(A_j)$.
Since $\deg(A_j)$ is also uniformly bounded below,
these bounds together with the bounds on $\norm{F(d_j)}_{L^2}$
then give uniform bounds on the $L^2$ norms of the
harmonic representatives
of the forms representing $c_1(A_j)$ in $H^2_{dR}(X)$,
and hence there is a convergent subsequence. Since the
image of $H^2(X,\m Z)$ in $H^2_{dR}(X)$ is discrete, this
convergent subsequence must be eventually constant, so
after passing to another subsequence, it can be assumed
that $c_1(A_j)$ is constant, $c$ say. Since $X$ is
K\"ahler, $Pic^c(X)$ is a compact torus, so after passing
to another subsequence, it can be supposed that
$\det A_j$ converges to a holomorphic line bundle
$L$ on $X$. Since $(\det A_j)^*\otimes \Lambda^a E_j$
has a non-zero holomorphic section for each $j$, so too does
$L^*\otimes \Lambda^a E_0$, and therefore $\Lambda^aE_0$ is strictly
destabilised by $L$. But by Theorem~2 of \cite{J}, semi-stability of
$E_0$ implies that of $\Lambda^aE_0$, giving the
desired contradiction. \quad \qed

\bigskip
\rem{\kCa}
It is worth noting that this result is rather delicate
in that the
K\"ahler class must be fixed. For example, every non-split
extension of the form $0\to {\cal O}(-1,1)\to E\to {\cal O}(1,0)
\to 0$ on $\m P_1\times \m P_1$ is strictly stable with
respect to $\omega_t := \pi_1^*\omega_0+t\,\pi_2^*\omega_0$
for $t>2$, is semi-stable but not polystable for $t=2$, and
is strictly unstable if $0<t<2$. The result also fails
in the non-K\"ahlerian case, at least when the degree
fails to be topological. For example, if $L$ is a non-trivial
holomorphic line bundle on an Inoue surface with $L\otimes L$ trivial,
the direct sum of $L$ with the trivial line bundle ${\bf 1}$ is
polystable with respect to every Gauduchon metric, every
small deformation is again a direct sum, and of these, the
generic one is strictly unstable with respect to every
Gauduchon metric. In this case, the
automorphism group $\Gamma$ of $E_0=L\oplus {\bf 1}$
acts trivially on $H^1(X,End\,E_0)$, so each of its
orbits is closed.

\bigskip
\secn{A neighbourhood of a polystable bundle.}

A holomorphic bundle $E$ that is a non-split extension
$0\to A\to E\to B\to 0$ by semistable bundles $A,B$ of
the same slope is semistable but not polystable, and
cannot be separated from the direct sum $A\oplus B$ in
the quotient topology on $\{integrable~semiconnections\}/\{
complex~gauge\}$. Given this, it makes sense to focus on
a neighbourhoods of polystable bundles, in which case
a great deal more can be said than in the previous section.
The same objects and definitions as in
the last section are used here, but now $E_0$ is assumed
to be a polystable holomorphic bundle.

The following is a minor generalisation of a
well-known result essentially due
to Kobayashi \cite{Ko1}.  Although its proof is elementary,
it is presented here for the reason that in some
respects, it is the pivotal result used in the paper.

\medskip

\medskip
\lem{\rCC} Let $d$ be a connection on a bundle $E$  with
$\wh F(d)=0$. If $s\in A^{0,0}(\et)$ satisfies
$\db s=0$, then $d s=0$.

\pf
The equation $\db s=0$ implies
$\d\<s,s\>=\<s,\d s\>$, using
the convention that $\<\cdot,\cdot\>$ is conjugate-linear
in the first variable. Therefore
$\dbd\<s,s\>=\<\d s\,{\buildrel\w\over,}\,\d s\>
+ \<s,F^{1,1}(d)s\>$. Applying $i\Lambda$, it follows
$\lap''|s|^2+|\d s|^2=0$, so integration
over $X$ gives $\norm{\d s}^2=0$. \quad\qed

\medskip
Note that it is not assumed that $d$ should be integrable.

For a connection $d_0$ on $\et$ with central curvature
$\wh F(d_0)$ that is a constant multiple of the identity, the
induced connection on $\eet$ has $\wh F$ identically zero,
so \rCC\ implies
\medskip

\cor{\rCCa} Suppose $d_0$ is a connection on $\et$ with
$i\wh F(d_0)=\lambda\,1$ for some scalar $\lambda$.
If $\sigma\in A^{0,0}(\eet)$ satisfies $\db_0\sigma=0$,
then $d_0\sigma=0$. \quad \qed

\medskip
This corollary yields the equivariance property of the function
$\Psi$ asserted at the end of \S2:

\medskip

\cor{\rDequivariant} Under the hypotheses of \rCCa, there
is a constant $\epsilon=\epsilon(d_0)>0$ with the property that
for any $\db_0$-harmonic $\alpha\in A^{0,1}(\eet)$ and
$\db_0$-closed $\gamma\in {\cal G}$ satisfying 
$\norm{\alpha}_{L^2}+\norm{\gamma\alpha\gamma^{-1}}_{L^2}
<\epsilon$, the form $\beta=\beta(\alpha)\in A^{0,2}(\eet)$
 of \rAG\ satisfies
$\beta(\gamma\alpha\gamma^{-1})=\gamma\beta(\alpha)\gamma^{-1}$.

\pf Given $a''\in A^{0,1}(\eet)$ with 
$\norm{a''}_{L^p_1}$ sufficiently small, \rAG\ guarantees the
existence of a unique $\beta = \beta(a'')\in  A^{0,2}(\eet)$
orthogonal to $\ker\db_0^*$ and with $\norm{\beta}_{L^p_2}
\le C\norm{a''}_{L^4}^2$ such that $\tilde a'' := a''+\db_0^*\beta$
satisfies 
$\db_0^*\big(\db_0\tilde a''+\tilde a''\w \tilde a'')=0$.
If $\gamma\in {\cal G}$ is $\db_0$-closed, then \rCCa\ implies
that $\d_0\gamma=0$ and therefore conjugation by $\gamma$
commutes with both $\db_0$ and $\db_0^*$. Thus if $\beta$ is orthogonal to 
$\ker\db_0^*$ so too is $\gamma\beta\gamma^{-1}$ and also
$\db_0^*\big(\db_0(\gamma\tilde a''\gamma^{-1})
+\gamma(\tilde a''\w \tilde a'')\gamma^{-1})=0$. Since
$\gamma\tilde a''\gamma^{-1}=\gamma a''\gamma^{-1}
+\db_0^*(\gamma\beta\gamma^{-1})$, it follows from the uniqueness
statement of \rAG\ that $\beta(\gamma a''\gamma^{-1})
=\gamma\beta(a'')\gamma^{-1}$ if both $a''$ and 
$\gamma a''\gamma^{-1}$ are sufficiently small in 
$L^p_1$. If $a''=\alpha$ is $\db_0$-harmonic, then so too is
$\gamma\alpha\gamma^{-1}$, and the $L^p_1$ norms of these
forms are uniformly bounded by a fixed multiple of their
$L^2$ norms. \quad \qed

\bigskip

The proof of \rCB\ combined with \rCC\ have the following
useful consequence,  which simplifies a number of
subsequent arguments:

\medskip

\prop{\rCE} Suppose $d_0$ is an integrable connection
with $i\wh F(d_0)=\lambda\,1$ defining a polystable
holomorphic structure $E_0$. There exists $\epsilon>0$
such that any integrable semi-connection $\db=\db_0+a''$
with $\norm{a''}_{L^p_1}<\epsilon$ defines a
semi-stable holomorphic bundle $E$ with the property that
any  subsheaf $A\subset E$ with $\mu(A)=\mu(E)$ and
with $E/A$ torsion-free is a sub-\ul{bundle} of $E$.

\pf If not, then from the proof of \rCB,
there is a sequence of integrable semi-connections
$\db_j=\db_0+a''_j$ with $\norm{a''_j}_{L^p_1}
\to 0$ such that the corresponding holomorphic bundle $E_j$ has a subsheaf
$A_j$ with $\rk A_j = a$ independent of
$j$, $\mu(A_j)=\mu(E_j)=\mu(E_0)$ for all $j$, and with
$E_j/A_j$ torsion-free but not locally free. Hence
$\Lambda^a A_j$ is a rank $1$ subsheaf of $\Lambda^a E_j$,
and the bundle $(\det A_j)^*\otimes \Lambda^aE_j$ has
a holomorphic section that has a zero. After
passing to a subsequence, the Hermite-Einstein connections on the
line bundles $(\det A_j)^*$ can be assumed to converge
to define a holomorphic line bundle $L$ on $X$, and
compactness of the
embedding of  $L^p_1$ in $C^0$ implies that a subsequence
of the integrable connections defining $(\det A_j)^*\otimes \Lambda^aE_j$
converges  uniformly in $C^0$ to the Hermite-Einstein
connection on $L\otimes E_0$.
After scaling the sections of $(\det A_j)^*\otimes \Lambda^aE_j$
so as to have $L^2$ norm $1$, the convergence of the
connections implies that a subsequence of the rescaled
sections converges weakly in $L^2_1$ (say) and strongly in $L^2$ to
a holomorphic section of  $L\otimes \Lambda^a E_0$
of $L^2$ norm $1$. Every section in the sequence has a
zero, and therefore so too does this limit:
the Cauchy integral formula for a holomorphic function
in a polydisk integrated over a poly-annulus  shows
that a sequence of holomorphic functions converging in
$L^1$ on the polydisk converges uniformly on compact
sub-polydisks. But the limiting section is a
holomorphic section of a bundle of degree zero that
admits a connection with $\wh F=0$, so it is either
identically zero or nowhere zero.
\quad\qed

\medskip

The conclusion of \rCC\ can be strengthened
to give a perturbed version; integrability of
$d_0$ is again not required.
\smallskip

\prop{\rCD} Suppose $d_0$ is a connection on
$\et$ with
$\wh F(d_0)=0$. There is a constant
$\epsilon=\epsilon(d_0)>0$ such that any connection
$d_a=d_0+a$ with $a=a'+a''$, $\db_0^*a''=0$ and
$\norm{a}_{L^p_1}<\epsilon$ has the property
that any section $s\in \Gamma(X,\et)$ with
$\db_0s+a''s=0$ must in fact satisfy $\db_0s=0=a''s$.

\pf Suppose $\db_0s+a''s=0$ for some section $s \in A^{0,0}(\et)$,
and write
$s=s_0+s_1$ where $\db_0s_0=0$ and $s_1$
is orthogonal in $L^2$ to the kernel of $\db_0$. By
\rCC, $d_0s_0=0$. Then
$\db_0s_1+a''s_1+a''s_0=0$, so after applying $\db_0^*=
-i\,\Lambda\d_0$ it follows that
$0=\lap_0'' s_1+i\,\Lambda (a''\w \d_0 s_1)=
\lap_0'' s_1+ \db_0^*(a''s_1)$.
Using the $L^2$ inner product, this implies
that $\norm{\db_0s_1}^2=-\<\db_0s_1,a''s_1\>
\le \sup|a''|\,\norm{\db_0s_1}\,\norm{s_1}$,
and therefore $\norm{\db_0s_1} \le \sup|a''|\,\norm{s_1}$.
Since
$p>2n$, there is a constant $C$ such that
$\sup|a''|\le C\norm{a''}_{L^p_1}$, and
since $s_1$ is orthogonal to $\ker\db_0$, there
is a constant $c=c(d_0)>0$ independent of $s$
such that $c\,\norm{s_1}\le
\norm{\db_0s_1}$. So if $C\,\norm{a''}_{L^p_1} < c$, then
$s_1$ must be $0$, giving $s=s_0\in \ker\db_0$, with
$0=\db_0s+a''s=\db_0s_0+a''s_0=a''s_0$. \quad \qed

\medskip
This proposition leads to the verification of
the assertion made in the introduction
that if $d_0$ is an Hermite-Einstein connection on $\et$,
then a neighbourhood of $[d_0]\in \Psi^{-1}(0)/Aut(E_0)$
is homeomorphic to a neighbourhood of $[d_0]$ in the
set of isomorphism classes of integrable connections near 
$d_0$ equipped with the quotient topology:

\medskip
\cor{\rDhomeo} Let $d_0$ be an Hermite-Einstein
connection defining a holomorphic structure $E_0$,
and let $\Psi$ be the function of {\rm \ePsi}, defined in 
a neighbourhood of zero in $H^{0,1}=H^{0,1}(\et,\db_0)$ with values
in $H^{0,2}$. Let ${\cal F}:= \{d_0+a\mid F^{0,2}(d_0+a)=0\}$,
equipped with the $L^p_1$ topology.
If $\Psi$ is restricted to a sufficiently
small neighbourhood of zero, then with respect to
the quotient topologies,  the natural 
map $\Psi^{-1}(0)/Aut(E_0)\to 
{\cal F}/{\cal G}$ is a 
homeomorphism onto a neighbourhood of $[d_0]$.

\pf The continuity of the map is clear. If 
$\alpha_0\in H^{0,1}$, ellipticity 
of the $\db_0$-Laplacian implies that 
$\norm{\alpha_0}_{L^p_1}\le C\norm{\alpha_0}_{L^2}$ for
some constant $C=C(d_0)$. Hence if $\norm{\alpha_0}_{L^2}$
is sufficiently small,
\rAG\ yields a unique $\beta_0\in A^{0,2}(\eet)$
orthogonal to $\ker\db_0^*$
with $\norm{\beta_0}_{L^p_2}\le C\norm{\alpha_0}_{L^2}^2$
such that $a_0'' := \alpha+\db_0^*\beta_0$ satisfies
$\db_0^*(\db_0a_0''+a_0''\w a_0'')=0$. If $a=a'+a''$ 
is such that $\norm{a_0-a}_{L^p_1}$ is so small that
\rAE\ applies to $d_0+a$, then there is a unique $g\in {\cal G}$ 
of the form $g=\exp(\varphi)$ with $\varphi \in(\ker\db_0)^{\perp}$
and $\norm{\varphi}_{L^p_2}\le C\norm{\db_0^*a''}_{L^p}$
such that $d_0+a_1 := g\cdot (d_0+a)$ satisfies 
$\db_0^*a_1''=0$. Since $p>2n$, the Sobolev embedding
theorem gives uniform $C^1$ estimates on $g$ and $g^{-1}$,
implying that $\norm{a''-a_1''}_{L^p_1}\le 
C\norm{\db_0^*a''}_{L^p}=C\norm{\db_0^*(a''-a_0'')}_{L^p}$
for some new constant $C=C(d_0)$. So if $a$ is sufficiently
close to $a_0$ in $L^p_1$, \rAG\ now gives a uniquely 
determined
$\alpha_1\in H^{0,1}$ and $\beta_1
\in (\ker\db_0^*)^{\perp}\subsett A^{0,2}(\eet)$ such that
$a_1''=\alpha_1+\db_0^*\beta_1$, with 
$\db_0^*\big(\db_0a_1''+a_1''\w a_1'')=0$ and with 
$\norm{\beta_1}_{L^p_2} \le C\norm{\alpha_1}^2$. Thus
for some new constant $C=C(d_0)$, 
$\norm{\alpha_1-\alpha_0}_{L^2}
\le \norm{a_1''-a_0''}_{L^2} \le 
\norm{a_1''-a''}_{L^p_1}+\norm{a''-a_0''}_{L^p_1}
\le C\norm{a''-a_0''}_{L^p_1}$. 

In summary, given $\alpha_0
\in H^{0,1}$ sufficiently close to $0$, for each $a\in A^{0,1}(\eet)$
sufficiently close to $a_0$ in $L^p_1$ there exists
$g\in {\cal G}$ and $\alpha_1\in H^{0,1}$
such that $d_0+a$ is of the
form $d_0+a = g^{-1}\cdot (d_0+a_1)$ for
$a_1''=\alpha_1+\db_0^*\beta_1$ with 
$\norm{\alpha_1-\alpha_0}_{L^2}$ uniformly bounded by a multiple
of $\norm{a-a_0}_{L^p_1}$. If ${\cal A}$ denotes
$A^{0,1}(\eet)$ equipped with the $L^p_1$ norm 
and $\Gamma \subsett {\cal G}$ is the group of
$\db_0$-closed (complex) automorphisms of $\et$,
this implies that the natural map 
$[\alpha]\mapsto [\alpha+\db_0^*\beta]$ from a neighbourhood 
of $[0] \in H^{0,1}/\Gamma$ to ${\cal A}/{\cal G}$ is
an open mapping, and therefore so too is the 
natural map $\Psi^{-1}(0)/\Gamma\to {\cal F}/{\cal G}$.

To see that the mapping is 1--1, suppose that   
$d_0+a_0$, $d_0+a_1$ are (integrable) connections
with $\db_0^*a_0''=0=\db_0^*a_1''$ 
that define isomorphic holomorphic structures; that is,
there exists $g\in {\cal G}$ with $d_0+a_1=g\cdot(d_0+a_0)$.
If $\norm{a_0}_{L^p_1}+\norm{a_1}_{L^p_2} < \epsilon$ 
with $\epsilon$
as in \rCD, apply that proposition 
to the connection on 
$\eet=\rm Hom(\et,\et)$ defined by $d_0+a_0$ on one
side and $d_0+a_1$ on the other, with $s$ being the
section $g\in\rm Hom(\et,\et)$. It follows from that
result that $d_0g=0$ (so $g\in \Gamma$) and that $a_1g=ga_0$. Writing 
$a_j''=\alpha_j+\db_0^*\beta_j$ with $\beta_j$ determined
by \rAG, orthogonality of the decomposition implies
that $\alpha_1g=g\alpha_0$; that is, $\alpha_0,\,\alpha_1
\in H^{0,1}$ represent the same point in the quotient 
$H^{0,1}/\Gamma$. \quad \qed

\smallskip
Note that again, the integrability of the connections
$d_0+a_j$ is not critical; the essential ingredient is 
the ``quasi-integrability" condition implicit in the
statement of \rAG. 

\bigskip
The main application of \rCD\
will be to the connections on $\eet$ induced by
connections $d_0$ on $\et$
with $\wh F(d_0)$ a scalar multiple of the identity,
as in following theorem.
Here $[x,y]=xy-yx$ is the usual Lie bracket on endomorphisms,
and the
reader is alerted to the fact that the hypotheses on
$\sigma$ differ slightly in 1.\ and 2., although the
conclusions are essentially the same.

\medskip

\thm{\rCK} Let $d_0$ be a connection on $\et$ with
$i\wh F(d_0)=\lambda\,1$, and let $d_0+a$ be
a connection with $a=a'+a''$, where $a'=-(a'')^*$,
$a''=\alpha+\db_0^*\beta$ for some
$\alpha\in A^{0,1}(\eet)$ satisfying $\db_0\alpha
=0=\db_0^*\alpha$ and $\beta\in A^{0,2}(\eet)$.
There is a constant
$\epsilon>0$ depending only on $d_0$ with the property
that if $\norm{a}_{L^p_1}<\epsilon$ then the following hold
for any  endomorphism $\sigma\in A^{0,0}(\eet)$:
\par
{\parindent=25 pt \advance\rightskip by 30pt \advance\leftskip by 10pt
\item{1. } If $\db_0\sigma+[a'',\sigma]=0$ then $d_0\sigma=0$
and $[\alpha,\sigma]=0=[\db_0^*\beta,\sigma]$. Furthermore,
$[\beta,\sigma]=0$
if $\beta\in (\ker\db_0^*)^{\perp}$.
\item{2. } If $\db_0\sigma+[\alpha,\sigma]=0$ and $\beta\in (\ker\db_0^*)^{\perp}$
and $\db_0^*(\db_0^{}a''+a''\w a'')=0$
then $d_0\sigma=0$ and $[\alpha,\sigma]=0=[\beta,\sigma]$.
\item{3. } If $\db_0\sigma+[a'',\sigma]=0$ and $i\Lambda(\alpha\w\alpha^*
+\alpha^*\w\alpha)$ is orthogonal to $\ker\db_0$ then
$d_0\sigma=0$ and $[\alpha,\sigma]$,
$[\db_0^*\beta,\sigma]$ and $[\alpha,\sigma^*]$ all vanish. Furthermore,
$[\beta,\sigma]=0$
if $\beta\in (\ker\db_0^*)^{\perp}$.
\item{4. } If $\db_0\sigma+[a'',\sigma]=0$ and $\db_0^*(\db_0^{}a''+a''\w a'')=0$
and $i\Lambda(\alpha\w\alpha^*+\alpha^*\w\alpha)\in(\ker\db_0)^{\perp}$ and
$\beta\in(\ker\db_0^*)^{\perp}$ then
$d_0\sigma=0$ and $[\alpha,\sigma]$,
$[\beta,\sigma]$, $[\alpha,\sigma^*]$ and $[\beta,\sigma^*]$ all vanish.
\par}

\medskip

\pf
{\bf 1.}~ The connection induced on $\eet$ by $d_0$ has
$\wh F=0$, so from \rCCa\ and \rCD,
$
d_0\sigma=0=[a'',\sigma]
=[\alpha,\sigma]+[\db_0^*\beta,\sigma]
=[\alpha,\sigma]+\db_0^*[\beta,\sigma]\;.
$
Applying $\db_0$ gives
$\db_0^{}\db_0^*[\beta,\sigma]=0$, from which it follows
that $\db_0^*[\beta,\sigma]=0=[\alpha,\sigma]$.
If $\beta$ is orthogonal to $\ker \db_0^*$, so too
is $[\beta,\sigma]$ by virtue of identity
$\<\psi,[\beta,\sigma]\>=\<[\psi,\sigma^*],\beta\>$,
and if $\psi\in \ker \db_0^*$ then so too is
$[\psi,\sigma^*]$ since $d_0\sigma^*=0$.

\proclaim{2. } If $\tau := \db_0a''+a''\w a''$, then
$\tau = \db_0^{}\db_0^*\beta+\alpha\w\db_0^*\beta
+\db_0^*\beta\w\alpha+\alpha\w\alpha+\db_0^*\beta\w\db_0^*\beta$.
Setting $\beta=0$ in 1.\ it follows from that
case that $d_0\sigma=0=[\alpha,\sigma]$. Since $\sigma$ commutes
with $\alpha$, it also commutes with $\alpha\w\alpha$, implying that
$$
[\sigma,\tau]=
\db_0^{}\db_0^*[\sigma,\beta]+\alpha\w\db_0^*[\sigma,\beta]
+\db_0^*[\sigma,\beta]\w\alpha +
\db_0^*[\sigma,\beta]\w\db_0^*\beta
+\db_0^*\beta\w\db_0^*[\sigma,\beta]\;. \myeqn
$$
As in the proof of 1., because $\beta$ is orthogonal to
$\ker\db_0^*$, so too is $[\sigma,\beta]$. On the other
hand, since $\db_0^*\tau=0$, so too is
$\db_0^*[\sigma,\tau]=0$. Thus $[\sigma,\beta]$ is orthogonal
to $[\sigma,\tau]$.  Taking the inner
product on both sides of \lasteq0with $[\sigma,\beta]$,
rearranging terms and estimating,
it follows that
$$
\big\|\db_0^*[\sigma,\beta]\big\|^2
\le C\big(\sup|\alpha|+\sup|\db_0^*\beta|\big)
\big\|\db_0^*[\sigma,\beta]\big\|\,\big\|[\sigma,\beta]\big\|\;.
$$
Since $[\sigma,\beta]$ is orthogonal to $\ker\db_0^*$,
there is a constant $c>0$ depending only on $d_0$ such
that $\big\|{\db_0^*[\sigma,\beta]\big\|}^2 \ge c\nnorm{[\sigma,\beta]}^2$.
Since $\alpha$ is $\db_0$-harmonic, $\sup|\alpha|$ is bounded
by a constant multiple of $\norm{\alpha}$, and by \rAG\ it can
be assumed that
$\sup|\beta|$ is bounded by a constant multiple of
$\norm{\alpha}^2$. It follows that there is
a new constant $\epsilon_1$ depending only on
$d_0$ for which, if $\norm{\alpha}<\epsilon_1$, then
necessarily $[\sigma,\beta]=0$.

\proclaim{3. } By 1., $d_0\sigma=0$ and $[\alpha,\sigma]=0=[\db_0^*\beta,\sigma]$.
Since $\d_0\sigma=0$ it follows that
$\db_0\sigma^*=0$.  Using the fact that $\sigma$ commutes
with $\alpha$, a short calculation gives
$$
\tr\big((\sigma^*\alpha-\alpha \sigma^*)^*\!
\w(\sigma^*\alpha-\alpha \sigma^*)\big)
=\tr\big((\sigma^*\sigma-\sigma\sigma^*)
(\alpha^*\!\w\alpha+\alpha\w\alpha^*)\big)\;.
$$
After applying $\omega^{n-1}\w$ to both sides and integrating
over $X$, the fact that $(\sigma^*\sigma-\sigma^*\sigma)$
lies in $\ker\db_0$ together
with the fact that $\Lambda(\alpha^*\!\w\alpha+\alpha\w\alpha^*)$
is orthogonal to $\ker\db_0$ imply that
$\big\|\,[\sigma^*\!,\alpha]\,\big\|_{L^2}^2=0$.

\proclaim{4. } The vanishing of $d_0\sigma$, $[\alpha,\sigma]$ and
$[\beta,\sigma]$ follows from 1.
Then from 3.\ it follows that $[\alpha,\sigma^*]=0$, so
applying 2.\ to $\sigma^*$ gives
$[\beta,\sigma^*]=0$.
\quad \qed

\bigskip

\rem{\kDa} Given $a''\in A^{0,1}(\eet)$ with $\norm{a''}_{L^p_1}$
small, \rAE\ yields a uniquely determined $g\in {\cal G}$
near $1$ such that $\db_0^*(g\cdot a'')=0$, where
$g\cdot a'' := ga''g^{-1}-\db_0gg^{-1}$. Also, \rAG\ yields
a uniquely determined $\beta\in A^{0,2}(\eet)$ for which
$\tau := a''+\db_0^*\beta$ satisfies
$\db_0\tau+\tau\w\tau \in \ker\db_0^*$. The two operations
do not commute in general, but if $a_0''\in A^{0,1}(\eet)$
satisfies $\db_0^*a_0''=0$ and
$\db_0^*\big(\db_0^{}a''+a''_0\w a_0''\big)=0$, then for
$\beta\in A^{0,2}(\eet)$ and
$g$ in ${\cal G}$ with both $\beta$ and $g-1$ sufficiently
close to zero in $L^p_2$, the
element $a_1''$ of $A^{0,1}(\eet)$ obtained by applying
the operation of \rAE\ to  $g\cdot a_0''+\db_0^*\beta$
followed by the operation of \rAG\ to the result of that
has the form $a_1'' = \gamma a_0''\gamma^{-1}$ for some
automorphism $\gamma\in \Gamma$ near $1$, where $\Gamma$
here and subsequently denotes the $\db_0$-closed 
elements of ${\cal G}$.
All of this follows using the analysis of \rCK.

\bigskip

The content of
\rCK\ implies and is implied by a corresponding result for
connections on bundles of degree zero. For the purposes of
transparency of the proof, the theorem was presented in terms
of endomorphisms, but the alternative result is the following:

\medskip
\cor{\rCL} Let $d_0$ be a connection on $\et$ with
$\wh F(d_0)=0$, and let $d_0+a$ be
a connection with $a=a'+a''$,  $\db_0^*a''=0$,
$a''=\alpha+\db_0^*\beta$ for some
$\alpha\in A^{0,1}(\eet)$ satisfying $\db_0\alpha
=0=\db_0^*\alpha$ and some $\beta\in A^{0,2}(\eet)$.
There is a constant
$\epsilon>0$ depending only on $d_0$ with the property
that if $\norm{a}_{L^p_1}<\epsilon$ then the following hold
for any  section $s$ in $A^{0,0}(\et)$:
\par
{\parindent=25 pt \advance\rightskip by 30pt \advance\leftskip by 10pt
\i If $\db_0s+a''s=0$ then $d_0s=0$
and $\alpha s=0=(\db_0^*\beta) s$. Furthermore $\beta\,s=0$ if
$\beta\in (\ker\db_0^*)^{\perp}$.
\i If $\db_0s+\alpha s=0$ and $\beta\in(\ker\db_0^*)^{\perp}$ and
$\db_0^*(\db_0^{}a''+a''\w a'')=0$ then
$d_0s=0$ and $\alpha\,s=0=\beta\,s$;
\i If $\db_0s+a''s=0$ and $i\Lambda(\alpha\w\alpha^*
+\alpha^*\w\alpha) \in (\ker\db_0)^{\perp}$ then
$d_0s=0$ and $\alpha s$,
$(\db_0^*\beta) s$ and $\alpha^*s$ are all zero. Furthermore,
$\beta\,s=0$ if $\beta\in (\ker\db_0^*)^{\perp}$.
\i If $\db_0s+a''s=0$ and $\db_0^*(\db_0^{}a''+a''\w a'')=0$
and $i\Lambda(\alpha\w\alpha^*+\alpha^*\w\alpha)
\in (\ker\db_0)^{\perp}$ and $\beta\in(\ker\db_0^*)^{\perp}$ then
$d_0s=0$ and $\alpha s$,
$\beta s$, $\alpha^*s$ and $\beta^*s$ are all zero.
Hence $d_0s=0=d_as$.
\par}

\medskip
\pf Let ${\bf 1}$ be the trivial line bundle equipped with
its standard flat connection, and apply \rCK\ to the
endomorphism $\sigma = \bmatrix{0&s\cr 0&0}$ of
$\et\oplus{\bf 1}$ using the direct sum connections on
this bundle. \quad \qed

\bigskip

When the connection $d_0+a$ of \rCK\ is integrable,
the second statement of the theorem implies that if
the holomorphic structure $E$ defined by $d_0+a$ is
not simple then the isotropy subgroup of $\Gamma$
at $\alpha$ has dimension greater than $1$; (since
$\Gamma$ is acting by conjugation, the constant
multiples of the identity are always in the isotropy
subgroup). The converse of this follows from the third statement.
As will be discussed in the next section, the condition that
$i\Lambda(\alpha^*\!\w\alpha+\alpha\w\alpha^*)$
should be orthogonal to $\ker\db_0$ implies that
the form $\alpha$ is an element of $H^{0,1}$ that is polystable with
respect to the action of $\Gamma$ in the sense of geometric
invariant theory, and in that context, the last statement
of the theorem
can also be interpreted as a holomorphic condition on $E$
determined by an algebraic condition on $\alpha$.
This is a manifestation of one half of the
``local" Hitchin-Kobayashi correspondence that is at the heart
of this paper, made precise in \rDB\ below. 
The other half, or ``converse" of this is the focus
of attention in \S6, \S7 and \S8.

\bigskip

\secn{$\Gamma$-stability and $\omega$-stability.}

This section commences with a summary of notions and
ideas from geometric invariant theory, in particular of those
related to stability, formulated to suit the purposes of this paper.
The primary references for the definitions used here are \cite{MFK},
\cite{Nes}, \cite{Kir} and \cite{Th} in which the definitions
are essentially consistent. It is a testament to the breadth of
the applicability of the ideas that the literature is not always
consistent in its definitions.

\smallskip

Recall that when a reductive Lie group $G$ acts linearly
on a finite-dimensional
complex vector space $V$, a point $v\in V\-\{0\}$ is {\it unstable\/}
for the action if $0\in \overline{G\cdot v}$, is {\it semi-stable\/} if
$0\not\in \overline{G\cdot v}$, is {\it polystable\/}
if $v$ is semi-stable
and $G\cdot v$ is closed, and is {\it stable\/} if $v$ is
polystable and
the isotropy subgroup of $v$ is finite.
Fixing a positive hermitian form on $V$, in the closure of
each orbit
there is a unique point of smallest norm.
For each $v\in V$, the derivative at $1\in G$
of  the function $G\owns g\mapsto \norm{g\cdot v}^2$
gives a function
$\mu\: V \to  \frak{g}^*$, the moment map, and its zeros
on the set of semistable points are precisely the points
of smallest norm in the closed orbits. If $G$ is the complexification
of a real compact group $K$ acting isometrically on $V$ and
if $W$ is a $G$-invariant subset contained in the open set
of stable points,
then $W/G$ is naturally identified with $(\mu^{-1}(0)\cap W)/K$, this
a consequence of the fundamental results of Kempf and Ness \cite{KN}.

The Hilbert-Mumford criterion states that a point $v\in V$ is
stable if and only if it is stable for every $1$-parameter subgroup
in $G$. A $1$-parameter subgroup is given by a homomorphism
$\chi\: \m C^*\to G$, giving a representation of $\m C^*$ on $V$.
The irreducible representations of $\m C^*$ are all $1$-dimensional,
so $V$ splits as a direct sum of $1$-dimensional subspaces $V_j$ on each
of which $\m C^*$ acts with a given weight $w_j\in \m Z$:
$\m C^*\times V_j
\owns (t,v) \mapsto \chi(t)\cdot v = t^{w_j}\,v$. If $v$ has
$1$-dimensional
orbit under $\chi$, it is clear that this orbit is closed if and only
if $\m C^*\owns t\mapsto \chi(t)\cdot v\in V$ is proper, which in turn
is equivalent to the condition that the maximum weight $w^{\sharp}$
and minimum weight $w^{\flat}$
of $\chi$ on the non-zero components of $v$ in its decomposition into
irreducibles should differ in sign. Since $t\mapsto \chi(t^{-1})$
is another $1$-parameter subgroup of $G$ for which the maximum
and minimum such weights are respectively $-w^{\flat}$
and $-w^{\sharp}$, the criterion reduces
to the condition that $v$ is a stable point for the
action of $G$ if and only if $w^{\flat}$ is negative for every
$1$-parameter subgroup of $G$. In practice, this is
the condition that $\limm_{t\to0}\big(\log\norm{\chi(t)\cdot v}/\log|t|\big)
<0$ for every $1$-parameter subgroup $\chi$ of $G$, which is
an analogue of the numerical condition in the definition of
stability for a holomorphic vector bundle on a compact K\"ahler
manifold.

\medskip
Consider now the situation discussed in the previous section:
$(X,\omega)$ is a compact K\"ahler $n$-manifold and  $\et$ is a
complex vector bundle over $X$ equipped with a fixed hermitian
structure. The group ${\cal G}$ of complex automorphisms
of $\et$ acts on the affine space ${\cal A}$ of hermitian
connections on $\et$, preserving the subspace of integrable
such connections. The group ${\cal G}$ is the complexification
of the group ${\cal U}$ of unitary gauge transformations.

A connection $d_0\in {\cal A}$ that is integrable and has
curvature $F(d_0)$ satisfying $i\wh F(d_0)=\lambda\,1$
is a minimum of the Yang-Mills functional. By \rCC, the group
$\Gamma\subset {\cal G}$ of complex gauge transformations
fixing $\db_0$ is the same group as the group that fixes
$d_0$; these are the holomorphic automorphisms
of the holomorphic structure $E_0$ defined by $d_0$. The group
$\Gamma = Aut(E_0)$ acts on the space $H^1(X,End\,E_0)$
of infinitesimal deformations of $E_0$ by conjugation, and
since each element of $\Gamma$ is covariantly constant
with respect to $d_0$, this action preserves the harmonic
subspaces $H^{0,q} = H^{0,q}(\db_0,\eet)$. From
\rDequivariant\ the function
$\Psi$ of \ePsi\ is equivariant
with respect to the action of $\Gamma$ on
$H^{0,1}$ (at least, near $0$), from which it follows that
$\Psi^{-1}(0)\subsett H^{0,1}$ is invariant under
$\Gamma$.
The
linearization of the action of $\Gamma$ on $H^{0,q}$
at $1\in \Gamma$ is given by the Lie bracket,
$H^{0,0}\times H^{0,q}\owns (\varphi,\tau)\mapsto
[\tau,\varphi]$.
The group
$\Gamma$ is the complexification of the subgroup
$U(\Gamma)$ of $d_0$-closed unitary automorphisms of
$\et$, so is a reductive Lie group.

\medskip
\lem{\rFA} A form $\alpha\in H^{0,1}$ is of minimal
norm in its orbit under $\Gamma$  if and only if
$i\,\Lambda(\alpha\w\alpha^*+\alpha^*\w\alpha)\in A^{0,0}(\eet)$ is
orthogonal to $\ker\db_0$.

\pf Given self-adjoint $\delta\in H^{0,0}$, let
$\gamma_t := e^{t\delta}$ for $t\in \m R$.
Then with $\alpha_t := \gamma_t\alpha
\gamma_t^{-1}$,
differentiation with respect to $t$ gives
$\dot\alpha_t = [\delta,\alpha_t]=e^{t\delta}[\delta,\alpha]e^{-t\delta}$,
using the fact that $\delta$ commutes with $e^{t\delta}$. Consequently
$$
\fract{d~}{dt}\norm{\alpha_t}^2
= 2\,\Re\big\<\alpha_t,[\delta,\alpha_t]\big\>
=2\,\Re\big\<e^{2t\delta}\alpha e^{-2t\delta},[\delta,\alpha]\big\>
$$
so
$$
\fract{d^2}{dt^2}\norm{\alpha_t}^2
= 4\,\Re\big\<[\delta,\alpha_t],[\delta,\alpha_t]\big\>\ge 0\;.
$$
Thus any critical point of $t\mapsto \norm{\alpha_t}^2$ is
a minimum, and since $\Re\<\alpha_t,[\delta,\alpha_t]\>
=
-\big\<\delta,i\Lambda(\alpha_t^{}\w\alpha_t^*+\alpha_t^*\w\alpha_t^{})\big>$
(using self-adjointness of $\delta$),
the result follows. \quad \qed

\medskip

The assignment
$$
H^{0,1}\owns \alpha\mapsto m(\alpha) :=
\Pi^{0,0}\,i\,\Lambda\big(\alpha\w \alpha^*+\alpha^*\!\w \alpha\big)
\in H^{0,0} \myeqn
$$
maps elements of the hermitian vector space $H^{0,1}$ into
the space of trace-free self-adjoint elements in $H^{0,0}$,
the latter being $i$ times the Lie algebra of the group
$SU(\Gamma)$
of $d_0$-closed unitary automorphisms of $\et$ with unit
determinant, which is canonically identified with its real dual.
In view of this lemma,  it is natural to presume that $m$ is
a moment map for the action of $SU(\Gamma)$ on $H^{0,1}$,
where the latter is equipped with the symplectic form
$$
H^{0,1}\times H^{0,1}\owns
(\alpha,\beta) \mapsto \underline{\omega}(\alpha,\beta)
:= \Int_X i\Lambda\,\tr\big(\alpha\w\beta^*
-\beta\w\alpha^*\big)\dv=\<\alpha,\beta\>-\<\beta,\alpha\>\in \m R\;,
$$

(with the complex structure $J(\alpha) = i\,\alpha$).
Using the definition in \cite{Nes}, to prove that this is the case,
it must
be shown that $m$ is equivariant with respect to the
action of $SU(\Gamma)$ on $H^{0,1}$ and that
$dm_{\delta}(\dot\alpha)=\underline{\omega}(X_{\delta},\dot\alpha)$
for each skew-adjoint $\delta\in H^{0,0}$ and
each $\dot\alpha\in H^{0,1}$, where $X_{\delta}$
is the vector field on $H^{0,1}$ determined by
$\delta$, this having the value $[\delta,\alpha]$ at
$\alpha\in H^{0,1}$.

\medskip
\lem{\rEeq} Suppose that $\int_X\omega^n=1$. Then for any
$u\in H^{0,0}$ and $v\in A^{0,0}(\eet)$ in $L^2$,
$\tr(u^*\Pi^{0,0}v) = \int_X\tr(u^*v)\dV$.

\pf
If $\varphi,\,\psi \in H^{0,0}$ are arbitrary,
they are both covariantly constant and therefore so too is
$\varphi^*\psi$. Hence the trace of this endomorphism is
constant, which implies that $\tr(\varphi^*\psi)
= \int_X\tr(\varphi^*\psi)\dV = \<\varphi,\psi\>$.
Consequently, if $e_1,\dots, e_m$ is an $L^2$-orthonormal
basis for $H^{0,0}$, the endomorphisms $e_1,\dots, e_m$
are pointwise orthonormal.
Writing $u = \sum_j\<e_j,u\>\,e_j$ and
$\Pi^{0,0}\,v = \sum_k\<e_k,v\>\,e_k$, it follows that
$$
\eqalign{
\tr (u^*\Pi^{0,0}v) = \sum_{j}\<u,e_j\>\,\<e_j,v\>
&= \sum_{j}\<u,e_j\>\,\int_X\tr(e_j^*v)\dV \cr
&= \int_X\sum_{j}\tr(u^*e_j)\,\tr(e_j^*v)\dV
= \int_X\tr(u^*v)\dV\;. \qquad \qed
}
$$

\smallskip
It follows from this lemma that if $v\in A^{0,0}(\eet)$
is in $L^2$ and $u\in SU(\Gamma)$, then
$\Pi^{0,0}(uvu^{-1})=u(\Pi^{0,0}v)u^{-1}$, which implies
that the map $m$ of \lasteq0 is indeed equivariant with respect
to the action of $SU(\Gamma)$. Furthermore, thinking of
$\alpha$ as depending differentiably on a parameter $t\in \m R$
and differentiating at $t=0$, for $\dot\alpha_0 := d\alpha/dt|_{t=0}$
it follows
$$
{d~\over dt}m(\alpha)\Big|_{t=0}
= \Pi^{0,0}i\Lambda\big(\dot\alpha_0\w\alpha^*+\alpha\w\dot\alpha_0^*
+\dot\alpha_0^*\w\alpha+\alpha^*\w\dot\alpha_0\big)\;.
$$
Hence for $\delta\in H^{0,0}$, using \rEeq\ it follows that
$$
\eqalign{
{d~\over dt}\tr\big(m(\alpha)\delta\big)\Big|_{t=0}
&= \int_X\tr i \big(\delta\,i\Lambda(\dot\alpha\w\alpha^*
+\alpha\w\dot\alpha^*
+\dot\alpha^*\w\alpha+\alpha^*\w\dot\alpha)\big) \cr
&= \<[\alpha,\delta^*],\dot\alpha_0\>+\<\dot\alpha_0,[\alpha,\delta]\> \cr
&=\underline{\omega}([\delta,\alpha],\dot\alpha_0)
\quad\hbox{if $\delta=-\delta^*$.}
}
$$
Thus
\smallskip
\cor{\rEmoment} The assignment
$H^{0,1}\owns \alpha\mapsto m(\alpha)
=\Pi^{0,0}i\Lambda\big(\alpha\w\alpha^*+\alpha^*\w\alpha\big)
\in H^{0,0}$ is the moment map for the action of $SU(\Gamma)$
on $H^{0,1}$. \quad\qed

\medskip

Combining \rFA\ with the results of the previous sections
gives some of the main results of this article:

\medskip
\thm{\rDB} Let $d_0$ be an integrable connection on $\et$
with $i\wh F(d_0)=\lambda\,1$, and let
$\Psi$ be the holomorphic function \ePsi\ defined in
a neighbourhood of zero in $H^{0,1}$ with values in $H^{0,2}$.
Let $\alpha\in\Psi^{-1}(0)$ be a non-zero form, and let
$E_{\alpha}$ be the corresponding holomorphic structure.
Then $\alpha$ is polystable
with respect to the action of $\Gamma$
if  $E_{\alpha}$ is polystable. Moreover, when $E_{\alpha}$ is polystable,
$\alpha$ is stable with respect to the action of $\Gamma$
if and only if $E_{\alpha}$ is stable.

\pf The proof is by induction on the rank $r$ of $\et$, the
case $r=1$ being elementary.

Suppose first that $E_{\alpha}$ is stable but
$\alpha$ is not $\Gamma$-polystable. Then
the orbit of $\alpha$ under $\Gamma$ is not closed, and
indeed, the infimum of $\norm{\gamma\cdot\alpha}^2$ over
$\gamma\in \Gamma$ is not attained in that orbit. Let
$\alpha_0\in H^{0,1}$ be a point of smallest
norm in the closure of the orbit of $\alpha$ under $\Gamma$,
unique up to conjugation by unitary elements of $\Gamma$.
So there is a sequence $\gamma_j\in\Gamma$
with $\det\gamma_j=1$ for all $j$ such that
$\norm{\gamma_j\cdot \alpha}^2$ is decreasing to $\norm{\alpha_0}^2$
but $\norm{\gamma_j}$ is not bounded. Let $\beta\in A^{0,2}(\eet)$
be the form orthogonal to $\ker\db_0^*$ determined by
\rAG\ such that
$\db_0+\alpha+\db_0^*\beta$ is an integrable semi-connection
defining the holomorphic structure $E_{\alpha}$
that is $L^p_1$-near to $E_0$,
and let $a'' = \alpha+\db_0^*\beta$. If $\alpha_j :=
\gamma_j\cdot \alpha$ then $\norm{\alpha_j}$ is
decreasing. From \rAG\ it follows that
$\norm{\beta_j}_{L^p_2}$ is uniformly bounded, where
$\beta_j := \gamma_j\cdot\beta$.  If $a_j''
:= \alpha_j+\db_0^*\beta_j = \gamma_j\cdot a''$, it
follows that $\norm{a_j''}_{L^p_1}$ is uniformly bounded,
so after passing to a subsequence if necessary, the
forms $a_j''$ converge weakly in $L^p_1$ and strongly
in $C^0$ to a limit $a_{0}''$. The limiting
connection $d_0+a_{0}$ is integrable, defining a
holomorphic structure $E_{\alpha_0}$ on $\et$.  After
rescaling $\gamma_j$ to $\tilde\gamma_j$ with $\norm{\tilde\gamma_j}
=1$, (a subsequence of) the automorphisms
$\tilde\gamma_j$ converges to a
non-zero limit $\tilde\gamma_0$ with
$\norm{\tilde\gamma_0}=1$ and $\det\tilde\gamma_0=0$,
this defining a holomorphic map from the holomorphic
structure $E_{\alpha}$ defined by $d_0+a$ to the holomorphic structure
defined by $d_0+a_{0}$. Since the latter can be assumed
to be semistable (by \rCB), there is a non-zero holomorphic
map from $E_{\alpha}$ to a semistable bundle of the same
degree and rank that is not an isomorphism, and this
contradicts the assumption that $E_{\alpha}$ is stable.
Therefore $\alpha$ is $\Gamma$-polystable.
If $\alpha$ is polystable but not stable with respect to
the action of $\Gamma$, the isotropy subgroup $\Gamma_{\alpha}
\subset \Gamma$ of
$\alpha$ has dimension greater than one. Then it follows
from 2.\ of \rCK\ that $E_{\alpha}$ is polystable but
not stable, a contradiction.

Suppose now that $E_{\alpha}$ is polystable but not stable.
Then $E_{\alpha}$ splits into a direct sum of stable subbundles,
each of the same slope. In terms of connections, the
bundle-with-connection
$(\et,d_0+a)$ has a unitary splitting into a direct sum of
irreducible unitary bundles-with-connection. Orthogonal projection
onto any of these subbundles (followed by inclusion)
is a holomorphic endomorphism of $E_{\alpha}$, and by
\rCK, such an endomorphism is in fact covariantly constant
with respect to $d_0$ and also commutes with $\alpha$
and $\beta$. Thus the bundle-with-connection
$(\et,d_0)$ has a unitary splitting into a direct sum of
subbundles-with-connection, each of which defines a
polystable subbundle of $E_0$ of the same slope. If
$(\et,d_0)=\sum_j({\rm B}_{j},d_{0,j})$ is this last
splitting, then $i\wh F(d_{0,j})=\lambda\,1$ and
for some skew-adjoint $a_j\in A^1({\rm End\,}B_j)$,
$(\et,d_0+a)=\sum_j({\rm B}_{j},d_{0,j}
+ a_j)$ corresponds to the splitting of $E_{\alpha}$ into
stable components. The compatibility of the splittings
implies that each $a_j$ is of the form $a_j=a_j'+a_j''$
with $a_j''=\alpha_j+\db_{0,j}^*\beta_j$ where
$\alpha_j$ is a $\db_{0,j}$-harmonic $(0,1)$-form with
coefficients in ${\rm End\, B_j}$. By the inductive hypothesis,
each $\alpha_j$ is stable with respect to the action of
$\Gamma_j$, the group of $\db_{0,j}$-closed automorphisms
of $\rm B_j$. Hence there exists such an automorphism
$\rho_j$ such that $\rho_j\cdot\alpha_j=: \hat\alpha_j$ is a zero of
the moment map, which means that
$\Pi^{0,0}_ji\Lambda(\hat\alpha_j^{}\w\hat\alpha_j^*+
\hat\alpha_j^*\w\hat\alpha_j^{})
=0$, where $\Pi^{0,0}_j$ is the orthogonal projection onto
$\ker\db_{0,j}$. But this implies that there is an
automorphism $\sigma\in \Gamma$ such that $\hat\alpha
:= \sigma\cdot\alpha$ satisfies
$\Pi^{0,0}i\Lambda(\hat\alpha\w\hat\alpha^*+
\hat\alpha^*\w\hat\alpha)=0$. Therefore
$\alpha$ is polystable.

The last statement of the proposition follows immediately
from \rCK. \quad \qed

\bigskip
\rem{\kEa}
If $\alpha\in H^{0,1}$ satisfies $\Pi^{0,0}i\Lambda(\alpha\w\alpha^*
+\alpha^*\w\alpha)=0$ and if $\Gamma_{\alpha}\subsett\Gamma$
is the isotropy subgroup of $\alpha$, then it follows from
2.\ of \rCK\ that $\gamma$ commutes with $\beta(t\alpha)$
for $t$ sufficiently small, where $\beta(-)$ is the
function defined implicitly in \rAG. Then by 3.\ of
\rCK\ it follows that $\gamma^*$ also commutes with $\alpha$,
and indeed, it also commutes with $\beta(t\alpha)$. It then
follows from Proposition~1.59 of \cite{K} that $\Gamma_{\alpha}$
is itself a complex reductive group.

The group $\Gamma_{\alpha}$ acts fibrewise on $\et$, splitting
each fibre into a direct sum of
irreducible $\Gamma_{\alpha}$-invariant subspaces. These subspaces
together form subbundles, namely the intersections of the eigen-bundles
associated with the elements of $\Gamma_{\alpha}$.
The fact that $\Gamma_{\alpha}$
is closed under adjoints implies that the splitting of $\et$ into
$\Gamma_{\alpha}$-irreducible subbundles is a {\it unitary\/} splitting.
For $|t|$ sufficiently
small and $|s|\le1$, the connections $d_{t,s}=d_0+a_{t,s}$
given by $a_{t,s}''=t\alpha+s\db_0^*(t\alpha)$ preserve these
splittings, and restrict to irreducible unitary connections on each
of the $\Gamma_{\alpha}$-irreducible subbundles.

\medskip

As mentioned at the end of \S4, \rDB\ is one half of a
local version of the Hitchin-Kobayashi correspondence:
$\omega$-(poly)stability of $E_{\alpha}$ implies
$\Gamma$-(poly)stability of $\alpha$, where the
latter term means polystable with respect to the action
of $\Gamma$ in the sense of geometric invariant
theory. The more difficult
task is to establish the other half; that is, the converse,
and this effectively involves solving differential
equations. This will be the subject of the next three sections.

\bigskip

\secn{Connections with constant central curvature.}

As in previous sections, let $(X,\omega)$ be a compact K\"ahler
$n$-manifold
and let $\et$ be a fixed complex $r$-bundle equipped with a fixed
hermitian metric; all conventions and notations from
previous sections are also retained.
In this section, the study of
\S2 into an $L^p_1$ neighbourhood of a given (hermitian)
connection will be
continued but
the focus is now is on the  central component
$\wh F=\Lambda F$ of the curvature $F$ rather than the $(0,2)$-component.
As previously, ${\cal G}$ is the group of complex automorphisms
of $\et$, with ${\cal U}\subset {\cal G}$ the subgroup preserving
the given hermitian metric. Let ${\cal A}$ denote the space of
 hermitian connections $d=\d+\db$ on $\et$, so the action
of ${\cal G}$ on ${\cal A}$ is given by
$$
{\cal G}\times {\cal A}\owns (g,d)\;\mapsto\; g\cdot d
:= g^*{}^{-1}\circ \d\circ g^*+g\circ \db\circ g^{-1}
\;=\; d+g^*{}^{-1}\d g^*-\db gg^{-1}\;\in\; {\cal A}\;.
\myeqn
$$
Unless otherwise stated, elements of ${\cal G}$ are assumed to lie
in $L^p_2$ and elements of ${\cal A}$ to lie in $L^p_1$.
Projection to the central component of the curvature
defines a function $\Phi$ on ${\cal G}\times {\cal A}$ with values
in the space of self-adjoint endomorphisms of $\et$ lying in
$L^p$ given by
$$
\Phi \: {\cal G}\times {\cal A}\to A^{0,0}(\eet)\;, \quad
\Phi(g,d) := i\Lambda F(g\cdot d)\;, \myeqn
$$\eqtag{\eFtwo}
and it is the properties of this function and its derivatives
with respect to each of its arguments on which the analysis
concentrates in this section.

\smallskip
Let $d_0=\d_0+\db_0$ be a connection on $\et$, and let
$a=a'+a''$ be an element of $A^1(\eet)$ with $a'=-(a'')^*$.
If $d_a= d_0+a$, the curvature of this connection is
$$
F(d_a)=F(d_0)+d_0a+a\w a\;. \myeqn
$$\eqtag{\eSixthree}
It follows that if $a=a(t)$ depends differentiably on the real
parameter $t$, then
$$
{d~\over dt}\big[F(d_a)\big]= d_0\dot a+\dot a\w a+a\w \dot a
= d_a\dot a\;, \myeqn
$$
where $\cdot$ denotes differentiation with respect to $t$.

The action of
${\cal G}$ on ${\cal A}$ has the explicit form
$$
g\cdot d_a = d_0 + \big(g^*{}^{-1}a'g^*+g^*{}^{-1}\d_0g^*\big)
+\big(ga''g^{-1}-\db_0gg^{-1}\big)\;, \quad g\in {\cal G}\;. \myeqn
$$
If $g=g(t)$ also depends differentiably on $t$, then by direct
calculation from \lasteq0, it follows that
$$
{d~\over dt}\big[g\cdot d_a\big] =
\big(g^*{}^{-1}\dot a'g^*+\d_{g\cdot d_a}(g^*{}^{-1}\dot g^*)\big)
+
\big(g\dot a''g^{-1}-\db_{g\cdot d_a}(\dot gg^{-1})\big)\;.
\myeqn
$$\eqtag{\eFaz}
Hence, from \lasteq2\ and \lasteq0, it follows that
$$
\fract{d~}{dt}\big[F(g\cdot d_a)\big] =
d_{g\cdot d_a}\big(g^*{}^{-1}\dot a' g^*+
g\dot a''g^{-1}+\d_{g\cdot d_a}(g^*{}^{-1}\dot g^*)
-\db_{g\cdot d_a}(\dot g g^{-1})\big)\;. \myeqn
$$
Applying $i\Lambda$ to both sides and recalling that
$\d^*=i\Lambda\db$ and $\db^*=-i\Lambda\d$ on $1$-forms,
$$
\eqalignno{
\fract{d~}{dt}\big[i\wh F(g\cdot d_a)\big]
&= i\Lambda\db_{g\cdot d_a}(g^*{}^{-1}\dot a' g^*)
+i\Lambda\d_{g\cdot d_a}(g\dot a''g^{-1})
+i\Lambda\db_{g\cdot d_a}\!\d_{g\cdot d_a}(g^*{}^{-1}\dot g^*)
-i\Lambda\d_{g\cdot d_a}\!\db_{g\cdot d_a}(\dot g g^{-1}) \cr
&= \d^*_{g\cdot d_a}(g^*{}^{-1}\dot a' g^*)
-\db^*_{g\cdot d_a}(g\dot a''g^{-1})
+\d^*_{g\cdot d_a}\!\d^{}_{g\cdot d_a}(g^*{}^{-1}\dot g^*)
+\db^*_{g\cdot d_a}\!\db^{}_{g\cdot d_a}(\dot g g^{-1})
&\eqn
}
$$
Since $\Lambda(\db\d+\d\db)=\wh F$, writing $\sigma := \dot gg^{-1}$
and decomposing into self-adjoint and skew-adjoint components
gives the following conclusion, which will be used frequently
in this section:


\medskip
\lem{\rFagn} Let $a=a(t)\in A^1(\eet)$ be a differentiable
$1$-real parameter family of skew-adjoint forms, and let
$g=g(t)$ be a differentiable $1$-real parameter family of
complex automorphisms of $\et$. If $d_0=\d_0+\db_0$ is a
connection on $\et$ and $d_a=d_0+a$, then
$$
\fract{d~}{dt}\big[i\wh F(g\cdot d_a)\big]=
\d^*_{g\cdot d_a}(g^*{}^{-1}\dot a' g^*)
-\db^*_{g\cdot d_a}(g\dot a''g^{-1})
+\lap_{g\cdot d_a}\sigma_+ -[i\wh F(d_{g\cdot d_a}),\sigma_-]\;, \myeqn
$$\eqtag{\eEaa}
where $\sigma := \dot gg^{-1}$ and
$\sigma_{\pm} := {1\over 2}(\sigma\pm\sigma^*)$.
\quad \qed

\medskip

Taking $a\in A^1(\eet)$ to be independent of $t$, it
follows from \lasteq0 that the linearization of the
map ${\cal G}\owns g\mapsto i\wh F(g\cdot d_a)= \Phi(g,d_a)$
at a connection $g_0\cdot d$ is
$$
(D_1\Phi)_{(g_0,d_a)}(\sigma) =
\lap_{g_0\cdot d_a}\sigma_+-[i\wh F(g_0\cdot d_a),\sigma_-],\,
\quad \sigma\in A^{0,0}(\eet)\;.
$$
In particular, if $g_0=1$ and $a=0$,
the linearization at $(1,d_0)\in
{\cal G}\times {\cal A}$ of the action of ${\cal G}$ on the
space of $L^p_1$ (hermitian) connections
is an isomorphism from the space of $L^p_2$ self-adjoint sections
of $\eet$ that are orthogonal to the $d_0$-closed
sections to the space of self-adjoint $L^p$ sections of
$\eet$ that are again orthogonal to the kernel of $d_0$.

\smallskip
From now on, let $d_0$ be a connection with
$i\wh F(d_0)=\lambda\,1$, so by \rCCa, $\ker d_0=\ker\db_0$.
In general, the map ${\cal G}\times A^1(\eet)\owns (g,a)\mapsto
\Phi(g,d_a)= i\,\widehat F(g\cdot d_a)$
takes values in the self-adjoint endomorphisms of $\et$, but it
does not map into
the space orthogonal to $\ker \db_0$. However, if $\Pi^{0,0}$ is the
$L^2$ projection of
$A^{0,0}(\eet)$ onto $\ker \db_0$ and
$\Pi_{\;\perp}^{0,0} = 1-\Pi^{0,0}$
is the projection onto
the orthogonal complement, then for any skew-adjoint $L^p_1$ section
$a\in A^1(\eet)$, the composition
$$
A^{0,0}(\eet) \owns \varphi \mapsto
\Pi_{\;\perp}^{0,0}\,i\,\widehat F\big(\exp(\varphi)\cdot(d_0+a)\big)
\in A^0(\eet)
$$
maps the Banach space of self-adjoint $L^p_2$ sections 
in $A^{0,0}(\eet)$
orthogonal to $\ker\db_0$ into the Banach space of
such sections lying in $L^p$,
and when $a=0$, this map has the same linearization at
$0$ as the earlier map.  The implicit function theorem
for Banach spaces then yields:

\penalty-100


\prop{\rEA} There exists $\epsilon>0$ depending
only on $d_0$ with the property that for each skew-adjoint
section $a\in A^1(\eet)$ satisfying $\norm{a}_{L^p_1} < \epsilon$
there is a unique self-adjoint $\varphi \in
(\ker \db_0)^{\perp}\subset A^{0,0}(\eet)$
for which
$\Pi^{0,0}_{\perp}\,i\wh F\big(\exp(\varphi)\cdot(d_0+a)\big)=0$.
Furthermore, there is a constant $C$ depending only on $d_0$
such that $\varphi$ satisfies
$\norm{\varphi}_{L^p_2}\le
C\,\nnorm{\Pi_{\;\perp}^{0,0}\,\Lambda(d_0a+a\w a)}_{L^p}$.
\quad\qed \rm

\medskip
The power series expansion in $\varphi$ at $\varphi=0$ of
$i\,\wh F(\exp(\varphi)\cdot d_a)$ is, to first order,
$$
\eqalignno{
i\wh F\big(&\exp(\varphi)\cdot(d_0+a)\big)
= \;i\wh F(d_0)+i\Lambda(d_0a+a\w a)
-[i\wh F(d_0)+i\Lambda(d_0a+a\w a),\varphi_-] \cr
&\; + \lap_0\varphi_++2i\Lambda\big(a''\w
\d_0\varphi_++\d_0\varphi_+\w a''-
a'\w\db_0\varphi_+-\db_0\varphi_+\w a'\big)+[d_0^*a,\varphi_+] \cr
&\qquad+i\Lambda\big([a',\varphi_+]\w a'' + a''\w[a',\varphi_+]
-a'\w[a'',\varphi_+]-[a'',\varphi_+]\w a'\big)
+R_2(\varphi)\;,&\eqn
}
$$
where the term $R_2(\varphi)$ involves products of $\varphi$ and
its first and second derivatives with respect to $d_a$
with at least two such factors, but where the second-order
derivatives appear linearly and the first-order derivatives
appear at most quadratically. Consequently, if $a$ satisfies
the hypotheses of \rEA\ and if $\varphi\in A^{0,0}(\eet)$
satisfies the inequality in the statement of that proposition,
then $\norm{R_2(\varphi)}_{L^p} \le C\,\norm{\Lambda(d_0a+a\w a)}_{L^p}^2$
for some constant $C = C(d_0)$.

If $\varphi$ is self-adjoint and if $\db_0^*a''=0$, the
formula \lasteq0 simplifies somewhat. In this case, it is useful
to project both sides into $\ker\db_0$ and
$(\ker\db_0)^{\perp}$, respectively giving
$$
\eqalignno{
\Pi^{0,0}i\wh F\big(\exp(\varphi)\cdot d_a\big)
&= \lambda\,1 + i\,\Pi^{0,0}\,\Lambda\big(a'\w a''+a''\w a'\big)
 &\eqn{\rm (a)}\cr
&\qquad +i\,\Pi^{0,0}\Lambda\big([a',\varphi]\w a''+a''\w[a',\varphi]
  -[a'',\varphi]\w a'-a'\w [a'',\varphi]\big)+\Pi^{0,0} R_2(\varphi)\;,\cr
\strutt{13}0 \global\advance\eqcount by -1
  \Pi^{0,0}_{\;\perp}i\wh F\big(\exp(\varphi)\cdot d_a\big)
&= \Pi^{0,0}_{\;\perp}\,i\,\Lambda\big(a'\w a''+a''\w a'\big)  &\eqn{\rm (b)}\cr
&\;\;\quad+\lap_0\varphi+2i\Lambda\big(a''\w \d_0\varphi+\d_0\varphi\w a''
 -a'\w\db_0\varphi-\db_0\varphi\w a'\big) \cr
&\quad\quad
  +\Pi^{0,0}_{\;\perp}i\Lambda\big([a',\varphi]\w a''+a''\w[a',\varphi]
  -[a'',\varphi]\w a'-a'\w [a'',\varphi]\big)
   +\Pi^{0,0}_{\;\perp}R_2(\varphi)\;.
}
$$\eqtag{\eEaz}
If $\varphi=\varphi(a)$ is the endomorphism of \rEA, the left-hand
side of \lasteq0(b) vanishes, giving the equation that
effectively determines
$\varphi(a)$. Since $\Lambda d_0a=0$ now, the uniform estimate
on $\varphi$ provided by that proposition then implies
$\norm{\varphi(a)}_{L^p_2}\le C\norm{a}_{L^{2p}}^2$, so
the remainder term $R_2(\varphi)$ appearing in \lasteq0
is uniformly
bounded by a constant multiple of
$\norm{a}_{L^{2p}}^4$, something that is
also true of the other terms on the second line
of \lasteq0(a).
If $a''=\alpha+\db_0^*\beta$ for $\db_0$-harmonic
$\alpha\in A^{0,1}(\eet)$ and $\beta\in(\db_0^*)^{\perp}
\subsett A^{0,2}(\eet)$,
the (negative of the) terms involving
$a'\w a''+a''\w a'$ in \lasteq0 expand to
$$
\eqalignno{
i\,\Lambda(\alpha^*\!\w\alpha+\alpha\w\alpha^*)
+i\,\Lambda(&\d_0^*\beta^*\!\w\db_0^*\beta+\db_0^*\beta\w\d_0^*\beta^*) \cr
&+i\,\Lambda(\alpha^*\!\w\db_0^*\beta+\d_0^*\beta^*\!\w\alpha
+\alpha\w\d_0^*\beta^*+\db_0^*\beta\w\alpha^*)\;. &\eqn
}
$$
If $\sigma\in A^{0,0}(\eet)$ is $\db_0$-closed, then $\d_0\sigma=0$ by
\rCCa,
and
$$
\big\<\sigma,i\,\Lambda(\alpha^*\!\w\db_0^*\beta+\d_0^*\beta^*\!\w\alpha
+\alpha\w\d_0^*\beta^*+\db_0^*\beta\w\alpha^*)\big\>
=
\<[\alpha,\sigma],\db_0^*\beta\>+\<[\db_0^*\beta,\sigma]),\alpha\>\;.
\myeqn
$$\eqtag{\eEk}
The first term on the right vanishes because $\db_0[\alpha,\sigma]= 0$,
and the second term vanishes
because $[\db^*_0\beta,\sigma]= \db^*_0[\beta,\sigma]$; consequently
the projection of the term in the second line of \lasteq1 on
$\ker\db_0$ is zero. To summarize the discussion so far:

\medskip
\lem{\rFzz} Suppose that $a=a'+a''$ satisfies the hypotheses of
{\rm\rEA} and $a''=\alpha+\db_0^*\beta$ where
$\alpha$ is $\db_0$-harmonic. If $\varphi=\varphi(a)$ is the
endomorphism of that proposition, then
$$
\eqalignno{
i\wh F\big(\exp(\varphi)\cdot&(d_0+a)\big) -\lambda\,1 =
   -i\,\Pi^{0,0}\,\Lambda(\alpha^*\w\alpha+\alpha\w\alpha^*)
   -i\,\Pi^{0,0}\,
   \Lambda(\d_0^*\beta^*\w\db_0^*\beta+\db_0^*\beta\w\d_0^*\beta^*) \cr
+\;i\,&\Pi^{0,0}\Lambda\big([a',\varphi]\w a''+a''\w[a',\varphi]
   -[a'',\varphi]\w a'-a'\w [a'',\varphi]\big)\;+R(a) &\eqn
}
$$\eqtag{\eEC}
where $\norm{R(a)}_{L^p} \le C\norm{a}_{L^{2p}}^4$ for some constant
$C=C(d_0)$. \quad \qed

\medskip
The term $\Pi^{0,0}\Lambda(\alpha^* \w \alpha+\alpha\w \alpha^*)$
is $O(\norm{\alpha}^2)$ in general, whereas if
$\beta$ is as in \rAG,
all the other terms on the right of \lasteq0 are $O(\norm{\alpha}^4)$.
But if $\Pi^{0,0}\Lambda(\alpha^* \w \alpha+\alpha\w \alpha^*)$
vanishes (as considered in \S5) then
the whole of the right-hand side of \lasteq0 is $O(\norm{\alpha}^4)$
and the connection $\exp(\varphi(a))\cdot (d_0+a)$ is very
close to having central curvature equal to $-i\,\lambda\,1$.
Given that $\varphi(a)$ is orthogonal to $\ker\db_0$, it
can be hoped that a small perturbation by an element of
$\ker\db_0$ will yield a connection with $i\wh F \equiv\lambda\,1$.
For $\gamma\in\Gamma$, $\gamma\cdot (d_0+a)=d_0+\gamma\cdot a
:= d_0+\gamma^*{}^{-1}a'\gamma^*+\gamma a''\gamma^{-1}$, so
since $\exp(\varphi+\delta)$
is close to $\exp(\varphi)\exp(\delta)$ for small
$\varphi\in (\ker\db_0)^{\perp}$
and small $\delta
\in \ker\db_0$,  an alternative is to perturb $\alpha$ by
conjugation with an element
of $\Gamma$ close to $1$.
Such an argument will involve an application of the inverse
function theorem in finite dimensions, for which purpose
the variation in $\wh F$ as $\alpha$ is
varied in this way must be determined.

With  `$sk$' denoting skew-adjoint and `$sa$'
denoting self-adjoint, consider first the function
$G$ defined on $A^{0,0}(\eet)\times A^1_{sk}(\eet)$
with values in $A^{0,0}_{sa}(\eet)$ defined by
$$
G(\psi,a) := i\wh F(\exp(\psi)\cdot(d_0+a)) \qquad
\hbox{($=\Phi(\exp(\psi),d_0+a)$ for $\Phi$ as in \eFtwo).}
$$
If $G_0 := \Pi^{0,0}G$ and $G_1 := \Pi^{0,0}_{\;\perp}G$,
the conclusion of \rEA\ is that for $a\in A^1_{sk}(\eet)$
with $\norm{a}_{L^p_1}<\epsilon$, $G_1(\varphi(a),a)\equiv 0$,
where $a\mapsto \varphi(a)$ is the function specified
in that proposition, the existence of which is guaranteed by
the implicit function theorem.

If $a$ moves in a differentiable
$1$-parameter family $a(t)$, then it follows that
$$
0\equiv {d~\over dt}\Big[G_1(\varphi(a(t)),a(t))\Big]
= \big((D_1G_1)_{(\varphi(a),a)}\circ (D\varphi)_{a}\big)(\dot a)
+ (D_2G_1)_{(\varphi(a),a)}(\dot a)\;,
$$
where $D_1G_1$, $D_2G_1$ are respectively the partial derivatives
of $G_1$ with respect to its first and second arguments and
$\dot a$ denotes the derivative with respect to $t$ as before.
The implicit function theorem implies that if
$(\psi,a)$ is sufficiently
close to $(0,0)$, then $(D_1G_1)_{(\psi,a)}$ is
an isomorphism from the space of self-adjoint elements of $A^{0,0}(\eet)$
orthogonal to $\ker\db_0$ lying in $L^p_2$ to the space of
self-adjoint elements of $A^{0,0}(\eet)$ orthogonal to
$\ker\db_0$ lying in $L^p$, so
$$
(D\varphi)_{a}(\dot a)
= -\big((D_1G_1)_{(\varphi(a),a)}^{-1}\circ
(D_2G_1)^{}_{(\varphi(a),a)}\big)(\dot a)\;.
$$
The variation in $G(\varphi(a),a)$ at
$a$ is therefore given by
$$
\eqalignno{
{d~\over dt}\Big[G(\varphi(a),a)\Big]
&= \hphantom{-}\big((D_1G)_{(\varphi(a),a)}\circ
(D\varphi)_{a}\big)(\dot a)
+ (D_2G)_{(\varphi(a),a)}(\dot a) &\eqn \cr
&=-\big((D_1G_0)_{(\varphi(a),a)}\circ
(D_1G_1)_{(\varphi(a),a)}^{-1} \circ
(D_2G_1)_{(\varphi(a),a)}\big)(\dot a)
+ (D_2G_0)_{(\varphi(a),a)}(\dot a)\;.
}
$$\eqtag{\eFvar}


The partial derivatives appearing here
are determined by the formula \eEaa:
setting $d_b := \exp(\varphi(a))\cdot (d_0+a)$
and substituting $g=e^{\varphi(a)}$
into \eEaa\ gives
$$
\eqalignno{
(D_2G)_{(\varphi(a),a)}(\dot a)
&= \d_b^*(e^{-\varphi(a)}\dot a' e^{\varphi(a)})
-\db_b^*(e^{\varphi(a)}\dot a'' e^{-\varphi(a)}) \cr
&= e^{\varphi(a)}\d_a^*\big(e^{-2\varphi(a)}\dot a' e^{2\varphi(a)}
\big)e^{-\varphi(a)}
-e^{-\varphi(a)}\db_a^*\big(e^{2\varphi(a)}\dot a'' e^{-2\varphi(a)}\big)
e^{\varphi(a)}\,, &\eqn
}
$$\eqtag{\eEg}
where $\dot a=\dot a' +\dot a''$,
with the derivatives $D_2G_0$ and $D_2G_1$ obtained by projecting into
$\ker\db_0$ and its orthogonal complement respectively.
Similarly, the partial derivative of $G$ with respect to
its first variable $\psi$ at $(\varphi(a),a)$ is obtained from \eEaa\
by substituting $\dot a=0$ and
$\sigma=({de^{\psi}}\!/{dt})\,e^{-\psi}\big|_{\psi=\varphi(a)}$ for
a $1$-parameter family $\psi(t)$  into that formula, so
$$
(D_1G)_{(\varphi,a)}(\sigma)=
\db_b^*\db_b^{}\sigma+\d_b^*\d_b^{}\sigma^*
=\lap_b\,\sigma_+-[i\wh F(d_b),\sigma_-]\,,
\myeqn
$$\eqtag{\eEp}
again with the derivatives of $G_0$ and $G_1$ obtained by
taking $L^2$ projection into $\ker\db_0$ and its orthogonal
complement.
Note that
$\sigma$ is not self-adjoint in general, but satisfies
$\sigma^* = e^{-\varphi(a)}\sigma e^{\varphi(a)}$. Thus the
second term on the right of \lasteq0 is not zero in general,
unless $i\wh F(d_b)=\lambda\,1$.

\medskip
\lem{\rEC} Under the hypotheses of {\rm\rEA}, suppose in addition
that $\db_0^*a''=0$. Then there is a constant $C=C(d_0)$
such that, for skew-adjoint $\dot a\in A^1(\eet)$,
$$
\nnorm{\big((D_1G_0)_{(\varphi(a),a)}\circ
(D_1G_1)_{(\varphi(a),a)}^{-1} \circ
(D_2G_1)^{}_{(\varphi(a),a)}\big)(\dot a)}_{L^2}
\le C\,\norm{a}_{L^p_1}^3\,\big(\norm{\db_0^*\dot a}_{L^p}
+\norm{\dot a}^{}_{L^p}\big)\;. \myeqn
$$\eqtag{\eEcc}

\pf
The assumption
that $\db_0^*a''=0$ implies that $\Lambda d_0a=0$, so
the bound of \rEA\ implies that
$\norm{\varphi}_{L^p_2}\le C\norm{a}_{L^{2p}}^2$
for some uniform constant $C$. By the Sobolev embedding
theorem, there is a similar such bound on the $C^1$ norm
of $\varphi$, so (given that $\norm{a}_{L^p_1}$ is sufficiently
small),  an arbitrary endomorphism $\psi\in A^{0,0}(\eet)$
will satisfy a pointwise bound of the form
$\big|\psi-e^{\varphi}\psi e^{-\varphi}\big|
\le C\norm{a}_{L^{2p}}^2\,|\psi|$, which can be seen
by orthogonally diagonalizing $\varphi$ at the point
in question.  Since
$(\db_0e^{\varphi})e^{-\varphi}$ is uniformly bounded
in $C^0$ by a multiple of $\norm{a}_{L^{2p}}^2$, it
follows that
$b''=e^{\varphi}a''e^{-\varphi}-
(\db_0e^{\varphi})e^{-\varphi}$ also satisfies
a pointwise bound of the form $|b''-a''|\le
C\norm{a}_{L^{2p}}^2$, and therefore
$(\db_{a}e^{\varphi})e^{-\varphi}
=a''-b''$ satisfies this bound; similarly,
$e^{-\varphi}\db_ae^{\varphi}$ also satisfies such a bound.
Consequently, for any $\chi\in A^{0,1}(\eet)$,
there is a uniform pointwise
bound of the form
$$
\big|\,e^{-\varphi}\db_a^*\big(e^{2\varphi}
\chi e^{-2\varphi}\big)
e^{\varphi}
- \db_a^*\chi\,\big|
\le C\norm{a}_{L^{2p}}^2\big(|\db_a^*\chi|+|\chi|\big)\;,
$$
and since $\db_a^*\chi =\db_0^*\chi -i\Lambda(a'\w\chi+\chi\w a')$,
$$
\big|\,e^{-\varphi}\db_a^*\big(e^{2\varphi}
\chi e^{-2\varphi}\big)
e^{\varphi}
- \db_a^*\chi\,\big|
\le C\norm{a}_{L^{p}_1}^2\big(|\db_0^*\chi|+|\chi|\big) \myeqn
$$\eqtag{\eEdd}
for some new uniform constant $C$.  Taking $\chi=\dot a$ in \eEg,
it follows that
$$
\big|(D_2G)_{(\varphi,a)}(\dot a)\big|
\le C\norm{a}^2_{L^{p}_1}\big(|\db_0^*\dot a|+|\dot a|\big) \quad
{\rm pointwise.} \myeqn
$$
It follows from this that the $L^2$ norm of
$(D_2G)_{(\varphi,a)}(\dot a)$
is uniformly bounded above by a constant multiple of
$\norm{a}^2_{L^p_1}\big(\norm{\db_0^*\dot a}_{L^2}+\norm{\dot a}_{L^2}\big)$,
which implies the same such bound for its orthogonal projection
onto $\ker\db_0$. Since $\ker\db_0$ is finite dimensional,
the $L^2$ norm on this space is equivalent to any other, so it follows
from \lasteq0 that there is
uniform bound of the form
$$
\nnorm{(D_2G_1)_{(\varphi,a)}(\dot a)}_{L^p}
\le C\norm{a}^2_{L^p_1}
\big(\norm{\db_0^*\dot a}_{L^p}+\norm{\dot a}_{L^p}\big)\;.
\myeqn
$$\eqtag{\eFfour}

From the proof of \rEA\ using the implicit function theorem, the operator
$(D_1G_1)_{(\varphi,a)}$
is an isomorphism from the space of
self-adjoint $L^p_2$
sections of $A^{0,0}(\eet)$ orthogonal to $\ker\db_0$
to the same such space of $L^p$ sections, so
$$
\nnorm{(D_1G_1)_{(\varphi,a)}^{-1}
\big((D_2G_1)^{}_{(\varphi,a)}\big)(\dot a)\big)}_{L^p_2}
\le C\nnorm{(D_2G_1)_{(\varphi,a)}(\dot a)}_{L^p} \myeqn
$$
for some new constant $C=C(d_0)$.

For a section $\sigma\in A^{0,0}(\eet)$, \eEp\ gives
$(D_1G_0)_{(\varphi,a)}(\sigma)=
\Pi^{0,0}\big(\db_b^*\db_b^{}\sigma+\d_b^*\d_b^{}\sigma^*\big)$,
where $d_b$ is the connection $d_b= e^{\varphi(a)}\cdot d_a$.
Since $\db_b^*\db_b^{}\sigma
= \db_0^*\db_b^{}\sigma-i\Lambda(b'\w\db_b\sigma
+\db_b\sigma\w b')$ for which the first term is
annihilated by $\Pi^{0,0}$, estimation of
the second term gives
$\nnorm{\Pi^{0,0}\db_b^*\db_b^{}\sigma}_{L^2}
\le C\big(\norm{b}_{L^2}\norm{\sigma}_{L^2_1}
+ \norm{b}_{L^4}^2\norm{\sigma}_{L^2}\big)$.
The a priori $L^p_1$ bound on $a$ and the estimate on
$\varphi(a)$ from \rEA\ implies that
$\norm{b}_{L^{2p}}$ is uniformly bounded above by
a multiple of $\norm{a}_{L^{2p}}$, so
$\nnorm{\Pi^{0,0}\db_b^*\db_b^{}\sigma}_{L^2}
\le C\norm{a}_{L^{2p}}\norm{\sigma}_{L^2_1}$.
Taking adjoints, this same estimate applies to
$\Pi^{0,0}\d_b^*\d_b^{}\sigma^*$ to give
$\nnorm{(D_1G_0)_{(\varphi(a),a)}(\sigma)}_{L^2}
\le C\norm{a}_{L^{2p}}\norm{\sigma}_{L^2_1}$. Combining this
with the estimates of \lasteq0 and \lasteq1 gives \eEcc.
\quad \qed

\medskip

In accordance with the strategy outlined following the statement
of \rFzz, \rEC\ provides sufficient information to analyse the variation
of $i\wh F(e^{\varphi(a)}\cdot d_a)$ as $a''=\alpha+\db_0^*\beta(\alpha)$ 
is varied
according to  $H^{0,1}\owns \alpha\mapsto \gamma\alpha\gamma^{-1}$ for
$\gamma\in \Gamma$, at least when the connections $d_a$ are
sufficiently near to $d_0$; here the notion of ``near" here depends
on the class $\alpha\in H^{0,1}$, which is the explanation for the 
reference to {\it near\/}
in the title of the next section.

\bigskip

\secn{$\Gamma$-stable implies nearly $\omega$-stable.}

Retaining all of the notation and definitions of the previous
section, suppose now that $\gamma_t\in\Gamma$ is a family
depending differentiably on the real variable $t$,
with $\gamma_0=1$ and $\dot\gamma_t|_{t=0}=\delta\in
H^{0,0}$. Suppose that $a=a'+a''$ satisfies the hypotheses
of \rEA\ as well as $\db_0^*a''=0$,
and let $a_t  := \gamma_t\cdot a = \gamma_t^{} a''\gamma_t^{-1}
+\gamma_t^*{}^{-1}a'\gamma_t^*$. Then at $t=0$,
$\dot a''_0 := \dot a''_t|_{t=0}=-[a'',\delta]$, so
$\db_0^*\dot a_0''=0$ and \rEC\ gives an estimate
of the contribution of the first term on the right
of \eFvar\ to the variation in $i\wh F$ in terms of
$\nnorm{[a'',\delta]}_{L^p}$. But since $\delta$
is $d_0$-closed, its $C^0$ norm is bounded by a uniform
multiple of its $L^2$ norm so
 $\norm{[a'',\delta]}_{L^p} \le C\norm{a}_{L^p}
\norm{\delta}$ for $C=C(d_0)$, and therefore that contribution is
uniformly bounded above by a constant multiple of
$\norm{a}_{L^p_1}^4\norm{\delta}$, (the fourth power
coming from \eEcc).
It follows that if
$\norm{a}_{L^p_1}$ is sufficiently small, the dominant
term in the variation \eFvar\
of $i\wh F$ is that coming from
$(D_2G_0)_{(\varphi_0,a)}\big([a'',\delta]_+\big)$, given by the projection
of \eEg\ onto $\ker\db_0$, provided that this is appropriately
non-degenerate as a function of $\delta$.

As in the proof of \rEC,
\eEdd\ gives a bound
$$
\nnorm{\,
e^{-\varphi}\db_a^*\big(e^{2\varphi}[a'',\delta] e^{-2\varphi}\big)
e^{\varphi}-\db_a^*[a'',\delta]\,}_{L^2}
\le C\,\norm{a}_{L^{p}_1}^2\,\nnorm{[a'',\delta]}\;.
$$
Noting that $\dot a''=-[a'',\delta]$
and $\db_a\delta=[a'',\delta]$, it follows that
$$
\eqalign{
\big\<\delta,e^{-\varphi}\db_a^*\big(e^{2\varphi}[a'',\delta] e^{-2\varphi}\big)
e^{\varphi}\big\>
&\ge \big\<\delta,\db_a^*[a'',\delta]\big\>
- C\norm{a}^2_{L^{p}_1}\,\nnorm{[a'',\delta]}\,\nnorm{\delta} \cr
&= \nnorm{[a'',\delta]}^2
- C\norm{a}^2_{L^{p}_1}\,\nnorm{[a'',\delta]}\,\nnorm{\delta}\cr
&\ge {1\over 2}\nnorm{[a'',\delta]}^2-C'\nnorm{a}_{L^p_1}^4\,\nnorm{\delta}^2\;.
}
$$
By taking adjoints and using $\dot a' = +[a',\delta^*]$,
the same estimates apply
to the other term in \eEg\ with $\delta$ replaced by
$\delta^*$. Then if $\delta$ is taken to be self-adjoint,
by combining these estimates with those of \rEC\ the
following conclusion is reached:

\medskip
\prop{\rED} Under the hypotheses of {\rm\rEA}, assume in addition
that $\db_0^*a''=0$. Let
$\wh \Phi$ be the function defined in a neighbourhood
of $1\in \Gamma$ with values in $H^{0,0}$ given
by $\wh \Phi(\gamma) :=
 i\wh F(\varphi(\gamma\cdot a)\cdot\gamma\cdot (d_0+a))$,
where $\varphi(-)$ is the function of {\rm \rEA}.
Then there is a constant $C=C(d_0)>0$ such that
$$
\nnorm{\,[a'',\delta]\,}^2
\le 2\,\big\<\delta,(D\wh \Phi)(\delta)\big\>
+C\norm{a}_{L^p_1}^4\norm{\delta}^2
$$
for any self-adjoint $\delta\in H^{0,0}$. \quad \qed

\medskip

With $a''$ and $\delta$ as in this proposition,
write $a''=\alpha+\db_0^*\beta$
where $\alpha\in H^{0,1}$ and $\beta\in A^{0,2}(\eet)$.
Given that $d_0\delta=0$, $[a'',\delta]=
[\alpha,\delta]+[\db_0^*\beta,\delta]=
[\alpha,\delta]+\db_0^*[\beta,\delta]$, and
since this is an orthogonal decomposition, it follows
that $\nnorm{[a'',\delta]}^2=\nnorm{[\alpha,\delta]}^2
+\nnorm{\db_0^*[\beta,\delta]}^2$.
Assuming that $\beta$ satisfies the uniform bound of \rAG, from
\rEA\ it follows that there is a bound of the form
$\norm{a}_{L^p_1} \le C\,\norm{\alpha}$ for some
constant $C=C(d_0)$. Hence for some new constant $C=C(d_0)$,
the bound of \rEC\ implies a bound of the form
$$
\nnorm{[\alpha,\delta]}^2 +\nnorm{\db_0^*[\beta,\delta]}^2
\le 2\,\big\<\delta,(D\wh\Phi)(\delta)\big\>
+C\norm{\alpha}^4\norm{\delta}^2
\quad\hbox{for $\delta=\delta^*\in H^{0,0}$.} \myeqn
$$\eqtag{\eEx}

\medskip
This estimate has been derived under the assumption that
$\beta$ satisfies the uniform bound given in \rAG. In
particular, it applies if $\beta=\beta_0(\alpha)$ where
$\beta_0(\alpha)$ is the unique element of $A^{0,2}(\eet)$
specified in that result, but more generally, it also
applies if $\beta=s\,\beta_0(\alpha)$ where $s\in [0,1]$.
This observation facilitates some homotopy arguments to
follow.

\smallskip
The leading term on the right of \eEC\ is the
negative of
$m(\alpha) := \Pi^{0,0}i\Lambda(\alpha\w\alpha^*+\alpha^*\w\alpha)$.
Fixing $\alpha$ temporarily, this gives a map
$\Gamma \to H^{0,0}$ give by $\Gamma\owns \gamma\mapsto m(\gamma\alpha\gamma^{-1})$,
and the derivative of this map
at $\gamma=1$ is given by
$$
\big(D_{\gamma}\big[m(\gamma\alpha\gamma^{-1})\big]\big)
\big|_{\gamma=1}(\delta)=
\Pi^{0,0}i\Lambda\big([\delta,\alpha]\w\alpha^*
+\alpha\w[\alpha^*,\delta^*]+
[\alpha^*,\delta^*]\w\alpha+\alpha^*\w[\delta,\alpha]\big)\;,
~~ \delta\in H^{0,0}\;. \myeqn
$$\eqtag{\eDeriv}
If $m_{\alpha}\: H^{0,0}\to H^{0,0}$ is the $\m R$-linear function
on the right of \lasteq0 and $L_{\alpha}\: H^{0,0}\to
H^{0,1}$ is the $\m C$-linear function
$L_{\alpha}(\delta) := [\alpha,\delta]$, it is clear that
the kernel of $L_{\alpha}$ is contained in that of $m_{\alpha}$.
In fact, by direct calculation, for $\delta\in H^{0,0}$
with $\delta=\delta_++\delta_-$ and $(\delta_{\pm})^*=\pm\delta_{\pm}$,
$$
\Re\<\delta,m_{\alpha}(\delta)\> =\<\delta_+,m_{\alpha}(\delta)\>
= -2\nnorm{[\alpha,\delta_+]}^2+\big\<[\delta_-,\delta_+],m(\alpha)\big>\;.
\myeqn
$$
Thus when acting on the self-adjoint elements of $H^{0,0}$,
$m_{\alpha}$ and $L_{\alpha}$ have the same kernel. Since
$m_{\alpha}$ is itself self-adjoint as an $\m R$-linear
map $H^{0,0}\to H^{0,0}$, it follows that $m_{\alpha}$
maps the space of self-adjoint elements of $H^{0,0}$
that are orthogonal to $\ker L_{\alpha}$ isomorphically
into the same space.

At this point, arguments are simplified if it is
assumed that $\alpha$ is stable with respect to the action of
$\Gamma$, not just polystable. Given that $\alpha$ is
$\Gamma$-polystable, the assumption of stability is equivalent
to the condition that $\ker L_{\alpha}=span\,\{1\}$. For
$g\in {\cal G}$, \eAA\ implies $\tr i\wh F(g\cdot (d_0+a))
= \tr i\wh F(d_0+a)+i\,\dbd\log\det(g^*g)
= r\lambda+i\Lambda d\,\tr a+i\,\Lambda\dbd\log\det(g^*\!g)$, so if
$\db_0^*a''=0$  then $\tr i\wh F(g\cdot (d_0+a))
= r\lambda+\lap'\log\det(g^*\!g)$. If $g=\exp(\varphi)\gamma$
for self-adjoint $\varphi\in A^{0,0}(\eet)$ orthogonal to $\ker\db_0$
and $\gamma\in \Gamma$,
$\log\det(g^*\!g)=2\,\tr\varphi+\log\det(\gamma\gamma^*)$, so
$\tr\big(i\wh F(g\cdot (d_0+a))-\lambda\,1\big)=\lap\tr\varphi$.
Hence for $\db_0^*$-closed $a''$ and $g=\exp(\varphi)\gamma$,
by taking the trace-free part of $\varphi$ and rescaling
$\gamma$ to have determinant $1$, it can be assumed
without loss of generality that $\tr i\wh F(d_0+a)
= r\lambda=\tr i\wh F(g\cdot(d_0+a))$, with $\det g\equiv 1
=\det\gamma$.

If $V$ is a hermitian vector space, the real vector space
$End^{\,sa}(V)$
of self-adjoint endomorphisms has a canonically induced
orientation. This can be seen by induction on $\dim V$,
given that the space of self-adjoint endomorphisms of
$V\oplus \m C$ is canonically isomorphic to
$End^{\,sa}(V)\oplus V\oplus \m R$. So the space of self-adjoint
automorphisms of $V$, which is a real closed submanifold
of $Aut(V)$, is orientable. Similarly, the space of trace-free
self-adjoint endomorphisms of $V\oplus\m C$ is canonically isomorphic
to $End^{\,sa}(V)\oplus V$, so the space of self-adjoint automorphisms
of determinant $1$, which is a real closed submanifold of
$Aut(V)$, is orientable.
The function
$\Gamma\owns\gamma\mapsto m(\gamma\alpha\gamma^{-1})$
restricted to the space $\Gamma^{sa}_0$ of
self-adjoint elements of $\Gamma$
of determinant $1$ maps $\Gamma^{sa}_0$ into the space
$(H^{0,0})^{sa}_0$
of trace-free self-adjoint elements of $H^{0,0}$; that is,
into its tangent space at $1\in \Gamma^{sa}_{0}$, and
(given that $\alpha$ is $\Gamma$-stable)
its derivative at
$1$ is an isomorphism between $T_{1}\Gamma^{sa}_0$ and
$(H^{0,0})^{sa}_0$.

Replacing
$\alpha$ in \eEx\
with $t\alpha$ for $t>0$
sufficiently small (depending on $\alpha$), it follows
from that estimate that once $t$ is sufficiently
small, the
kernel of the $\m R$-linear function $(D\wh\Phi)_1$
acting on the trace-free self-adjoint elements of $H^{0,0}$ is
zero, and this holds for any such
$t$ and any $s\beta(t\alpha)$ with
$|s|\le 1$ where $\beta(\alpha)\in A^{0,2}(\eet)$
is the section specified by \rAG.

Now view $\wh\Phi$ in \rED\ as a function
of $\gamma\in \Gamma$ and $\db_0^*$-closed $a''\in A^{0,1}(\eet)$,
taking values in the trace-free self-adjoint elements of $H^{0,0}$.
When restricted to those
$a''\in A^{0,1}(\eet)$ of the form
$a''=t\alpha+s\,\db_0^*\beta_0^{}(t\alpha)$ for
$0<t\le 1$ and $0\le s\le 1$ and self-adjoint $\gamma$
near $1$ of determinant $1$, there
is an induced map $\Gamma^{\rm sa}_0
\to (H^{0,0})^{sa}_0$, for which, from \eEx, the derivative
with respect to its first variable $\gamma$ at $1$ is injective
once $t$ is sufficiently small, where ``sufficiently"
depends upon $\alpha$; more specifically, on the first non-zero
eigenvalue of $L^*_{\alpha}L^{}_{\alpha}$.

\smallskip
\thm{\rClose} Suppose $\alpha\in H^{0,1}$ satisfies
$\Pi^{0,0}\Lambda(\alpha\w\alpha^*+\alpha^*\w\alpha)
=0$. Then for $t>0$ sufficiently small (depending on $\alpha$)
and $s\in [0,1]$
there exists self-adjoint $\delta\in \ker\db_0$ and self-adjoint
$\varphi\in(\ker\db_0)^{\perp}$
(both depending on $s$ and $t$) such
that $i\wh F\big(\exp(\varphi)\cdot \exp(\delta)\cdot
(d_0+a_{s,t})\big)=\lambda\,1$,
where $a_{s,t}^{}=a'_{s,t}+a''_{s,t}$, $a''_{s,t}=t\alpha+
s\db_0^*\beta(t\alpha)$,
and
$\beta(-)\in A^{0,2}(\eet)$ is determined by {\rm\rAG}.
Furthermore, there is a constant $C=C(d_0,\alpha)$ such that
$\norm{\delta}\le Ct^2$ and $\norm{\varphi}_{L^p_2}
\le Ct^2$ for $t>0$ sufficiently small (depending
on $\alpha$).

\pf The assumption on $\alpha$ implies that it is polystable
with respect to the action of $\Gamma$. Assume initially that
$\alpha$ is stable with respect to this action.

If $t>0$ is sufficiently small and $s\in [0,1]$,
$(D_1\wh\Phi)_{(1,a''_{s,t})}$
gives an isomorphism between the tangent space
to $\Gamma^{\rm sa}_0$
at $1$ and $(H^{0,0})^{sa}_0$. The manifold
$\Gamma^{\rm sa}_0$ is orientable,
and the degree of $\wh\Phi(-,a''_{s,t})$ at $\lambda\,1$
is independent of $s$ and sufficiently small $t$, and is therefore
equal to the degree at $\lambda\,1$ for
such $t$ and $s=0$. If $\varphi_t :=
\varphi(a_{0,t})$, then by \eEC, \rFzz\  and the remarks following
that lemma, for $\delta\in H^{0,0}$ near $0$,
there is a function  $R_4=R_4(t,\delta)$ with $\norm{R_4(t,\delta)}_{L^p}
\le Ct^4$ such that, for $\gamma_{\delta} := \exp(\delta)$,
$$
\eqalignno{
\wh\Phi(\gamma_{\delta},a''_{0,t})
&=\Pi^{0,0}i\wh F\big(\exp(\varphi(\gamma_{\delta}\cdot a_{0,t}))
\cdot (d_0+\gamma_{\delta}\cdot a_{0,t})\big) \cr
&= \lambda\,1+\lap_0\varphi(a_{0,t})-i\,t^2\Pi^{0,0}\Lambda(\alpha\w\alpha^*
+\alpha^*\w\alpha)+t^2\,m_{\alpha}(\delta)+R_4(t,\delta)\cr
&=\lambda\,1+t^2\,m_{\alpha}(\delta)+R_4(t,\delta)\;. &\eqn
}
$$
Since both $\wh\Phi-\lambda\,1$ and $m_{\alpha}$ take their values in the
space of trace-free self-adjoint elements of $H^{0,0}\!$, so too does
$R_4$. Since $m_{\alpha}$ is an isomorphism on this space
and $R_4(t,\delta)/t^4$ is uniformly bounded as $t\to 0$, it follows that
$\lambda\,1$ is in the range of $\wh\Phi(-,a''_{0,t})$
for $t$ sufficiently small, and indeed that the degree of
$\wh\Phi(-,a''_{0,t})$ at $\lambda\,1$ is precisely $1$.
Therefore, for $s\in [0,1]$ and $t>0$ sufficiently small
(depending on $\alpha$), there
exists $\gamma\in\Gamma$ such that
$i\wh F\big(\exp(\varphi(\gamma\cdot a_{s,t}))\cdot\gamma\cdot
(d_0+a_{s,t})\big)=\lambda\,1$.

The invertibility of
$m_{\alpha}$ on the space $(H^{0,0})^{sa}_0$ of trace-free self-adjoint elements
of $H^{0,0}$ implies that there is a solution $\gamma=
\exp(\delta(t))$
to $i\wh F\big(\exp(\gamma\cdot a_{s,t}))\cdot \gamma
\cdot (d_0+a_{s,t})\big)=\lambda\,1$ depending continuously
on $t$, and from \lasteq0, the section $\delta(t)$ satisfies
a bound of the form $\norm{\delta(t)}\le Ct^2/c_{\alpha}$
where $c_{\alpha}$ is the lowest eigenvalue of $L_{\alpha}^*L_{\alpha}^{}$
acting on $(H^{0,0})^{sa}_0$
and $C=C(d_0,\alpha)$.
Then  uniform bounds on $\exp(\delta(t))$ give estimates
on $\exp(\delta(t))\cdot\alpha$, and together with the estimates
of \rEA, a uniform $L^p_2$ bound on $\varphi$ of the
form $Ct^2/c_{\alpha}$ follows for some new constant $C$
depending on $d_0$ and $\alpha$.

If now $\alpha$ is assumed to be polystable but not stable,
then the isotropy subgroup $\Gamma_{\alpha}$ has dimension
greater than one. Hence there are non-zero trace-free elements
of $H^{0,0}$ commuting with $\alpha$, and by 2.\ of \rCK,
so too do their adjoints, as they all do with
$\beta(t\alpha)$ for any (small) $t>0$. Hence there are non-zero trace-free
self-adjoint elements of $H^{0,0}$ commuting with $\alpha$
and $\beta(t\alpha)$ for any such $t$, and therefore there are non-zero
trace-free
self-adjoint endomorphisms of $\et$ that are covariantly
constant with respect to $d_0+a_{s,t}$ for all $s,t$. Any
such endomorphism gives a unitary splitting of the bundle
and connections, and these splittings are all compatible
with one another, including the splitting of $(\et,d_0)$.
With respect to such a splitting, the form $\alpha$ splits
into a collection of endomorphism-valued $(0,1)$-forms
on $X$, each of which is $\db$-harmonic with respect to
the induced connection. The ``off-diagonal" components
of $\alpha$ with respect to such a splitting are zero,
since $\alpha$ commutes with the trace-free endomorphisms
determining the splitting.
Since $\alpha$ is of minimal norm in its orbit under
$\gamma$, each of the ``diagonal" components of $\alpha$ must be of
minimal norm under the action of the
subgroup of $\Gamma$ that is the automorphism group of
the corresponding component of $(\et,d_0)$. Hence each of these
components defines an element of the corresponding space that
is polystable with respect to the action of the corresponding
automorphism group, so it follows by induction on the rank
$r$ that for $t$ sufficiently small, for each of these new bundles
with connection, once $t>0$ is sufficiently small there is
a complex gauge transformation that gives a new connection with
$i\wh F$ a scalar multiple of $1$, (the case $r=1$ being
elementary).  Since the splitting of
$(\et,d_0)$ is unitary and $i\wh F(d_0)=\lambda\,1$, the
relevant scalar in all cases is $\lambda$. For each
summand, the estimates on the corresponding endomorphisms
$\delta$ and $\varphi$ imply an estimate of the required
form on the direct sum connection, verifying the last statement of
the theorem.
 \quad \qed

\medskip
A direct proof of \rClose\ that does not use induction
on rank appears to be possible, but raises a number of interesting
representation-theoretic questions.

\smallskip

\bigskip
\cor{\rEint} Under the hypotheses of {\rm \rClose}, if $t>0$ is
sufficiently small and  $d_0+a_{s,t}$
is integrable, then the corresponding holomorphic bundle is
polystable, and is stable only if $\alpha$ is stable with respect
to the action of $\Gamma$.

\pf By \rClose\ and the results of Kobayashi and L\"ubke, the
holomorphic bundle defined by $d_0+a_{s,t}$ is polystable. If
it is not stable, then there exists a non-zero trace-free
holomorphic endomorphism of this bundle. By \rCK, this endomorphism
is covariantly constant with respect to $d_0$ and commutes with
$\alpha$, these facts contradicting the assumption that $\alpha$
is stable with respect to the action of $\Gamma$. \quad \qed

\bigskip
\rems{\kGw}

\proclaim{1. } \advance\hfuzz by 3pt
From the viewpoint of deformation theory, an unobstructed
infinitesimal deformation is (poly)stable with respect to the action
of $\Gamma$ if and only if there is a $1$-complex parameter
family of (poly)stable holomorphic structures whose tangent at
$E_0$ is the given infinitesimal deformation. Of course, in
the case that the latter is polystable but not stable, there
may be families of bundles that are semi-stable but not polystable with that
tangent vector, which will often be the case if $H^2(X,End\,E_0)$
vanishes.
(Direct sums of non-isomorphic line bundles of degree zero
on a torus provide an example when this is not the case.)

\proclaim{2. }
A relatively straightforward application of the continuity
method applied to the assignment $t\mapsto g_t\in {\cal G}$
solving $i\wh F(g_t\cdot(d_0+a_{s,t})=\lambda\,1$ gives
a more quantitative version of \rClose,
namely that the dependence of $t$ on $\alpha$ stated in the
theorem can be made explicit if the constant $C$ there
is replaced by 
by $C_0\norm{\alpha}^2/c_{\alpha}$ where $C_0$ is a constant
depending only on $d_0$ and where $c_{\alpha}$ is the
first non-zero eigenvalue of $L_{\alpha}^*L_{\alpha}^{}\:
H^{0,0}\to H^{0,0}$. In the interests of brevity, an explicit
proof will not be given.
\medskip

The following proposition is the companion uniqueness result to \rClose\
(existence). The proof of the first statement
is based on the proof of Corollary~9 in \cite{Do3}:

\prop{\rEuniq} Suppose $\alpha\in H^{0,1}$,
$a''=\alpha+\db_0^*\beta$, $a=-(a''{})^*+a''$,  and
$\norm{a}_{L^p_1} <\epsilon$ where $\epsilon>0$ is as in {\rm \rCK}.
If $i\wh F(g_1\cdot (d_0+a))=\lambda\,1=i\wh F(g_2\cdot(d_0+a))$,
then $g_2=u\,g_1\gamma$ for some $u\in {\cal U}$ 	and
$\gamma\in \Gamma_{\alpha}$.
Furthermore, there exists
$g_0\in {\cal G}$ with $i\wh F(g_0\cdot(d_0+a))=\lambda\,1$
satisfying the conditions that $g_0^{}=g^*_0$ is positive,
$\det g_0\equiv 1$, and $\Pi^{0,0}(g^*_0g_0^{})\in H^{0,0}$ is
orthogonal to the space of trace-free elements of $\ker L_{\alpha}$,
with these conditions determining $g_0$  uniquely up to conjugation
by unitary elements of $\Gamma_{\alpha}$.

\pf
Let $d_a := d_0+a$ and set $d_b :=
g_1\cdot d_a$, so for $g := g_2^{}g_1^{-1}$ it follows
that $g_2\cdot d_a=g\cdot d_b$ with $i\wh F(d_b)=\lambda\,1
=i\wh F(g\cdot d_b)$. After a unitary change of gauge applied
to $g_2\cdot d_a$, it
can be supposed that $g$ is positive self-adjoint, with
$g=\exp(v)$ for some self-adjoint $v$. If $y\in \m R$
and with $d_y := \exp(y\,v)\cdot d_b$, by \eEaa\
the function
$\m R\owns y\mapsto \<i\wh F(\exp(y\,v)\cdot d_b)-\lambda\,1,v\>\in\m R$ has
derivative $\<\lap_yv,v\>=\norm{d_yv}^2\ge 0$ and is therefore
a non-decreasing function on $\m R$. Since it attains the
value $0$ at both $y=0$ and $y=1$, it must be constant
on $[0,1]$ with derivative identically $0$.
Hence $d_bv=0$, which implies that $\db_a(g_1^{-1}g^{}_2)=0$.
By 1.\ of \rCK, $\gamma := g_1^{-1}g_2^{}$ is $d_0$-covariantly
constant and commutes with $\alpha$.
Thus $ug_2=g_1\gamma$ for some
$u\in {\cal U}$ and $\gamma\in \Gamma_{\alpha}$.

To prove the second statement, note that
$\Gamma$ acts freely on ${\cal G}$ by right multiplication, as does
the closed subgroup $\Sigma_{\alpha}\subset \Gamma$ of elements
in $\Gamma_{\alpha}$ of unit determinant. Given a fixed
$g\in {\cal G}$ there is a constant $c=c(g)$ such
that $c\norm{\gamma}^2 \le \norm{g\gamma}^2\le c^{-1}\norm{\gamma}^2$,
so there exists $\gamma_0\in \Sigma_{\alpha}$ minimising
$\norm{g\gamma}^2$ over all $\gamma\in\Sigma_{\alpha}$.
The Euler-Lagrange equation for this functional on $\Sigma_{\alpha}$
is $\Pi(g_0^*g_0^{})=0$, where $\Pi$ is $L^2$-orthogonal
projection onto the space of trace-free elements in $\ker L_{\alpha}$,
the Lie algebra of $\Sigma_{\alpha}$.

Suppose now that $g_1, \, g_2\in {\cal G}$ are as in the
statement of the proposition, with $g_2=ug_1\gamma$ for some
$u\in{\cal U}$ and some $\gamma\in \Gamma_{\alpha}$.
Suppose in addition that both $g_1$ and $g_2$ have
unit determinant, are both positive and self-adjoint, and
that $\Pi^{0,0}(g_j^*g_j^{})$ is orthogonal to the
trace-free elements of $\ker L_{\alpha}$ for $j=1,2$. Then for
any trace-free $\phi\in\ker L_{\alpha}$, and using the fact
that the trace of a covariantly constant endomorphism is
constant,
$$
0=\<g_2^*g_2^{},\phi\>
=\<\gamma^*g_1^*g_1^{}\gamma,\phi\>
=\<g_1^*g_1^{},\gamma\phi\gamma^*\>
=\<g_1^*g_1^{},{\tr\gamma\phi\gamma^*\over r}\,1\>
=\norm{g_1}_{L^2}^2\,{\tr\gamma\phi\gamma^*\over r}
= \norm{g_1}^2_{L^2}{\<\gamma^*\gamma,\phi\>\over r}\;.
$$
Therefore $\gamma^*\gamma$ is a multiple of $1$, and this
multiple must be $1$ since $1=\det u\,\det\gamma$. Thus
$\gamma\in U(\Gamma_{\alpha})$, the group of unitary
elements in $\Gamma$ commuting with $\alpha$. Then since
$g_1$ and $g_2$ are both self-adjoint,
$g_2^2=g_2^*g_2^{}=\gamma^*g_1^*g_1^{}\gamma=
(\gamma^{-1}g_1\gamma)^2$, and positivity implies
$g_2=\gamma^{-1} g_1 \gamma$. From $g_2=ug_1\gamma$,
it then follows that $u=\gamma^{-1}$. \quad \qed

\medskip
\rem{\kGbb} In general, if $g\in {\cal G}$ satisfies
$i\wh F(g\cdot d_a)=\lambda\,1$, the particular form of $a
=a'+a''$ for $a''=\alpha+\db_0^*\beta$ imposes certain
conditions on $h := g^*g$, at least if $\alpha$ is sufficiently
small. Specifically, if $s\in A^{0,0}(\eet)$ is
$\db_a$-closed, then by \rCK\ so too is $s^*$ and both
commute with $\alpha$ and $\beta$. Then the fact that
$g^{-1}sg$ is $\db_b$-closed for $d_b := g\cdot d_a$
together with the condition that $i\wh F(d_b)=\lambda\,1$
implies that $h^{-1}s^*h$ is $\db_a$-closed and indeed,
$d_a$-closed. The same applies to $s^*$, so
$H^{0,0}\owns s\mapsto h^{-1}sh\in H^{0,0}$ is an
automorphism of $H^{0,0}$, one that is self-adjoint and
positive. Any eigenvector with eigenvalue that is not $1$
must be a $d_0$-closed nilpotent endomorphism of $\et$ that
commutes with $\alpha$ and $\beta$ and their adjoints.
The
bundle-with-connection $(\et,d_a)$ has a unitary splitting
into bundles-with-connection that are irreducible with
respect to the action of $\Gamma_{\alpha}$, these splittings
being compatible with splittings of $(\et,d_0)$ into
a direct sum of polystable components. If $\Pi$
is
orthogonal projection onto the direct sum of all
components isomorphic to a given component, then
$h^{-1}\Pi h$ must map that sum to itself, but
is not necessarily self-adjoint.  In sum, this structure raises
a number of interesting questions associated with the
representation theory of $\Gamma_{\alpha}$.

\medskip

\rClose\ gives a condition under which a
connection near $d_0$ has a
connection with central component of the curvature equal to a
scalar multiple of the identity in its orbit under ${\cal G}$,
but the deficiency
of the result is that how near to $d_0$ the connection must be
depends on the connection itself,  dictated by the relative
sizes of the eigenvalues of $L_{\alpha}^*L_{\alpha}^{}$.
This issue is addressed in the next section.

\bigskip
\secn{The local Hitchin-Kobayashi correspondence.}

As stated at the end of the previous section, the objective of
this section is to remove the dependency of \rClose\ on $\alpha$
other than through $\norm{\alpha}$. That is, retaining all of
the notion of that section, the objective is to prove the following
result:

\medskip

\thm{\rThere} Let $d_0$ be a connection on $\et$ with
$i\wh F(d_0)=\lambda\,1$. Then there is a constant
$\epsilon=\epsilon(d_0)$ with the following property:
If $\alpha\in H^{0,1}$
is polystable with respect to the action of $\Gamma$
and $\norm{\alpha}<\epsilon$, and if
$\beta\in A^{0,2}(\eet)$
is as in
{\rm \rAG},
then there exists $g\in {\cal G}$ with
$i\wh F(g\cdot (d_0+a))=\lambda\,1$, where
$a=a'+a''$ for $a''=\alpha+\db_0^*\beta$.

\medskip
The approach to proving this result is to ensure that the analysis
is performed in a sufficiently small neighbourhood of $d_0$ where
the connections are well-approximated by their linearizations,
which has the effect of reducing the problem to a finite-dimensional
question that is naturally attacked using the methods of
classical geometric invariant theory.
Before commencing the proof of the theorem, there are several
remarks and observations that simplify matters considerably.

\smallskip
First, consider the case in which the rank $r$ of the bundle $\et$ is $1$.
Then $\alpha$ is a harmonic $(0,1)$-form on $X$, $\beta$ must be
zero since $\alpha\w\alpha=0$, and the connection
$d_a=d_0+(\alpha-\alpha^*)$ has curvature $F(d_0)$, which already
satisfies the condition $i\wh F=\lambda$. Thus $g\equiv 1$ solves
the equation. In the general case, if $i\wh F(g\cdot(d_0+a))
=\lambda\,1$, then on taking the trace of both sides it
follows that
$i\Lambda\big(\tr(F(d_0)+d_0a+a\w a)+\ddb\log\det(g^*\!g))=r\,\lambda$,
which implies that $i\Lambda\,\ddb\log\det(g^*\!g)\equiv 0$ and hence that
$\det(g^*\!g)$ is constant. After rescaling $g$ by a constant, it
can therefore be assumed that $|\det g|\equiv 1$.

Second, given that $\alpha$ is polystable with respect to the action
of $\Gamma$, it may be assumed without loss of generality that
$\alpha$ is of minimal norm in its orbit under $\Gamma$, and
therefore $i\Lambda(\alpha\w\alpha^*+\alpha^*\w\alpha)$ is
orthogonal to $\ker\db_0$, by \rFA.

Third, as it was for the proof
of \rClose,  if $\alpha$ is polystable but not
stable, precisely the same argument using induction on $r$
that was employed at the
end of the proof of \rClose\ reduces the problem to
the case when $\alpha$ is stable with respect to the
action of $\Gamma$. Given this, the uniqueness result
\rEuniq\ implies that the only freedom in choice of $g$
is that of replacing $g$ by $ug$ for $u\in {\cal U}$.

\smallskip
Hitherto, little use has been made of unitary gauge
freedom ${\cal U}\owns u\mapsto u\cdot d$ for a connection
$d$, as this is subsumed into the complex gauge freedom
${\cal G}\owns g\mapsto g\cdot d$. But since the
equation $i\wh F(d)=\lambda\,1$ is invariant under
unitary gauge transformations of the connection $d$, it is helpful to make use
of the opportunity to place connections in good (unitary)
gauges:

\medskip
\lem{\rFgauge} There are constants $\epsilon>0$ and $C$ depending
only on $d_0$ with the property that if $d_0+b$ is a connection
with $\norm{b}_{L^p_1}<\epsilon$ then there is a unique
skew-adjoint section $\psi\in A^{0,0}(\eet)$ orthogonal to
$\ker d_0$ for which
$d_0^*\big(e^{\psi}\cdot (d_0+b)-d_0\big)=0$, with
$\norm{\psi}_{L^p_2} \le C\norm{d_0^*b}_{L^p}$.

\pf The linearization of the function
${\cal U}\owns u\mapsto d_0^*\big(ubu^{-1}-d_0u\,u^{-1}\big)$
at $u=1$ and $b=0$ is $A^{0,0}(\eet)\owns \sigma
\mapsto \lap_0\sigma$, which is an isomorphism from the
space of  skew-adjoint elements of $A^{0,0}(\eet)$ orthogonal
to $\ker d_0$ lying in $L^p_2$ to the same such space of elements
lying in $L^p$. Since the original function takes
values in the latter space, an application of the
implicit function theorem implies that there is
a number $\epsilon>0$ such that the equation
$d_0^*\big(e^{\psi}be^{-\psi}-(d_0e^{\psi})e^{-\psi}\big)
=0$ has a unique skew-adjoint solution $\psi\in (\ker d_0)^{\perp}
\subsett A^{0,0}(\eet)$
if $\norm{b}_{L^p_1}<\epsilon$, and moreover
$\norm{\psi}_{L^p_2} \le C\norm{d_0^*b}_{L^p}$ for some
$C=C(d_0)$. \quad\qed

\medskip

Turning now to the proof of \rThere\ and retaining most
of the notation of the previous section,
suppose that  $\alpha\in H^{0,1}$ is stable with respect to
the action of $\Gamma$ and is of minimal norm
in its orbit under this action, with  $\norm{\alpha}=1$.
For $t >0$ sufficiently small that \rAG\ is valid, let
$\beta_t := \beta(t\alpha)$ and let $a_t=a_t'+a_t''$ for
$a_t'':= t\alpha+\db_0^*\beta_t$, with $d_t := d_0+a_t$.
Fix a number $\epsilon_0\in (0,1]$, the precise value of which
will be fixed later, but for the moment should satisfy
the condition that for any $t\in (0,\epsilon_0]$, $a_t$ satisfies
the hypotheses of \rAE, \rAG, \rCK, \rCL\ and \rFgauge. Now let
$$
\eqalign{
S:= \Big\{t_0\in (0,\epsilon_0]\;\;\vrule height 10pt depth 5pt width .2pt\;\;
   &\hbox{\sl for every $t\in (0,t_0]$ there exists $g\in {\cal G}$
   with $\;\nnorm{(\db_{t}g)g^{-1}}_{L^p_1}<t$} \cr
\noalign{\vskip-6pt}
   &\hbox{\sl \quad for which $g\cdot d_t =: d_0+b_t$ satisfies
	$d_0^*b_t^{}=0\;$ and  $\;i\wh F(d_0+b_t)=\lambda\,1$.}\Big\}
}
$$
By \rClose, for $t>0$ sufficiently small (depending on $\alpha$)
there exist trace-free self-adjoint
$\delta\in \ker\db_0$ and $\varphi\in (\ker\db_0)^{\perp}$
with $\norm{\delta}^2+\norm{\varphi}_{L^p_2}\le C_{\alpha}t^2$
such that $i\wh F(g\cdot d_t)=\lambda\,1$ for $g=\exp(\varphi)
\exp(\delta)$. Then $(\db_tg)g^{-1}
=(\db_0g)g^{-1}+a_t''-ga_t''g^{-1}
=(\db_0e^{\varphi})e^{-\varphi}+[a_t'',g]g^{-1}$. The first
term is bounded in $L^p_1$ by $C_{\alpha}t^2$, and
since $g^{-1}$ is uniformly bounded in $C^0$ whilst
$\norm{g-1}_{L^p_1} \le C_{\alpha}t$, the bound $\norm{a''_t}
\le Ct$ implies that
$\norm{(\db_tg)g^{-1}}_{L^p_1}\le C_{\alpha}t^2$ for some new
constant $C_{\alpha}$. Since
$\db_0^*a''_t=0$, it follows easily that
$
\nnorm{\db_0^*\big(-(\db_0g)g^{-1}+ga''_tg^{-1}\big)}_{L^p}
\le C_{\alpha}t^2
$, so by \rFgauge, after a unitary gauge transformation
$u\cdot g\cdot d_t = d_0+b_t$ so that $d_0^*b_t=0$,
the complex automorphism $\tilde g =ue^{\varphi}e^{\delta}\in{\cal G}$
satisfies the requirements for $t$ to lie in $S$ once
$t>0$ is sufficiently small. Thus $S$ is not empty.

The fact that $S$ is open (if $\epsilon_0$ is sufficiently
small) will be shown shortly, this being a straightforward
consequence of the implicit function theorem. The proof that
$S$ is closed is more involved, involving a priori estimates
on solutions.

\smallskip

To see that $S$ is open, suppose that
$t_0\in (0,\epsilon_0)\cap S$, and let $g_0\in{\cal G}$
satisfy $\nnorm{(\db_{t_0}g_0^{})g_0^{-1}}_{L^p_1}<t_0$,
$d_0^*b_0^{}=0$ for $d_0+b_0 := g_0\cdot d_{t_0}$ and
$i\wh F(d_{t_0})=\lambda\,1$.
The linearization
of the function ${\cal G}\owns g\mapsto i\wh F(g\cdot d_{t_0})$ at
$g_0\in {\cal G}$ is $A^{0,0}(\eet)\owns \sigma
\mapsto d_{b_0}^*d_{b_0}^{}\sigma_+ \in A^{0,0}(\eet)$. If $\sigma$ is in the
kernel of this map, then $\db_{b_0}\sigma_{{\!}_+\,}=0$, so
$(\db_0+ a_{t_0}'')(g_0^{-1}\sigma_{{\!}_+\,}g_0^{})=0$. Given that
$\norm{a_{t_0}}_{L^p_1} \le Ct_0$ and $t_0$ is sufficiently small,
it follows
from \rCK\ that the
endomorphism $g_0^{-1}\sigma_{{\!}_+\,}g_0^{}$
is covariantly
constant with respect to $d_0$ and commutes with
$\alpha$ as well as with $\beta_{t_0}$.
Since $\alpha$
is $\Gamma$-stable, this implies that $g_0^{-1}\sigma_{{\!}_+\,}g_0^{}$
is a scalar multiple of the identity, and therefore
so too is $\sigma_+$. Hence the kernel of $\lap_{b_0}=d_{b_0}^*d_{b_0}^{}$ 
acting on
the trace-free self-adjoint sections of $A^{0,0}(\eet)$ is
zero, and so an application of the implicit function theorem
implies that there is a small neighbourhood of $t_0$ in $(0,\epsilon_0)$
that lies in $S$, proving that $S$ is open.

\smallskip

It remains to show that $S$ is also closed. Suppose now that
$(0,t_0) \subset S$, and for each $t\in (0,t_0)$, let $g_t\in {\cal G}$
satisfy $\det g_t\equiv 1$, $i\wh F(g_t\cdot d_t)=\lambda\,1$ and 
$d_0^*b_t^{}=0$
for $g_t\cdot d_t=:d_0+b_t$, with
$\norm{(\db_t^{}g_t^{})g_t^{-1}}_{L^p_1} < t$.
Since $b_t''=-(\db_0^{}g_t^{})g_t^{-1}+
g_t^{}a_t''g_t^{-1} = -(\db_t^{}g_t^{})g_t^{-1}+a_t''$, it
follows that $\norm{b_t}_{L^p_1} \le Ct$ for some constant $C=C(d_0)$,
and therefore the preparatory results \rAE, \rAG, \rCK, and \rFgauge\ apply
to the connection $d_0+b_t$ if $\epsilon_0$ is sufficiently
small.

The equation $i\wh F(d_0+b_t)=\lambda\,1=i\wh F(d_0)$ implies that
$i\Lambda(d_0b_t+b_t\w b_t)=0$, which can be re-written as
$\d_0^*b_t'-\db_0^*b_t''=-i\Lambda(b_t'\w b_t''+b_t''\w b_t')$.
Since $0=d_0^*b_t=\d_0^*b_t'+\db_0^*b_t''$, it follows that
$2\db_0^*b_t''=i\Lambda(b_t'\w b_t''+b_t''\w b_t')$, implying
that $\norm{\db_0^*b_t''}_{L^p} \le Ct^2$ for some $C=C(d_0)$.
Note that since $\db_0^*a''_t=0$, this is a bound on the
$L^p$ norm of $\db_0^*\big((\db_tg_t)g_t^{-1}\big)$, and because of
the uniform bounds on $a_t$, this can also be seen as a uniform
bound on the $L^p$ norm of $\db_t^*\big((\db_t^{}g_t^{})g_t^{-1}\big)$.

By \rAE, there is a unique self-adjoint
$\varphi_t \in A^{0,0}(\eet)$ orthogonal
to $\ker\db_0$ such that
$$
\db_0^*\big(e^{-\varphi_t}\db_0e^{\varphi_t}+e^{-\varphi_t}b_t''e^{\varphi_t}\big)=0
\quad {\rm with} \quad \norm{\varphi_t}_{L^p_2}\le 
C\norm{\db_0^*b_t''}_{L^p}\le Ct^2\;. \myeqn
$$\eqtag{\eHa}
If $d_0+c_t := \exp(-\varphi_t)\cdot (d_0+b_t)
= (\exp(-\varphi_t)g_t)\cdot(d_0+a_t)$, then $\db_0^*c_t''=0
=\db_0^*a_t''$, so by the first statement of \rCL\
(using the connection on $\eet=\rm Hom(\et,\et)$
induced by $d_0+a_t$ and $d_0+c_t$ respectively),
$e^{-\varphi_t}g_t =: \gamma_t$
is $d_0$-covariantly constant and
$c_t''\gamma_t=\gamma_ta_t''$. Note that since
$\varphi_0$ must be trace-free (being orthogonal to $1$)
and $\det g_t =1$, it follows that $\det\gamma_t=1$.
Since $\gamma_t$ is invertible,
$c_t'' =t\,\gamma_t^{}
\alpha\gamma_t^{-1}+\gamma_t^{}\db_0^*\beta_t\gamma_t^{-1}$,
but since $t\le\epsilon_0$ is sufficiently small,
$\gamma_t^{}\db_0^*\beta_t\gamma_t^{-1}
=\db_0^*\beta(t\gamma_t^{}\alpha\gamma_t^{-1})$.
Since $\db_0+c_t''=\exp(-\varphi_t)\cdot(\db_0+b_t'')$
and $\norm{\varphi_t}_{L^p_2}\le Ct^2$ whilst
$\norm{b_t''}_{L^p_1}\le Ct$, it follows that $\norm{c_t''}_{L^p_1}
\le Ct$ and in particular, $\norm{c_t''}_{L^2}\le Ct$.
But since
$$
\norm{c_t''}^2_{L^2} =t^2\nnorm{\gamma_t^{}\alpha\gamma_t^{-1}}_{L^2}^2
+\nnorm{\db_0^*\beta(t\gamma_t^{}\alpha\gamma_t^{-1})}_{L^2}^2\;,
\myeqn
$$\eqtag{\eHb}
it follows that $\alpha_t := \gamma_t^{}\alpha\gamma_t^{-1}$ is uniformly
bounded in $L^2$ independently of $t$ and $\alpha$, and since $\alpha$
is stable with respect to the action of $\Gamma$,
the assignment
$\Gamma\owns \gamma\mapsto \norm{\gamma\alpha\gamma^{-1}}$
is proper on the elements $\gamma \in \Gamma$ that have unit determinant,
and therefore $\gamma_t$ converges to some $\gamma_0\in\Gamma$
as $t\to t_0$. Since $t_0$ is sufficiently small that
\rAG\ applies, $\gamma_t^{}\beta_t^{}\gamma_t^{-1}
=\beta(t\gamma_t^{}\alpha\gamma_t^{-1})$ converges in
$L^p_2$ to $\beta(t_0\gamma_0^{}\alpha\gamma_0^{-1})$ as
$t\to t_0$.

The uniform $L^p_2$ bounds on $\varphi_t$ imply that these
converge weakly in $L^p_2$ and
by the Sobolev embedding theorem, strongly in $C^1$ to
$\varphi_0\in (\ker\db_0)^{\perp}$, so the corresponding
automorphisms $g_t=\exp(\varphi_t)\gamma_t$ have the same
convergence.
By ellipticity of the $\db_0$-Laplacians on $A^{0,1}(\eet)$
and $A^{0,2}(\eet)$, the function $\beta$ of \rAG\ depends
smoothly on its argument, and therefore
$\beta(t\gamma_t^{}\alpha\gamma_t^{-1})$ converges smoothly
to $\beta(t_0\gamma_0^{}\alpha\gamma_0^{-1})$ as $t\to t_0$.
Ellipticity of the equations $\d_t^*(h_t^{-1}\d_t^{}h_t^{})
= \lambda\,1-h_t^{-1}F(d_t^{})h_t^{}$ for $h_t := g_t^*g_t^{}$
and smooth dependence of $a_t$ on $t$ then imply that
the family $(h_t)$ depends smoothly on $t$, and
the uniqueness of the unitary gauge stated in \rFgauge\
then gives smooth dependence of $(g_t)$ on $t$, so
$g_t$ converges smoothly to some $g_0\in {\cal G}$
as $t\to t_0$, with $d_0+b_0 := g_0\cdot(d_0+a_{t_0})$
satisfying $i\wh F(d_0+b_0)=\lambda\,1$,
$d_0^*b_0^{}=0$, and $\nnorm{(\db_{t_0}g_0^{})g_0^{-1}}_{L^p_1}
\le t_0$. The proof that $S$ is closed will be complete
if it can be shown that this is a strict inequality, but ensuring this
is the critical issue.

\smallskip
As above, $\alpha_t := \gamma_t\cdot \alpha = \gamma_t^{}\alpha\gamma_t^{-1}$
is uniformly bounded in $L^2$ and hence in $L^p_2$
independent of $t\in (0,\epsilon_0]$
and of $\alpha$. Since $c_t''=t\alpha_t +\db_0^*\beta(t\alpha_t)$,
$\norm{c_t''-t\alpha_t}_{L^2} \le C\norm{\beta(t\alpha_t)}_{L^p_1}
\le C\norm{t\alpha_t}_{L^2}^2\le Ct^2$ for some uniform
constant $C$, using here \rAG. Again using the $L^p_2$ bounds
on $\varphi_t$, it then follows that
$$
t^2\nnorm{i\Lambda\big(\alpha_t^{}\w\alpha_t^*+\alpha_t^*\w\alpha_t\big)}
\le \nnorm{i\Lambda\big(b''_t\w b'_t+b'_t\w b''_t\big)}
+ Ct^3 = Ct^3\;. \myeqn
$$\eqtag{\eHd}
Thus, if
$
m(\sigma) := \Pi^{0,0}i\Lambda\big(\sigma\w\sigma^*+\sigma^*\w\sigma\big)$
for $\sigma\in H^{0,1}$
is the moment map for the action of $\Gamma$ as in \S5,
then $\norm{m(\alpha_t)} \le Ct$ for some
constant $C$ that is independent of $\alpha$ and $t$. In fact,
it follows from \eEaz(a) that there is a better estimate,
namely $\norm{m(\alpha_t)}\le Ct^2$ for some new constant
$C$ independent of $\alpha$ and $t$, provided that $\epsilon_0$
is sufficiently small.

The relationship between
$\norm{\sigma}$ and $\norm{m(\sigma)}$ for $\sigma\in H^{0,1}$
(or more precisely, the behaviour of the
projectively-invariant function $\norm{m(\sigma)}/\norm{\sigma}^2$)
has been studied by a number of experts in the field,
lying at the heart of the interaction between geometric
invariant theory and symplectic geometry. It
plays a critical role in determining the
topology of quotients; see, for example, \cite{Kir}, \cite{Nee}.

\medskip
\lem{\rHc} For every $\epsilon>0$ there is a $\delta>0$ such
that $\norm{\gamma\alpha\gamma^{-1}}^2-1 <
\epsilon$ for every polystable $\alpha\in H^{0,1}$
with $m(\alpha)=0$ and $\norm{\alpha}=1$
for which $\norm{m(\gamma\alpha\gamma^{-1})} < \delta$.

\medskip
Assuming for the moment that this ``uniform continuity" result
holds, the proof of \rThere\ is easily completed. As before,
$$
\db_tg_t^{}g_t^{-1} =\db_0g_t^{}g_t^{-1}+a''_t-g_t^{}a_t''g_t^{-1}
= (\db_0e^{\varphi_t})e^{-\varphi_t}+a_t''-
e^{\varphi_t}\gamma_t^{}a_t''\gamma_t^{-1}e^{-\varphi_t}\;,
$$
so with the earlier estimates on $\norm{\varphi_t}_{L^p_2}$
and on $\norm{\beta}_{L^p_2}$,
$$
\norm{\db_tg_t^{}g_t^{-1}}_{L^p_1} \le
t\norm{\alpha_t-\alpha}_{L^p_1} + Ct^2
\le C't\norm{\alpha_t-\alpha}_{L^2} + Ct^2 \myeqn
$$\eqtag{\eHz}
for some uniform constants $C$ and $C'$.

The unitary gauge transformations provided by \rFgauge\ have
been applied to the connections with central curvature $-i\lambda\,1$,
but as yet, none has been applied to the connections $d_t$.
This is now done by writing $\gamma_t = p_tu_t$ for some uniquely
determined positive self-adjoint $p_t\in \Gamma$ and some
unitary $u_t\in \Gamma$. The convergence of $\gamma_t$ to
$\gamma_0$ implies convergence of $p_t$ and $u_t$ to some
positive $p_0$ and $u_0$ in $\Gamma$ respectively. Note that
all the estimates above apply equally with no change of constants
when $\alpha$ is replaced by $\alpha_t$ since they depended
on $\alpha$ only through $\norm{\alpha}$.
Since $\alpha_t' := u_t^{}\alpha u_t^{-1}$ is of minimal norm
in its orbit under $\Gamma$ (i.e., that of $\alpha$), it follows that
$$
\nnorm{\alpha_t-\alpha_t'}^2=
\norm{\alpha_t}^2-2\,\Re\big<\alpha_t,\alpha_t'\big\>+\norm{\alpha_t'}^2
=\norm{\alpha_t}^2-2\,\big\<p_t^{1/2}\cdot\alpha_t'\,,\,p_t^{1/2}\cdot\alpha_t'\big\>
+\norm{\alpha_t'}^2
\le \norm{\alpha_t}^2-1\;.
$$
From \lasteq0, $\norm{\db_tg_t^{}g_t^{-1}}_{L^p_1}$ will
be less than $t/2$ if both $t\le t_0$ and $t_0$ is
sufficiently small (depending only on $d_0$) and
$\norm{\alpha_t-\alpha_t'}$ is also sufficiently small. The
latter condition will hold if $\norm{\alpha_t}^2$ is sufficiently
close to $1$, and by \rHc, this in turn will hold if
$\norm{m(\alpha_t)}$ is sufficiently small. From \eHd,
that last condition will be satisfied provided that
$t$ is sufficiently small, where ``sufficiently small" is a condition
that depends only on $d_0$, and not on $\alpha$. Consequently,
provided that $\epsilon_0$ is chosen to be sufficiently small,
the set $S$ will be closed as well as open, and hence be
equal to $(0,\epsilon_0]$, completing the proof of \rThere. \quad\qed

\medskip
It remains to prove \rHc, which will be a consequence of the following:

\medskip
\lem{\rHskew} Suppose $\alpha\in H^{0,1}\-\{0\}$ satisfies
$m(\alpha)=0$. Then for any $\gamma\in\Gamma$,
$\nnorm{\gamma\alpha\gamma^{-1}}^2-\norm{\alpha}^2
\le C\nnorm{m(\gamma\alpha\gamma^{-1})}$ for some
constant $C$ depending only on $d_0$.

\pf
The polystable holomorphic bundle $E_0$ splits as a
direct sum $\bigoplus_{i=1}^m E_i$ of stable bundles all of the
same slope.
With respect to this splitting of $E_0$, a
form $\tau\in H^{0,1}$ corresponds to an $m\times m$ matrix
$[\tau_{j}^{\;i}]$ of $(0,1)$-forms, with
$\tau_j^{\;i}$ being $\db$-harmonic with respect to the
induced Hermite-Einstein connection on $Hom(E_j,E_i)$.
Then
$$
\norm{\tau}^2
= \sum_{i=1}^m\sum_{j=1}^m\norm{\tau_j^{\;i}}^2\;.
$$
Moreover,
$m(\tau)=\Pi^{0,0}i\Lambda(\tau\w\tau^*+\tau^*\w\tau)$
corresponds to an $m\times m$ matrix for which the
$i$-th diagonal entry is
$$
m(\tau)_i^{\;i}=\sqrt{-1}\sum_{j=1}^m\Pi^{0,0}_{\;i}
\Lambda\big((\tau_j^{\;i})^*\w\tau_j^{\;i}
+\tau_i^{\;j}\w (\tau_i^{\;j})^*\big)\;,
$$
where $\Pi^{0,0}_{\;i}$ is $L^2$-orthogonal projection onto
the $(i,i)$-component (of $\ker\db_0$). Since
$E_i$ is stable for each $i$, $Aut(E_i)=\m C^*$ and
so the projection $\Pi^{0,0}_{\;i}$ here is simply
given by integrating the trace over $X$. Thus
$$
m(\tau)_i^{\;i}=\sum^m_{j=1}\big(\norm{\tau_j^{\;i}}^2
-\norm{\tau_i^{\;j}}^2\big)\;.
$$
(More precisely, $m(\tau)_{i}^{\;i}$ is the number on
the right multiplied by the identity endomorphism of
$E_i$, but this fact only changes estimates by
combinatorial factors bounded by a combinatorial function
of $r$.)

Suppose now that $\alpha\in H^{0,1}$ satisfies $m(\alpha)=0$
and $\gamma\in\Gamma$.
Using a Cartan decomposition of $\Gamma$ into
$\Gamma= U\,T\,U$ where $T$ is a
maximal complex torus and $U=U(\Gamma)$, it follows from the
left and right unitary
invariance of the norm and the unitary equivariance of
$m$ that $\gamma$ may be assumed to lie
in $T$; that is, $\gamma=diag(t_1,\dots, t_m)$ for some
$t_j\in \m C^*$.

Instead of working directly with
$\alpha$, it is more convenient to work with $\tau := \alpha\gamma^{-1}$,
so $m(\tau\gamma)=0$ and it must be shown that
$\norm{\gamma\tau}^2-\norm{\tau\gamma}^2\le C\norm{m(\gamma\tau)}$.
This will follow if it can be shown that
$$
\sum_{j=1}^m\big(|t_j|^2\norm{\tau_j^{\;i}}^2
-|t_i|^2\norm{\tau_i^{\;j}}^2) = 0 \quad
\hbox{\sl for $i=1,\dots,m$ implies that} \myeqn
$$
$$
\sum_{i=1}^m\sum_{j=1}^m \big(|t_i|^2\norm{\tau_j^{\;i}}^2
-|t_j|^2\norm{\tau_j^{\;i}}^2\big)
\le C\,\sum_{i=1}^m\bigg|\sum_{j=1}^m\big(|t_i|^2\norm{\tau_j^{\;i}}^2
-|t_j|^2\norm{\tau_i^{\;j}}^2\big)\bigg|
\myeqn
$$
{\sl for some constant $C$\/},
using here the fact that the $\ell_1$ and $\ell_2$ norms
on $H^{0,0}$ are equivalent in this representation.

Observe that
$$
\big(|t_i|^2-|t_j|^2\big)\big(\norm{\tau_{j}^{\;i}}^2-
\norm{\tau_i^{\;j}}^2\big) =
\big(|t_i|^2\norm{\tau_{j}^{\;i}}^2-
|t_j|^2\norm{\tau_j^{\;i}}^2\big) +\big(|t_i|^2\norm{\tau_i^{\;j}}^2
-|t_j|^2\norm{\tau_j^{\;i}}^2\big)\;,
$$
and by \lasteq1, for each fixed $i$ the second term
on the right sums to zero on application of $\sum_j$. Similarly,
$$
\big(|t_i|^2-|t_j|^2\big)\big(\norm{\tau_{j}^{\;i}}^2+
\norm{\tau_i^{\;j}}^2\big)
= \big(|t_i|^2\norm{\tau_j^{\;i}}^2-|t_j|^2\norm{\tau_i^{\;j}}^2\big)
+\big(|t_i|^2\norm{\tau_{i}^{\;j}}^2-|t_j|^2\norm{\tau_j^{\;i}}^2\big)\;,
$$
and again \lasteq1 implies that the second term on the right
sums to zero on application of $\sum_j$. So \lasteq0 is equivalent to
$$
\sum_{i=1}^r\sum_{j=1}^r
\big(|t_i|^2-|t_j|^2\big)\big(\norm{\tau_{j}^{\;i}}^2-
\norm{\tau_i^{\;j}}^2\big)
\le C\,\sum_{i=1}^r\bigg|\sum_{j=1}^r
\big(|t_i|^2-|t_j|^2\big)\big(\norm{\tau_{j}^{\;i}}^2+
\norm{\tau_i^{\;j}}^2\big)\bigg|\;.
$$
After renumbering, it can be assumed that $|t_i|\ge |t_j|$
if $i<j$.
The summand on the left is symmetric under interchange
of $i$ and $j$, so it can be written as
$2\sum_{i=1}^r\sum_{j=i+1}^r
\big(|t_i|^2-|t_j|^2\big)\big(\norm{\tau_{j}^{\;i}}^2-
\norm{\tau_i^{\;j}}^2\big)$. Then the desired inequality
will certainly follow if it can be shown that
$$
\sum_{i=1}^r\sum_{j=i+1}^r
\big||t_i|^2-|t_j|^2\big|\big(\norm{\tau_{j}^{\;i}}^2+
\norm{\tau_i^{\;j}}^2\big)
\le C\,\sum_{i=1}^r\bigg|\sum_{j=1}^r
\big(|t_i|^2-|t_j|^2\big)\big(\norm{\tau_{j}^{\;i}}^2+
\norm{\tau_i^{\;j}}^2\big)\bigg|\;.
$$
That this is true is a consequence of the following:

\medskip
\lem{\rHatlast}
Let $S = (s_{ij})$ be a skew-symmetric $m\times m$ matrix with
$s_{ij}\ge 0$ if $i<j$. Then
$$
\sum_{i=1}^m \sum_{j=i+1}^ms_{ij} \le
2^{m-1}\sum_{i=1}^m\bigg|\sum_{j=1}^m s_{ij}\bigg|\;.
$$

\pf It will be shown inductively that for $k\in \{1,\dots, m\}$,
$$
\sum_{i=1}^k\sum_{j=i+1}^ms_{ij}
\le 2^{k-1}\sum_{i=1}^k\bigg|\sum_{j=1}^m s_{ij}\bigg|\;. \myeqn
$$
For $k=1$, the inequality clearly holds since $s_{1j}\ge 0$ for
every $j$. Suppose that
the inequality has been shown to hold for $k=1,\dots,\ell-1$. Then for
$k=\ell$, the new term on the left is
$L_{\ell} := \sum_{j=\ell+1}^m s_{\ell j}$, and
on the right the new term is $R_{\ell} := \big|\sum_{j=1}^m s_{\ell j}\big|$.
Let $A_{\ell} := -\sum_{j=1}^{\ell-1} s_{\ell j}\ge 0$ and
$B_{\ell} := \sum_{j=\ell+1}^m s_{\ell j}\ge 0$, so $R_{\ell} =
\big|B_{\ell}-A_{\ell}\big|$ and  by inspection, $B_{\ell}=L_{\ell}$.

By skew-symmetry of $S$,
$A_{\ell}=\sum_{j=1}^{\ell-1}s_{j\ell}=\sum_{i=1}^{\ell-1}s_{i\ell}
\le \sum_{i=1}^{\ell-1}\sum_{j=\ell}^ms_{ij}
\le 2^{\ell-2}\sum_{i=1}^{\ell-1}\big|\sum_{j=1}^m s_{ij}\big|$,
using the inductive hypotheses.
Now if $B_{\ell}-A_{\ell}\le 0$, then
$
L_{\ell} = B_{\ell}\le A_{\ell} \le
2^{\ell-2}\sum_{i=1}^{\ell-1}\big|\sum_{j=1}^m s_{ij}\big|\;,
$
so
$$
\sum_{i=1}^{\ell}\sum_{j=i+1}^ms_{ij}
= B_{\ell}+\sum_{i=1}^{\ell-1}\sum_{j=i+1}^m s_{ij}
\le A_{\ell}+\sum_{i=1}^{\ell-1}\sum_{j=i+1}^m s_{ij}
\le 2^{\ell-1}\sum_{i=1}^{\ell-1}\bigg|\sum_{j=1}^m s_{ij}\bigg|
\le 2^{\ell-1}\sum_{i=1}^{\ell}\bigg|\sum_{j=1}^m s_{ij}\bigg|\;.
$$
On the other hand, if $B_{\ell}-A_{\ell}\ge 0$, then
$$
B_{\ell}=(B_{\ell}-A_{\ell})+A_{\ell}
= \big|B_{\ell}-A_{\ell}\big|+A_{\ell}
\le \big|B_{\ell}-A_{\ell}\big|+2^{\ell-2}
\sum_{i=1}^{\ell-1}\bigg|\sum_{j=1}^m s_{ij}\bigg|\;,
$$
so by the inductive hypothesis again,
$$
\sum_{i=1}^{\ell}\sum_{j=i+1}^ms_{ij}
= B_{\ell}+\sum_{i=1}^{\ell-1}\sum_{j=i+1}^m s_{ij}
\le \big|B_{\ell}-A_{\ell}\big|+2^{\ell-1}
\sum_{i=1}^{\ell-1}\bigg|\sum_{j=1}^m s_{ij}\bigg|
\le 2^{\ell-1}\sum_{i=1}^{\ell}\bigg|\sum_{j=1}^m s_{ij}\bigg|\;.
$$
This completes the proof of the lemma, and with it, the proof
of \rHskew\ and hence of \rHc. \quad \qed


\bigskip
\rem{\kHz}
There is a alternative proof of \rHc\ that is less direct but
which ties in well with several deeper results in the context of
geometric invariant theory. Namely, one can study the
downwards gradient flow for the function
$H^{0,0}\owns \alpha \mapsto \norm{\alpha}^2$, which is a
flow that preserves the orbits of $\Gamma$, a fact easily
checked by solving $\dot\gamma_t^{}\gamma_t^{-1}
=[m(\alpha_t),\alpha_t]$ for
$\alpha_t=\gamma_t^{}\alpha\gamma_t^{-1}$. Modulo reparameterisation,
this flow turns out to be
the same as the downwards gradient flow for $H^{0,0}
\owns \alpha \mapsto \norm{m(\alpha)}^2$. The latter covers
the downwards gradient flow for $\m P(H^{0,0})\owns [\alpha]
\mapsto \norm{m(\alpha)}\big/\norm{\alpha}^2$,
$m(\alpha)/\norm{\alpha}^2$ being the moment map for the
action of $\Gamma$ on $\m P(H^{0,0})$ (\cite{Nes}).
An unpublished
theorem of Duistermaat using the \L{}ojasiewicz inequality
presented in \cite{Le} shows that this flow
defines a (strong)
deformation retract of the set of $\Gamma$-polystable points
onto the zero set of the moment map (analogous to the
result of Neeman \cite{Nee} in the algebraic setting), and
\rHc\ follows quite easily from this. The fact that it is
possible to obtain the very precise estimate given
in \rHskew\ is perhaps a reflection of the quasi-linearity of the
Yang-Mills equations, and does not ordinarily hold in
the fully non-linear setting; cf.\ \cite{CS}.

\bigskip

\penalty-1000

\secn{Non-stability.}

\rThere\
is a version of the Hitchin-Kobayashi
correspondence for bundles in an $L^p_1$ neighbourhood of
a polystable bundle,
but it does not provide much detail in the case of
connections and/or classes that are not polystable.
Whereas non-zero elements of $H^{0,1}$ may be unstable with
respect to the action of $\Gamma$---that is, zero is in the
closures of their orbits, \rCB\ states that
there are no strictly unstable bundles near
$E_0$, so the correspondence between the two different
notions of stability is not perfect. However, it is
nevertheless true that the interrelation between
the two notions goes further than just that described
by \rThere, as will be seen in this section.
All notation from earlier sections continues to be retained.

\medskip
In general, if ${\cal E}$  is an arbitrary torsion-free
semi-stable sheaf that is not stable, there is
a non-zero subsheaf ${\cal S}
\subset {\cal E}$ with $\mu({\cal S})=\mu(\cal E)$ and with
torsion-free quotient ${\cal Q}={\cal E}/{\cal S}$ for which
$\mu({\cal Q})=\mu({\cal E})$.
Both ${\cal S}$ and ${\cal Q}$ are necessarily semi-stable,
and if ${\cal S}$ is of maximal rank, then
${\cal Q}$ is stable. Iterating this process yields a
filtration of ${\cal E}$, $0={\cal S}_0\subset {\cal S}_1
\subset {\cal S}_2\subset\cdots \subset {\cal S}_k = {\cal E}$
such that the successive quotients are all torsion-free
and stable. Any such filtration is known as a
{\it Seshadri filtration\/} or sometimes a
Jordan-H\"older filtration, and although it is not
unique,  the graded object $Gr({\cal E})=
\bigoplus_{j=1}^k({\cal S}_j/
{\cal S}_{j-1})$ is unique after passing to the
double-dual.

In the current setting of holomorphic structures $E$ near to
$E_0$,
\rCB\ states that
$E$ is semistable whilst
\rCE\
states that any
destabilising subsheaf of $E$ is a subbundle. In this case
therefore, there is a Seshadri filtration of $E$ defined
by subbundles, so the graded object $Gr(E)$ associated
to $E$ is a polystable holomorphic structure
on $\et$ close to $d_0$ in $L^p_1$.

Recall from the proof of \rCA\ that
if $A$ is a holomorphic subbundle of $E$ with quotient
$B$, then in a unitary frame for $A$ and $B$, a hermitian
connection $d_E$ on $E$ and its curvature $F_E$
have the form
$$
d_E=\bmatrix{d_A& \beta\cr -\beta^*& d_B}\,,
\qquad
F_E= \bmatrix{F_A-\beta\w\beta^* & d_{{}_{BA}}\beta\cr -d_{{}_{AB}}\beta^*
& F_B-\beta^*\!\w \beta}\;,
$$
where now $\beta\in A^{0,1}(Hom(B,A))$ is a $\db$-closed
form representing the extension $0\to A\to E\to B\to 0$
and where $d_A$ and $d_B$ are the
connections on $A$ and $B$ induced by the hermitian structure
and $d_E$. If $\rk A=a$ and $\rk B=b$, and if $t>0$,
let
$
h_t=\bmatrix{t^b&0\cr  0 & t^{-a}}
$
so $h_t$ is covariantly constant with respect to the
direct sum connection $d_{A\oplus B}$ on $A\oplus B$,  $\det h_t=1$, and
$h_t\cdot d_E$ has the same form as $d_E$ with $\beta$
replaced by $t^{a+b}\beta$. So as $t\to 0$, $h_t\cdot d_E
\to d_{A\oplus B}$. Proceeding inductively, it follows easily that
there exist $g_t\in {\cal G}$ such that $g_t\cdot (d_0+a)$
converges in $C^{\infty}$ to the Hermite-Einstein
connection on $Gr(E)$.

Note that by \rCD, every holomorphic
endomorphism of $Gr(E)$ is also $\db_0$-closed, and is
therefore covariantly constant with respect to $d_0$, by
\rCCa. Thus the automorphisms $g_t\in {\cal G}$ can even
be taken to lie in $\Gamma$. The following result gives
something of a converse to this observation, albeit in
a rather special case. In its hypotheses,  how small is
``sufficiently small" is determined by the connection $d_0$,
so that \rAD\ is applicable.

\medskip
\lem{\Approx}
Let $d_0+a$ be an integrable connection on
$\et$ with $\norm{a}_{L^p_1}$ sufficiently small, and suppose that
$a''=\alpha+\db_0^*\beta$
for some  $\alpha\in H^{0,1}$ and $\beta\in A^{0,2}(\eet)$.
Then the following are equivalent:
\smallskip

{\advance\parindent by 35pt\parskip=0pt \sl
\item{1. }  For any $\epsilon>0$ there exists
$\gamma\in\Gamma$ such that $\norm{\gamma\alpha\gamma^{-1}}<\epsilon$;
\item{2. } For any $\epsilon>0$ there exists $g\in {\cal G}$
such that $\norm{g\cdot(d_0+a)-d_0}_{L^p_1} <\epsilon$.
\par}

\pf
The implication 1.\ $\Rightarrow$ 2.\ follows immediately
from \rAD.
For the converse, assume $\epsilon>0$ is smaller than the
number specified in \rAE\ and let $g\in {\cal G}$ be an
automorphism such that
$\norm{g\cdot(d_0+a)-d_0}_{L^p_1}<\epsilon$. Applying
\rAE\ to the semi-connection $g\cdot(\db_0+a'')$ yields
a unique $\varphi\in A^{0,0}(\eet)$ orthogonal to
$\ker\db_0$ such that
$d_0+\tilde a:= \exp(\varphi)\cdot g\cdot (d_0+a)$
satisfies $\db_0^*\tilde a''=0$, with $\norm{\tilde a''}$
bounded by a fixed multiple of $\epsilon$. Applying \rCD\ to the
connection on $\rm Hom(\et,\et)$ induced by $d_0+\tilde a$
and $d_0+a$ and the section $\exp(\varphi)g$ of this bundle,
it follows that if $\epsilon$ is sufficiently small then
$\exp(\varphi)g =: \gamma$ must be $\db_0$-closed with
$a''\gamma=\gamma\tilde a''$. Then if $\tilde a'' = \tilde\alpha
+\db_0^*\tilde\beta$, orthogonality of the decompositions
gives $\gamma^{-1}\alpha\gamma = \tilde\alpha$, and
$\norm{\tilde\alpha}$ is bounded by a fixed multiple of $\epsilon$
since $\norm{\tilde a''}_{L^p_1}$ is. \quad \qed

\medskip
Consider now an integrable connection
$d_0+a$, with $\norm{a}_{L^p_1}$
assumed to be appropriately small and with
$a''=\alpha+\db_0^*\beta$ for some $\alpha\in H^{0,1}$
and some $\beta\in A^{0,2}(\eet)$ orthogonal to the
kernel of $\db_0^*$.
Under the action of $\Gamma$ on $H^{0,1}$,
there is a
point $\bar\alpha\in H^{0,1}$ of
smallest norm in the
closure of the orbit of $\alpha$ unique up to conjugation by
unitary elements in $\Gamma$, and this is a
$\Gamma$-polystable point (if not zero). Since $\bar\alpha$ is in the
closure of the orbit of $\alpha$ and each element
near $0$
in this orbit lies in the analytic set
$\Psi^{-1}(0)$, there is a unique
section $\bar\beta\in A^{0,2}(\et)$ such that
$\db_0+\bar a'' := \db_0+\bar\alpha+\db_0^*\bar\beta$ is integrable, so
by \rThere\ the corresponding holomorphic bundle $\bar E$
is polystable. The following is the main result of
this section:

\medskip
\thm{\nonpoly} With the preceding definitions, let
$E$ be the holomorphic structure defined by $d_0+a$.
Then $\bar E \simeq Gr(E)$.

\smallskip
The proof, which is principally by induction on the rank $r$ of $\et$
(with the initial case $r=1$ being self-evident)
proceeds in several stages, corresponding to three
cases: 1. {\sl That $\alpha=0$\/}; ~ 2.
{\sl That $\alpha$ is not zero and is not $\Gamma$-semistable\/};
and ~ 3. {\sl That $\alpha$ is $\Gamma$-semistable\/}.
The first
is the totally degenerate case for which $\alpha=0$:

\medskip
\prop{\zero} Let $d_0+a$ be an integrable connection on
$\et$ with
$\db_0^*a''=0$, and let $E$ be the corresponding
holomorphic structure. If $\norm{a}_{L^p_1}$ is sufficiently small
then $\Pi^{0,1}a''=0$ if and only
if $E\simeq E_0$.

\pf If $\Pi^{0,1}a''=0$, then it follows from \rAD\ that
$a''=0$, and therefore $a=0$. Conversely, if
$E\simeq E_0$, then there exists $g\in {\cal G}$ such
that $g\cdot d_0=d_0+a$, or  equivalently, $\db_0g+a''g=0$.
Applying \rCD\ to the connection on $Hom(\et,\et)$ induced
by $d_0$ and $d_0+a$, it follows that $d_0g=0=a''g$, so
$a''=0$. \quad \qed

\bigskip

Let $d_0+a$ be as above with $a''=\alpha+\db_0^*\beta$.
Choose a sequence $(\gamma_j)$ in $\Gamma$
with $\det\gamma_j=1$ for every $j$ such that
$\norm{\gamma_j^{}\alpha\gamma_j^{-1}}^2 \to
\inf\limits_{\gamma\in\Gamma} \norm{\gamma\alpha\gamma^{-1}}^2$
as $j\to\infty$, so
after passing to a subsequence if necessary, it can be
assumed that $\alpha_j := \gamma_j^{}\alpha\gamma_j^{-1}$
converges to $\bar\alpha\in H^{0,1}$.

If $\beta_j := \gamma_j^{}\,\beta\,\gamma_j^{-1}$ and
$a_j'' := \alpha_j+\db_0^*\beta_j$, then $d_0+a_j
=\gamma_j\cdot(d_0+a)$ is an integrable connection defining
a holomorphic structure isomorphic to $E$, with
$\db_0^*a_j''=0$. By \rAD, $\norm{a''_j}_{L^p_1}$ is
uniformly bounded independent of $j$, so after passing
to another subsequence if necessary, the connections
$d_0+a_j$ can be assumed to converge weakly in $L^p_1$
and strongly in $C^0$ (say) to a limiting connection
$d_0+\bar a$, with $a\in L^p_1$. Elliptic regularity
combined with integrability of the
connection together with the equation $\db_0^*\bar a''=0$
imply that $\bar a$ is in fact smooth. Indeed, using
the analysis of \S1, the forms $\beta_j$ can be
assumed to be converging in $L^p_2$ to a limit
in $A^{0,2}(\eet)$ that is orthogonal to $\ker\db_0^*$,
and by the uniqueness statement of \rAG, this limit
must be the form $\bar\beta$ mentioned earlier,
with $\bar a'' =\bar\alpha+\db_0^*\bar\beta$.

Since $\det\gamma_j=1$ for every $j$, it follows that
if $\norm{\gamma_j}$ is uniformly bounded then a subsequence
can be found converging to some $\gamma_0\in \Gamma$,
and then $d_0+\bar a=\gamma_0\cdot (d_0+a)$. This is the
case considered in the previous section when $\alpha\in H^{0,1}$ is
a $\Gamma$-polystable point. So it may be supposed that
$\norm{\gamma_j}$ is not uniformly bounded, and after
rescaling $\gamma_j$ to $\eta_j := \gamma_j/\norm{\gamma_j}$,
these may be assumed to converge to some $\gamma_0\in H^{0,0}$
with $\norm{\gamma_0}=1$ and $\det\gamma_0=0$. It may also
be assumed without loss of generality that $\gamma_j$ is
self-adjoint and positive for each $j$, so $\gamma_0$ is
also self-adjoint and non-negative.

The equation $\gamma_j\cdot (d_0+a)=d_0+a_j$ is equivalent
to $\db_j\gamma_j=0$ where $d_j$ is the connection on
$Hom(\et,\et)$ induced by $d_0+a$ and $d_0+a_j$, these
connections converging to the connection on this bundle induced
by $d_0+a$ and $d_0+\bar a$. So $\gamma_0$ defines a non-zero
holomorphic map from  $E$ to $\bar E$,
this map having determinant $0$.  Since $\gamma_0$ must be of
constant rank on $X$, its kernel $K$ is a holomorphic subbundle
of $E$, necessarily of the same slope as that of $\et$.
Thus $E$ may be expressed as an extension by holomorphic
semi-stable bundles $0\to K\to E\to Q\to 0$, where
$Q := E/K$.

\smallskip
Consider now Case 2.\ of \nonpoly, namely
when $\alpha$ is non-zero and is not $\Gamma$-semistable. By definition,
zero is in the closure of the orbit of $\alpha$ under $\Gamma$,
so $\bar\alpha=0$ and therefore $\bar E=E_0$ by \zero.
For notational convenience, set $E_1 := K = \ker\gamma_0$
and $E_2 := Q = E/K$. Since $\gamma_0$ is self-adjoint,
$E_2$ can be identified with $E_1^{\perp} \subset \et$ as a
unitary bundle. The holomorphic structures on $E_1$ and $E_2$
are those induced from $E$ as holomorphic sub- and quotient bundles.
But since $E_1=\ker\gamma_0$ and $\gamma_0$ is a $d_0$-closed
self-adjoint endomorphism of the holomorphic bundle $E_0$, both $E_1$ and $E_2$
have hermitian connections induced from $d_0$, with respect
to which the connections are Hermite-Einstein with the same
Einstein constant as $E_0$. These connections will be
denoted by $d_{0,1}$, $d_{0,2}$ respectively, so
$d_0=d_{0,1}\oplus d_{0,2}$ using self-evident notation.

The limit $\gamma_0$ of the (rescaled) automorphisms
$\gamma_j$ is $d_0$-closed and satisfies
$\gamma_0\alpha=0=\gamma_0\beta$ and also
$\gamma_0a''=0$, so in terms of the splitting
$\et=E_1\oplus E_2$,
$$
\gamma_0=\bmatrix{0&0\cr0&\bar\gamma}\,,\quad
\alpha=\bmatrix{\alpha_{11}&\alpha_{12}\cr 0 & 0}\,,\quad
\beta=\bmatrix{\beta_{11}&\beta_{12}\cr 0 & 0}\,, \quad
a''=\bmatrix{a_{11}'' & a_{12}''\cr 0 &0}\,, \myeqn
$$\eqtag{\conforms}
where $\bar\gamma=\bar\gamma^*$ has non-zero determinant.
Since $\gamma_0$ is $d_0$-closed, the connection on $E_2$ induced
by the connection $d_0+a$ (i.e., as a quotient of $E$) is the
same as the connection on this bundle induced by $d_0$
(i.e., as a subbundle of $E_0$), so the holomorphic bundle
$E_2$ is isomorphic to a direct sum of stable summands of $E_0$.

The connection on $E_1$ induced by $d_0+a$ is identified with
$d_{0,1}+a_{11}$, with $a_{11}=\alpha_{11}+\db_{0,1}^*\beta_{11}$.
By the inductive hypothesis (of \nonpoly), under the action of
$\Gamma_1 = Aut(E_1(d_{0,1}))$, there are connections
in the orbit of $d_{0,1}+a_{11}$ that are arbitrarily close in
$L^p_1$ to the Hermite-Einstein connection $d_G$ on $Gr(E_1)$.
Then using automorphisms of the form $h_t$ as described earlier,
the off-diagonal term $a''_{12}$ in \conforms\
can be made arbitrarily small whilst leaving the diagonal terms fixed,
from which it follows
that there exist $g_t\in {\cal G}$ for $t>0$ such that
$g_t\cdot(d_0+a)\to (d_{G}\oplus d_{0,2})$ as $t\to 0$.

Thus the two Hermite-Einstein connections $d_0$ and $d_{G\oplus E_2}$
lie in the closure of the orbit of $d_0+a$ under ${\cal G}$, but
up to unitary isomorphism, there is at most one such connection
since the space of Yang-Mills connections modulo unitary gauge is
Hausdorff, by the reasoning of \S6 of \cite{AHS}. So $Gr(E_1)\oplus
E_2\simeq E_0$, which implies that $Gr(E) \simeq E_0$.
(Note that this argument has proved that $\alpha_{11}\in H^{0,1}(E_{0,1})$
is not $\Gamma_1$-semi-stable, which is not self-evident.)

\medskip
It remains to complete the proof of \nonpoly\ in Case 3., this
being when
$\alpha\in H^{0,1}$ is $\Gamma$-semistable but not $\Gamma$-polystable.
With the same objects as earlier, this is the case that $\bar\alpha\not=0$,
so $\bar E$ is not isomorphic to $E_0$ (by \zero), but $\bar E$ is
polystable, by \rThere.

By \rCK, every endomorphism of $\et$ that is holomorphic with respect
to $\db_0+\bar a''$ is in fact covariantly constant with respect to
$d_0$ and commutes with $\bar a$, so $\bar\Gamma := Aut(\bar E)$ 
is the subgroup of
$\Gamma=Aut(E_0)$ commuting with $\bar a$. 
The connection $d_0+\bar a$ defines
an $\omega$-polystable point, and
$\gamma_j\cdot (d_0+a)\to d_0+\bar a$.
From \Approx, once $d_0+a$ has been placed in the good
complex gauge $d_0+\tilde a$
of \rAE\ with respect to the Hermite-Einstein
connection $\bar d$ inducing $\bar E$, there are
automorphisms $\bar\gamma_j$ that are $\bar d$-closed
such that $\bar\gamma_j\cdot(d_0+\tilde a)\to \bar d$. But
now Case 2.\ applies with $\bar E$ replacing $E_0$, for which 
the conclusion is that
the bundle $Gr(E)$ is isomorphic to $\bar E$,
as desired. Consequently, the proof of \nonpoly\
is complete. \quad \qed

\bigskip
\secn{Conclusion.}

\quad We end this paper with several concluding remarks.

\medskip
\proclaim{1.} The results presented here appear to be of some significance
even in the case of compact Riemann surfaces. When the
degree and the rank of $\et$ are coprime, the moduli space
of stable holomorphic structures on $\et$ is a smooth
compact manifold, of considerable interest in its own right,
not least because this space carries a natural hyper-K\"ahler metric
of Weil-Petersson type.
When the rank and degree of $\et$ are not coprime,
the stable bundles are naturally compactified by adding the
polystable bundles, and the results here provide a
description of neighbourhoods of boundary points.
The degeneration  of the metric at the boundary
of the stable moduli space can be interpreted in terms
of quotient singularities arising from the K\"ahler quotient of
a smooth metric on a framed moduli space.

The case of compact Riemann surfaces is also helpful
for obtaining a better understanding of several of the results
presented here. There are no integrability conditions to
be considered, and the only singularities occurring in moduli
spaces result from quotient singularities which can be
viewed in the light of isotropy for the action of $\Gamma$ on classes
in $H^{0,1}$. The cases of genus $0$, $1$ and $2$ for
$E_0$ being the trivial bundle of rank $2$, or the direct sum
of a non-trivial line bundle of degree $0$ with the trivial
line bundle all provide considerable insight.

\medskip

\proclaim{2.} In the case $n=2$, relatively explicit examples of
moduli spaces
can be computed, particularly when $X$ is a ruled surface and
even more particularly when $X=\m P_1\times\m P_1$. Using
monads, moduli
spaces of $2$-bundles have been computed explicitly in
\cite{Bu1}, including an explicit description of the
space of deformations of the bundle ${\cal O}(1,-1)\oplus
{\cal O}(-1,1)$, which again illustrates many of the
results here; cf.\ the \kCa\ following \rCB.
Monads also feature in Donaldson's paper \cite{Do2}, which
presents another illustration of the interrelation between
the notions of stability in K\"ahler geometry and geometric
invariant theory, one that is not independent of the results
in this paper.

\medskip
\proclaim{3. } Because the results here focus on a {\it neighbourhood\/}
of a fixed polystable bundle, it is reasonable to expect that
they will hold mutatis mutandis on arbitrary compact complex
manifolds equipped with Gauduchon metrics. However, in light
of the \kCa\ following the proof of \rCB, there may be some
unforeseen subtleties. For the sake of brevity and
simplicity, we have considered only the K\"ahler case.

\medskip
\proclaim{4.} It is evident from the analysis that the assumption
of integrability for connections is not nearly as important
as might be expected, as the $(0,2)$ and $(2,0)$ components
of the curvature are well-controlled by \rAG, given that the
calculations are local to $d_0$. This highlights the interesting
class of solutions $d$ of the Yang-Mills equations on a compact
K\"ahler manifold for which $\d F^{0,2}(d)=0=\db\wh F(d)$ (which
includes some
{\it self\/}-dual solutions on compact surfaces), these bearing
some formal similarities to solutions of the Seiberg-Witten
equations.

\medskip
\proclaim{5.} At its heart, the proof of \rCB\ is a manifestation of
a very coarse compactness property of stable bundles, a
 desirable property used to great effect in
gauge theory. Moduli spaces of stable holomorphic bundles
on a K\"ahler surface can fail to be compact in two ways,
one reflecting the degeneration from stable to polystable
and the other in terms of the concentration of curvature
of Hermite-Einstein connections. The former is the
subject of this paper, whereas the latter is considered in
\cite{Bu3}.
Although the failure of moduli spaces of stable
bundles on compact K\"ahler surfaces can be controlled
to some extent as described in that reference, in higher dimensions there is
less control on the degeneration and one is forced to
consider compactifications in terms of sheaves (\cite{BS}). 
The Bogomolov inequality $(c_2-(r-1)c_1^2/2r)
\cdot \omega^{n-1} \ge 0$ for semi-stable sheaves and
bundles does not provide sufficient control on subbundles in
dimensions greater than $2$.

\medskip

\proclaim{6.} As alluded to in the introduction,
there are profound relationships between
the theory of stable holomorphic vector bundles on
compact K\"ahler manifolds and the theory of
constant scalar curvature K\"ahler metrics, these
relationships mediated by geometric invariant theory.
In that the former theory is a quasi-linear analogue
of the latter (in the sense of partial differential
equations),
it can be hoped that the results
here may provide useful directions for the further
investigation of moduli of compact complex manifolds and
their geometries.

\medskip

\proclaim{7.} To conclude on an even more speculative note,
in view of the critical importance of Yang-Mills theory and
of representation theory in contemporary physics, it might also
be hoped that the our results may provide deeper insight
into the nature of elementary particles and their interactions.


\bigskip
\bigskip
\centerline{\bf References}

\frenchspacing
\parindent=35pt
\parskip=2pt

\bigskip


\paper{AB}{M. F. Atiyah and R. Bott}{The Yang-Mills equations 
over Riemann surfaces}{Philos. Trans. Roy. Soc. London Ser. A}%
{308}{1983}{523--615}

\paper{AHS}{M. F. Atiyah, N. J. Hitchin, I. M. Singer}%
{Self-duality in four-dimensional Riemannian geometry}%
{Proc. Roy. Soc. London Ser. A}{362}{1978}{425--461}

\chap{BS}{S. Bando and Y.-T. Siu}{Stable sheaves and 
Einstein-Hermitian metrics}{Geometry and analysis on complex manifolds}
{World Sci. Publ.}{River Edge, NJ}{1994}{39--50}

\paper{Bu1}{N. P. Buchdahl}{Stable 2-bundles on Hirzebruch surfaces} 
{Math. Z.}{194}{1987}{143--152} 

\paper{Bu2}{N. P. Buchdahl}{Hermitian-Einstein connections
and stable vector bundles over compact complex surfaces}%
{Math. Ann.}{280}{1988}{625--648} 

\paper{Bu3}{N. P. Buchdahl}{Sequences of stable bundles over
compact complex surfaces}{J. Geom. Anal.}{9}{1999}{391--428} 

\paper{CS}{X. Chen and S. Sun}{Calabi flow, geodesic rays, and uniqueness 
of constant scalar curvature K\"ahler metrics}{Ann. of Math.}{180}{2014}%
{407--454}

\paper{Do1}{S. K. Donaldson}{A new proof of a theorem of 
Narasimhan and Seshadri}{J. Differential Geom.}{18}{1983}{269--277}

\paper{Do2}{S. K. Donaldson}{Instantons and geometric invariant theory}%
{Comm. Math. Phys.}{93}{1984}{453--460}
 
\paper{Do3}{S. K. Donaldson}{Anti self-dual Yang-Mills connections 
over complex algebraic surfaces and stable vector bundles}%
{Proc. London Math. Soc.}{50}{1985}{1--26}

\book{DK}{S. K. Donaldson and P. B. Kronheimer}%
{The geometry of four-manifolds}%
{Oxford University Press}{New York}{1990}

\paper{FK}{O. Forster and K. Knorr}{\"Uber die Deformationen von 
Vektorraumb\"undeln auf kompakten komplexen R\"aumen}{Math. Ann.}%
{209}{1974}{291--346}

\paper{FS}{A. Fujiki and G. Schumacher}%
{The moduli space of Hermite-Einstein bundles on a compact K\"hler manifold}%
{Proc. Japan Acad. Ser. A Math. Sci.}{63}{1987}{69--72} 

\paper{G}{P. Gauduchon}%
{Le th\'eor\`eme de l'excentricit\'e nulle}%
{C. R. Acad. Sci. Paris S\'er. A-B}{285}{1977}{A387--A390}

\item{[Hi]~~} N. J. Hitchin, Review of 
``{\sl La topologie d'une surface hermitienne d'Einstein\/}" 
by P. Gauduchon,
C. R. Acad. Sci. Paris S\'er. A-B 290 (1980), A509--A512.
In:  Mathematical Reviews MR571563 (81e:53052), 1981.\target{Hi}

\paper{J}{A. Jacob}{Existence of approximate Hermitian-Einstein 
structures on semi-stable bundles}{Asian J. Math.}%
{18}{2014}{859--883} 

\book{K}{A. W. Knapp}{Lie Groups Beyond an Introduction (Second
edition)}{Birkh\"auser}{Basel}{2002}

\chap{KN}{G. Kempf and L. Ness}{Lengths of vectors in 
representation spaces}{Lecture Notes in Mathematics Vol. 732}%
{Springer}{Berlin-Heidelberg-New York}{1978}{233--244}

\paper{Kim}{Hong-Jong Kim}{Moduli of Hermite-Einstein vector bundles}%
{Math. Z.}{195}{1987}{143--150}

\book{Kir}{F. C. Kirwan}{Cohomology of quotients in symplectic
and algebraic geometry}{Princeton University Press}{Princeton, NJ}{1984}

\paper{Ko1}{S. Kobayashi}{First Chern class and 
holomorphic tensor fields}{Nagoya Math. J.}{77}{1980}{5--11}

\paper{Ko2}{S. Kobayashi}{Curvature and stability of vector bundles}%
{Proc. Japan Acad. Ser. A Math. Sci.}{58}{1982}{158--162} 

\book{Ko3}{S. Kobayashi}{Differential geometry of complex vector bundles}%
{Princeton University Press}{Princeton, NJ}{1987}

\paper{Kur}{M. Kuranishi}{On the locally complete families of 
complex analytic structures}{Ann. of Math.}{75}{1962}{536--577}

\paper{Le}{E. Lerman}
{Gradient flow of the norm squared of a moment map}%
{Enseign. Math.}{51}{2005}{117--127}

\chap{LY}{J. Li and S.-T. Yau}{Hermitian-Yang-Mills connection 
on non-K\"ahler manifolds}{Mathematical aspects 
of string theory (San Diego, Calif., 1986)}
{World Sci. Publishing}{Singapore}{1987}{560--573}

\paper{L\"u}{M. L\"ubke}{Stability of Einstein-Hermitian vector bundles}%
{Manuscripta Math.}{42}{1983}{245--257}

\paper{MW}{J. Marsden and A. Weinstein}%
{Reduction of symplectic manifolds with symmetry}%
{Rep. Mathematical Phys.}{5}{1974}{121--130} 

\paper{Miy}{K. Miyajima}{Kuranishi family of vector bundles 
and algebraic description of the moduli space 
of Einstein-Hermitian connections}{Publ. Res. Inst. Math. Sci.}{25}%
{1989}{301--320}

\book{MFK}{D. Mumford, J. Fogarty, F. Kirwan}
{Geometric invariant theory (3rd ed.)}%
{Springer}{Berlin}{1994}

\paper{NS}{M. S. Narasimhan and C. Seshadri}%
{Stable and unitary vector bundles on a compact Riemann surface}%
{Ann. of Math.}{82}{1965}{540--567}

\paper{Nee}{A. Neeman}{The topology of quotient varieties}
{Ann. of Math.}{122}{1985}{419--459}

\paper{Nes}{L. Ness}{A stratification of the null cone via the moment map.
 With an appendix by David Mumford.}{Amer. J. Math.}{106}{1984}{1281--1329} 
 
\paper{S}{G. Sz\'ekelyhidi}{The K\"ahler-Ricci flow and K-polystability}%
{Am. J. Math.}{132}{2010}{1077--1090}

\paper{Ta}{F. Takemoto}{Stable vector bundles on algebraic surfaces}%
{Nagoya Math. J.}{47}{1972}{29--48}

\chap{Th}{R. P. Thomas}{Notes on GIT and symplectic reduction 
for bundles and varieties}{Surveys in differential geometry Vol. 10} 
{Int. Press}{Somerville, MA}{2006}{221--273}
 
\paper{UY}{K. K. Uhlenbeck and S.-T. Yau}%
{On the existence of Hermitian-Yang-Mills connections 
in stable vector bundles}{Comm. Pure Appl. Math.}{39}{1986}%
{S257--S293}


\bigskip
\bigskip

{

\obeylines \parskip=0pt\leftskip=.5in
School of Mathematical Sciences
University of Adelaide
Adelaide,  Australia 5005
\hskip-.25in{\it E-mail address}: {\tt nicholas.buchdahl@adelaide.edu.au}

\vskip.15in
Fachbereich Mathematik und Informatik
Philipps-Universit\"at Marburg
Lahnberge, Hans-Meerwein-Strasse
D-35032 Marburg, Germany
\hskip-.25in{\it E-mail address}: {\tt schumac@mathematik.uni-marburg.de}
}

\bye